\documentclass[11pt, a4paper, oneside]{article}
\usepackage[utf8]{inputenc}
\usepackage[austrian,british]{babel}
\usepackage[english]{varioref}
\usepackage[T1]{fontenc}
\usepackage[square, numbers]{natbib}
\usepackage{amsmath}			
\usepackage{amsthm}
\usepackage{amssymb}
\usepackage[hyperfootnotes=false]{hyperref}
\expandafter\def\expandafter\UrlBreaks\expandafter{\UrlBreaks
  \do\a\do\b\do\c\do\d\do\e\do\f\do\g\do\h\do\i\do\j%
  \do\k\do\l\do\m\do\n\do\o\do\p\do\q\do\r\do\s\do\t%
  \do\u\do\v\do\w\do\x\do\y\do\z\do\A\do\B\do\C\do\D%
  \do\E\do\F\do\G\do\H\do\I\do\J\do\K\do\L\do\M\do\N%
  \do\O\do\P\do\Q\do\R\do\S\do\T\do\U\do\V\do\W\do\X%
  \do\Y\do\Z}

\usepackage{multicol}
\usepackage{tabularx}
\usepackage{siunitx}
\usepackage{fp}
\usepackage{tikz}
\usetikzlibrary{calc, arrows, decorations.markings, decorations.pathmorphing, backgrounds, positioning, fit}
\usepackage{graphicx}			
\usepackage{verbatim} 		
\usepackage{dcolumn}			
\newcolumntype{d}[1]{D{.}{.}{#1}}
\usepackage{comment}
\includecomment{comment:figures}
\usepackage{xfrac}
\usepackage{tabto}
\usepackage{setspace}

\usepackage{diagbox}			
\usepackage{adjustbox}		
\usepackage{subfigure} 		
\usepackage{multirow}	
\usepackage[en-US]{datetime2}  

\newcommand*{\citePrefix}{}

\newcommand*{\DFRLP}{DFRLP}
\newcommand*{\MCFRLP}{MC\texorpdfstring{\underline{\;\;}}DFRLP}
\newcommand*{\LCFRLP}{LC\texorpdfstring{\underline{\;\;}}DFRLP}
\newcommand*{\MCLCFRLP}{MC+LC\texorpdfstring{\underline{\;\;}}DFRLP}
\newcommand*{\CFRLP}{C\texorpdfstring{\underline{\;\;}}DFRLP}
\newcommand*{\CMCFRLP}{C+MC\texorpdfstring{\underline{\;\;}}DFRLP}

\newcommand*{\Ie}{I.\,e.\ }
\newcommand*{\ie}{i.\,e.\ }
\newcommand*{\Eg}{E.\,g.\ }
\newcommand*{\eg}{e.\,g.\ }
\newcommand*{\etal}{et al.\ }
\newcommand*{\wlogen}{w.\,l.\,o.\,g.\ }
\newcommand*{\Wlogen}{W.\,l.\,o.\,g.\ }

\newcommand*{\rz}{{\mathbb{R}}}
\newcommand*{\nz}{{\mathbb{N}}}
\newcommand*{\zz}{{\mathbb{Z}}}

\newcommand*{\defeq}{\mathrel{\vcenter{\baselineskip0.5ex \lineskiplimit0pt
			\hbox{\scriptsize.}\hbox{\scriptsize.}}}%
	=}
\newcommand*{\eqdef}{=\mathrel{\vcenter{\baselineskip0.5ex \lineskiplimit0pt
			\hbox{\scriptsize.}\hbox{\scriptsize.}}}%
}

\newcommand*{\en}{\selectlanguage{british}}
\newcommand*{\de}{\selectlanguage{austrian}}

\newcommand*{\tabbingItem}[2]{\vspace*{-4ex}\item[] \begin{tabularx}{0.935\textwidth}[t]{p{3.168cm}X} #1 & #2\end{tabularx}}
\newcommand*{\tabbingItemBeforeParagraph}{\vspace*{-1ex}}
\newcommand*{\tabbingItemAfterItemize}{\vspace*{2.5ex}}

\makeatletter
\newenvironment{figurehere}
{\def\@captype{figure}}
{}
\makeatother

\newtheorem{axiom}{Axiom}

\newtheorem{definition}[axiom]{Definition}

\newtheorem{example}{Example}

\newcommand*{\generalXScale}{9.9}
\newcommand*{\generalYScale}{3}
\newcommand*{\axisAdditionalLengthPlus}{0.5}
\newcommand*{\axisAdditionalLengthMinus}{0.15}
\newcommand*{\axisLabel}{0.1}
\newcommand*{\crossSize}{0.075}
\newcommand*{\instanceScale}{13.61}

\tikzset{
    arrow/.style={decoration={markings, mark=at position 1 with
    {\arrow[scale=1.5,>=stealth]{>}}}, postaction={decorate}},
    arrow/.default=1
}

\begin{document}
	\sloppy
	
	\en{}
	
	\title{\bf Advanced optimization models for the location of charging stations in e-mobility}
	\author{
		\sc Anna Elisabeth Kastner{\footnotemark[1]}
		\and
		\sc Peter Greistorfer{\footnotemark[3]}
		\and
		\sc Rostislav Stan\v{e}k{\footnotemark[2]}  
		}
	\date{}
	\maketitle
	\renewcommand{\thefootnote}{\fnsymbol{footnote}}
	\footnotetext[1]{
		{\tt kastner.ae@gmail.com}.
		Department of Operations and Information Systems, University of Graz, Universit\"{a}tsstra{\ss}e 15/E3, 8010 Graz, Austria}
	\footnotetext[3]{
		{\tt peter.greistorfer@uni-graz.at}.
		Department of Operations and Information Systems, University of Graz, Universit\"{a}tsstra{\ss}e 15/E3, 8010 Graz, Austria}
	\footnotetext[2]{
		{\tt rostislav.stanek@unileoben.ac.at}.
		Department of Mathematics and Information Technology, Montanuniversit\"{a}t Leoben, Peter-Tunner-Stra{\ss}e 25/I, 8700 Leoben, Austria}
	\renewcommand{\thefootnote}{\arabic{footnote}}
	
	
	
	
	\begin{abstract}

For a reduction in environmental pollution and dependency on petroleum, electric vehicles (EV) present an advantageous alternative to traditionally fossil-fuel powered automobiles. Rapid growth in the number of EVs requires an urgent need to develop an adequate charging station infrastructure to stimulate and facilitate their usage. Due to restricted investments in the development of a sufficient infrastructure, locations have to be chosen deliberately.

In this paper, three extensions considering different objectives and various constraints to the deterministic flow refuelling location problem (DFRLP), described 2017 by \citeauthor{DeVriesDuijzer:IncorporatingDrivingRangeVariabilityInNetworkDesignForRefuelingFacilities}, are introduced. In the first extension we ask how many charging stations (CS) are necessary to cover a pre-specified number of EVs and therefore exchange the original objective function for a minimizing cost function. Secondly, our research shows that, when considering location-dependent construction costs, results heavily depend on the relations of said cost differences. Tests for different cost scenarios are carried out and policy implications are discussed. In the last extension, we consider the capacity of a CS to be limited. The DFRLP assumes an unlimited capacity, meaning it is always possible to refuel all EVs at all CSs, where they stop. In our model the capacity is put into relation to the total sum of demands generated by all EVs, passing a particular CS, which means that our model determines the placement and the sizes of all CSs simultaneously.

Finally, all extensions are evaluated using benchmarks instances based on test instances from the literature.

	\end{abstract}
	
	
	\medskip

	\section{Introduction}
		\label{section:Introduction}
In the 1970s, developed countries began undertaking research and development in the electric vehicle (EV) industry responding to the oil crisis and environmental pollution. As political pressure to reduce environmental pollution and dependency on petroleum became factors of great significance, the interest of searching for alternatives to traditionally fossil-fuel powered automobiles started to increase rapidly.\cite[\citePrefix][p.\ 129]{WuZhang:CanTheDevelopmentOfElectricVehiclesReduceTheEmissionOfAirPollutantsAndGreenhouseGasesInDevelopingCountries}

The earth experiences changes in climate since the beginning of time. However, in the last centuries especially anthropogenic factors increased the CO\textsubscript{2} level and has led to accumulated greenhouse gases in the atmosphere. These effects are mostly caused by combustion of fossil fuels. Considering \eg the greenhouse gas emissions of the transportation sector in the European Union, road transport accounts for around 71\% of them.\cite{EuropeanEnvironmentAgency:GreenhouseGasEmissionsFromTransportInEurope} Comparing internal combustion engine vehicles to EVs and considering their impact on CO\textsubscript{2} emissions, EVs can have a significant effect on reducing greenhouse gases even if it is well-known that EVs are no zero-emission vehicles, \eg air pollution caused by electricity production varies substantially in different countries.\cite[\citePrefix][p.\ 129]{WuZhang:CanTheDevelopmentOfElectricVehiclesReduceTheEmissionOfAirPollutantsAndGreenhouseGasesInDevelopingCountries}


\medskip

The increased use of alternative fuel vehicles, \eg battery electric vehicles, demands more focus regarding the impact of these alternative technologies on people's driving behaviour and patterns compared to traditional combustion engine vehicle. Most important challenges arising from the widespread use of EVs are posed by vehicles' recharging needs.\cite[\citePrefix][p.\ 1]{MeharSenouci:AnOptimizationLocationSchemeForElectricChargingStations} It is a classic chicken and egg problem. Especially in the first period after introducing EVs, the investments in charging station facilities were quite scarce. There was little opportunity to make money of them because the number of EV users was relatively small. Which in turn created little incentives to switch from a car running with a combustion engine to an EV.\cite[\citePrefix][p.\ 102]{DeVriesDuijzer:IncorporatingDrivingRangeVariabilityInNetworkDesignForRefuelingFacilities}

The inadequate development of infrastructure for alternative fuels is a major difficulty for consumer acceptance and the purchase of EVs. Still, the network for alternative fuel charging stations is insufficient compared to the infrastructure needed to enable market acceptance of these vehicles. In more detail, it is not because of an insufficient development of the electricity grid in most countries, it lies in the development of public charging stations.
Thus, a lot of effort was and is needed to install a sufficient infrastructure of charging stations, since it has a major influence on whether e-mobility will prevail in the long-term. Therefore, models are designed that depict reality as good as possible in order to obtain optimal locations for installing charging stations.


\medskip

The remainder of this paper is organized as follows: Section~\ref{subsection:RelatedLiterature} gives an overview about related literature. Our contribution is specified in Section~\ref{subsection:OurContribution}.
In Section~\ref{subsection:FormalProblemDefinitionAndTheBasicModel} the deterministic flow refuelling location problem (DFRLP) is formally introduced and three extensions to the DFRLP are described in Section~\ref{section:EnhancedModels}. Section~\ref{section:Evaluation} presents and analysis the numerical results. The final Section~\ref{section:ConclusionsAndFutureResearch} offers conclusions and states opportunities for future research.

		\subsection{Related literature}
			\label{subsection:RelatedLiterature}
Starting with a paper from 1990, written by \citet{Hodgson:AFlowCapturingLocationAllocationModel}, several researchers have been working on improving the allocation-location model for “flow-capturing” and later on “flow refuelling” problems.\cite[\citePrefix][p.\ 125]{KubyLim:TheFlowRefuelingLocationProblemForAlternativeFuelVehicles} In these problems a flow depicts a certain amount of EVs travelling cyclically from the same origin to the same destination point. These models aim for a deliberate choice of locations for charging stations, since budget for the construction of charging stations is limited in most cases.\cite[\citePrefix][p.\ 102]{DeVriesDuijzer:IncorporatingDrivingRangeVariabilityInNetworkDesignForRefuelingFacilities}
The difference between flow capturing location problems (FCLP) and flow refuelling location problems (FRLP) is that the first one assumes a flow to be covered if a flow makes use of a single facility, located somewhere along its path. This assumption is no longer adequate as the FCLP does not take the driving range of an EV into consideration and up to the present state of art the driving range of EVs is still a limiting factor. It might be necessary to stop more than once along a path to refuel the EV in order to complete a tour without running out of fuel. This is taken into consideration in the FRLP.\cite[\citePrefix][p.\ 125]{KubyLim:TheFlowRefuelingLocationProblemForAlternativeFuelVehicles}

Allocation-location models are normally designed to optimally choose among a set of potential locations and try to solve the question of where to open a set of facilities and allocate demand to these facilities.
In general, there is a distinction between two different ways of expressing demand in a network: 
\begin{itemize}
    \item Node-based: demand is expressed at fixed points in the network and facilities are located centrally to fulfil demand. Well-known problems are the $p$-median problem (see \cite{Hakimi:OptimumLocationsOfSwitchingCentersAndTheAbsoluteCentersAndMediansOfAGraph} and \cite{ReVelleSwain:CentralFacilitiesLocation}) and the location set-covering problem (see \cite{ToregasSwainReVelleBergman:TheLocationOfEmergencyServiceFacilities}). These models are based on the assumption that special-purpose trips are made between the points of demand and the locations of the facility in order to satisfy demand.
    \item Flow-based: \citeauthor{Hodgson:AFlowCapturingLocationAllocationModel} argued that the assumption of defining node-based demand is not adequate for planning a charging station infrastructure for EVs, because charging stations serve demand in form of traffic flows driving by these facilities. A traffic flow is meant to be a flow between an origin-destination pair with a certain flow volume representing the number of EVs travelling along this flow. EVs are recharged at any facility located along the shortest path from an origin point to a destination node.\cite[\citePrefix][p.\ 126]{KubyLim:TheFlowRefuelingLocationProblemForAlternativeFuelVehicles}
\end{itemize}

In order to take the driving range of an EV explicitely into account, \citet{KubyLim:TheFlowRefuelingLocationProblemForAlternativeFuelVehicles} developed a novel model formulated as a two-stage approach, which determines the maximum flow volume covered by using pre-generated combinations of facilities that do not exceed an EV’s driving range. These exogenously generated combinations serve as an input for the MILP solved in the second stage. For determining these combinations an algorithm described in \cite{KubyLim:TheFlowRefuelingLocationProblemForAlternativeFuelVehicles} can be applied. This algorithm generates all combinations of potential facility locations that can refuel a flow by considering a given driving range. This formulation is limited in its applicability, because the generation of feasible combinations for all paths is computationally burdensome, even for small instances.\cite[\citePrefix][p.\ 622]
{CaparKuby:AnEfficientFormulationOfTheFlowRefuelingLocationModelForAlternativeFuelStations}
Later on, \citet{LimKuby:HeuristicAlgorithmsForSitingAlternativeFuelStationsUsingTheFlowRefuelingLocationModel} provide some heuristic algorithms, including greedy and genetic algorithms, to overcome the two-stage approach of the time-consuming pre-generation in the first stage and solving a MILP to locate a certain number of refuelling stations in the second stage.\cite[\citePrefix][p.\ 624]
{CaparKuby:AnEfficientFormulationOfTheFlowRefuelingLocationModelForAlternativeFuelStations}

A radically new MILP formulation was later developed by  \citet{CaparKuby:AnEfficientFormulationOfTheFlowRefuelingLocationModelForAlternativeFuelStations}: it does
not require pre-generated combinations of possible facility locations as input. Their formulation is similar in logic and functionality to the original model of \citet{KubyLim:TheFlowRefuelingLocationProblemForAlternativeFuelVehicles}, but the underlying logic of the pre-generation of feasible combinations for each flow is incorporated into the constraints for the model formulation itself.\cite[\citePrefix][p.\ 626]{CaparKuby:AnEfficientFormulationOfTheFlowRefuelingLocationModelForAlternativeFuelStations}
This results in a much larger and more complex model, but the solving time for the FRLP remains stable or decreases compared to the heuristic algorithms published by \citet{LimKuby:HeuristicAlgorithmsForSitingAlternativeFuelStationsUsingTheFlowRefuelingLocationModel}. Nevertheless the driving range is still not considered explicitly in the model formulation.\cite[\citePrefix][p.\ 622]{CaparKuby:AnEfficientFormulationOfTheFlowRefuelingLocationModelForAlternativeFuelStations}

The basic FRLP assumes that an infinite number of EVs can recharge at a station. This assumption becomes impractical as the number of EVs increases. \citet{UpchurchKubyLim:AModelForLocationOfCapacitatedAlternativeFuelStations} were the first researchers addressing this concern in their model. They extended the original FRLP by \citet{KubyLim:TheFlowRefuelingLocationProblemForAlternativeFuelVehicles}, which requires pre-generated combinations of facilities as input, by limiting the capacity of EVs rechargeable at a refuelling station. \citet{HosseiniMirhassani:AHeuristicAlgorithmForOptimalLocationOfFlowRefuelingCapacitatedStation} defined capacity more accurately by estimating the quantity of consumed fuel at each refuelling station. For large-scale networks these authors proposed an effective heuristic algorithm to obtain good quality solutions in appropriate time. \citet{WangLin:LocatingMultipleTypesOfRechargingStationsForBatteryPoweredElectricVehicleTransport} consider a limiting budget constraint and multiple types of charging stations in their capacitated FRLP.

Finally, \citet{DeVriesDuijzer:IncorporatingDrivingRangeVariabilityInNetworkDesignForRefuelingFacilities} worked out a MILP formulation, which takes the driving range of an EV explicitly as an input parameter into account. The authors introduce two FRLPs in their paper: one deterministic and one stochastic. Since our research is based on the deterministic one, we call it the {\em deterministic flow refuelling location problem} ({\em DFRL}) in the reminder of this paper.

		\subsection{Our contribution}
			\label{subsection:OurContribution}
			Our research is based on the deterministic flow refuelling location problem (DFRLP) formulation of  \citet{DeVriesDuijzer:IncorporatingDrivingRangeVariabilityInNetworkDesignForRefuelingFacilities} and focuses on investigating the model while considering different objectives and taking various constraints into account. 

The first model extension, discussed in Section 2.1, addresses the problem of driving range anxiety by installing a minimum number of charging stations in order to guarantee that a prespecified proportion of the total flow volume is covered. This usually (as demonstrated in Section~\ref{subsubsection:naMinCoverage}) results in covering of frequently used routes and thus
increasing the general acceptance of e-mobility.

The second model extension, introduced in Section~\ref{subsection:LocationDependentCost}, deals with location-dependent costs while maximizing the total flow volume covered (CFV). It returns useful information on the effects of cost differences concerning the construction costs for charging stations. This allows governments to plan their subsidy measures accordingly to control and support investments in installing adequate infrastructures in certain locations. 

Due to the significantly rising number of EVs, it becomes necessary to think about capacity limits at charging stations and therefore not only the location, but also the size of charging stations. The enhanced model, described in Section~\ref{subsection:StationSize}, considers both, the location and size of charging stations, and again maximizes the total flow covered. In comparison to existing models, capacity is defined as the quantity of energy available at a charging pole. Available energy at charging stations is not actually a limiting factor, but it is implicitly used as a measure for the duration of loading. 



\subsection{Formal problem definition and the basic model}
		\label{subsection:FormalProblemDefinitionAndTheBasicModel}
The objective of the DFRLP (see \cite{DeVriesDuijzer:IncorporatingDrivingRangeVariabilityInNetworkDesignForRefuelingFacilities}) is to maximize the total number of EVs which can complete their trip without running out of fuel by optimally locating an exogenously given number of charging stations. Given an undirected graph $G(L, E)$, $L$ is a set of nodes, \ie locations, and $E$ is a set of edges, \ie streets, between these locations. The set of locations $L$ is the union of three disjunct subsets: set of driving origins $O$, set of driving destinations $D$ and the set of potential facility locations $K$ for charging stations (CSs) with no capacity limitations.

The overall amount in the traffic network is defined by set of cyclic {\em flows} $F$, where cyclic means that every EV after visiting its driving destination returns back home to its origin by using the same path. Each flow $f \in F$ is sufficiently defined by its origin $O_f \in O$, its integer flow volume $v_f \in \mathbb{N}$, its destination $D_f \in D$, and by the desired path between the origin-destination pair in $G$. The nodes on this path are potential facility locations and define the set $K_f \subseteq K$. All in all, $L_f = \{O_f\} \cup K_f \cup \{D_f\}$. We say a flow is {\em covered} if the driving distance between consecutively used CSs along a round-trip does not exceed the driving range of the EV.
Like $G$ and $F$, the fixed driving range $R$ is a parametric input.

Based on the cyclic property of our flows, similarly to \cite[\citePrefix][p.\ 103f]{DeVriesDuijzer:IncorporatingDrivingRangeVariabilityInNetworkDesignForRefuelingFacilities}, we define sub-trips from one CS to another as cycle segments:
\begin{definition} \label{definition:def1}
A cycle segment of the flow $f$ is identified by two nodes $k$ and $l$ and has corresponding distances $t_{k l} \in \mathbb{R}$ as defined below: 
	\begin{itemize}
		\item If $k = O_f$ and $l \in K_f$, the cycle segment defined by these two nodes is the path $l \rightarrow O_f \rightarrow l$ and its distance $t_{k l}$ is given by the distance from $l$ via $O_f$ to $l$, both along $f$.
		
		\item For all $k, l \in K_f$, where $k$ occurs before $l$ in the flow $f$ on the way from $O_f$ to $D_f$, the cycle segment defined by these two nodes is the path $k \rightarrow l$ and its distance $t_{k l}$ is given by the distance from $k$ to $l$, both along $f$.
			
		\item If $k \in K_f$ and $l = D_f$, the cycle segment defined by these two nodes is the path $k \rightarrow D_f \rightarrow k$ and its distance $t_{k l}$ is given by the distance from $k$ via $D_f$ to $k$, both along $f$.
		
		\item For $k = O_f$ and $l = D_f$ we do not define cycle segment, but set $t_{k l} = \infty$.
	\end{itemize}
\end{definition}

\begin{example}
	Consider a flow $f$ corresponding to the desired path along the following nodes $10 \rightarrow 7 \rightarrow 12 \rightarrow 21$, \ie having the origin $O_f = 10$ and the destination $D_f = 21$. Assume that the distances of the connecting edges are $7$, $2$, and $3$, respectively. Then the six possible cycle segments have distances: $t_{10, 7} = 14$, $t_{10, 12} = 18$, $t_{7, 12} = 2$, $t_{7, 21} = 10$, and $t_{12, 21} = 6$. Further, we have $t_{10, 21} = \infty$, but do not define a cycle segment for $k = 10$ and $l = 21$.
\end{example}


It is important to mention, that no charging station is possible at an origin or destination node in the DFRLP. However, to enable this property, one must only replace the original $O_f$ and/or $D_f$ with a dummy facility location, which is connected to the original $O_f$ and/or $D_f$ by a zero-distance edge.

Our notation and the DFRLP model are based on~\citet{DeVriesDuijzer:IncorporatingDrivingRangeVariabilityInNetworkDesignForRefuelingFacilities}.

			\tabbingItemBeforeParagraph\paragraph{Parameters:}
				\begin{itemize}\tabbingItemAfterItemize
					\tabbingItem{$F$}{set of flows}
					\tabbingItem{$O$ ($O_f \in O$)}{set of origins (origin of flow $f$)}
					\tabbingItem{$K$ ($K_f \in K$)}{set of potential facility locations (along flow $f$)}
					\tabbingItem{$D$ ($D_f \in D$)}{set of destinations (destination of flow $f$)}
					\tabbingItem{$v_f \in \mathbb{N}$}{volume of flow $f$}
					\tabbingItem{$L$ ($L_f \in L$)}{set of locations, \ie $L = O \cup K \cup D$ (set of locations along flow $f$, \ie $L_f = \{O_f\} \cup K_f \cup \{D_f\}$)}\tabbingItemAfterItemize
					\tabbingItem{$E$}{set of edges between locations}
					\tabbingItem{$L_{k f}^- \subsetneq L$ ($L_{k f}^+ \subsetneq L$)}{set of locations along flow $f$ passed before (after) location $k$ on a trip from $O_f$ to $D_f$}\tabbingItemAfterItemize
			
					\tabbingItem{$p \in \mathbb{N}$}{number of new facilities to locate}
					\tabbingItem{$R \in \mathbb{N}$}{driving range}
					\tabbingItem{$t_{k l} \geq 0$}{length of the cycle segment identified by locations $k$ and $l$}
				\end{itemize}
			
			\tabbingItemBeforeParagraph\paragraph{Decision variables:}
				\begin{itemize}\tabbingItemAfterItemize
					\tabbingItem{$x_k \in \{0, 1\}$}{$1$ if a facility is placed at location $k$ and $0$ otherwise}
					\tabbingItem{$y_f \in \{0, 1\}$}{$1$ if flow $f$ is covered and $0$ otherwise}
					\tabbingItem{$i_{k l f}  \in \{0, 1\}$}{$1$ if cycle segment $k$, $l$ is used in flow $f$ and $0$ otherwise}
				\end{itemize}
			
			\paragraph{DFRLP:}
				\allowdisplaybreaks[1]
				\begin{alignat}{5}
					\label{equation:BasicModel:OV}
					\mbox{max} \ 
						& \sum_{f \in F} v_f y_f \hspace*{-3cm}
							&&
								&
									&&\\
					\label{equation:BasicModel:st:NumberOfFacilities}
					\mbox{s.t.} \ 
						&& \sum_{k \in K} x_k \ 
							&& = \ 
								& p
									&&\\
					\label{equation:QP:st:DrivingRangeConstraint}
						&& \sum_{l \in L_{k f}^+} i_{k l f} t_{k l} - (1 - y_f) M \ 
							&& \leq \ 
								& R \quad
									&& f \in F,\ k \in \{O_f \cup K_f\}\\
					\label{equation:QP:st:eq4}
						&& \sum_{l \in L_{k f}^+} i_{k l f} \ 
							&& = \ 
								& x_k \quad
									&& f \in F,\ k \in K_f\\
					\label{equation:QP:st:eq5}			
						&& \sum_{l \in L_{O_f f}^+} i_{O_f l f} \ 
							&& = \ 
								& 1 \quad
									&& f \in F\\
					\label{equation:QP:st:eq6}
						&& \sum_{k \in L_{l f}^-} i_{k l f} \ 
							&& = \ 
								& x_l \quad
									&& f \in F,\ l \in K_f\\
						\label{equation:QP:st:eq7}
						&& \sum_{k \in L_{D_f f}^-} i_{k  D_f f} \ 
							&& = \ 
								& 1 \quad
									&& f \in F\\
						\label{equation:QP:st:DecVarI}
						&& i_{k l f} \ 
							&& \in \ 
								& \{0, 1\} \quad
									&& f \in F,\ k \in \{O_f \cup K_f\},\ l \in L_{k f}^+\\
						\label{equation:QP:st:DecVarXY}
						&& x_k,\ y_f \ 
							&& \in \ 
								& \{0, 1\} \quad
									&& k \in K,\ f \in F
				\end{alignat}
				\allowdisplaybreaks[0]

The objective \eqref{equation:BasicModel:OV} of the DFRLP is to maximize the total number of EVs covered by optimally locating an exogenously given number of charging stations stated by constraint \eqref{equation:BasicModel:st:NumberOfFacilities}. 
Constraint \eqref{equation:QP:st:DrivingRangeConstraint} ensures that a flow is covered if the length of each cycle segment used along the path of flow $f$ does not exceed the driving range of an EV. Constraints \eqref{equation:QP:st:eq4}--\eqref{equation:QP:st:eq7} are flow constraints, which also link variables $x$ and $i$.
Constraints \eqref{equation:QP:st:DecVarI}--\eqref{equation:QP:st:DecVarXY} define the decision variables of the model.


	\section{Enhanced models}
		\label{section:EnhancedModels}
Most of the research focusing on optimal charging station placements for EVs consider some parameter and ignore others. Usually, only the travel distance and flow volume along paths between potential facility locations are considered. The aim of the following enhanced models is to investigate the DFRLP and to take various constraints into account, such as guaranteeing coverage for a certain proportion of EVs, considering construction costs of charging stations as well as limited capacity at charging stations in order to depict reality as good as possible.
    
		\subsection{Minimum flow volume coverage (\MCFRLP)}
		    	\label{subsection:MinCoverage}
Most of the existing research on the deployment of an infrastructure for EVs does not consider the driver convenience issue as a hard constraint when locating charging stations. The formulation of this extended model requires a minimum coverage of all EVs travelling within a network. Our ideas are based on the concept of governmental concessions awarded to construction companies, which are able to guarantee a minimum coverage level within a pre-specified time period. 
Therefore the objective of this \MCFRLP{} extension is to choose a minimum number of optimally located refuelling stations in order to guarantee that a given proportion of all EVs can complete their trip without running out of fuel. In consistence with the basic model, we assume that the capacity at all refuelling stations is unlimited. 

				\tabbingItemBeforeParagraph\paragraph{Additional parameter:}
				\begin{itemize}\tabbingItemAfterItemize
                    \tabbingItem{$C \in (0, 1]$}{minimum coverage level as proportion of the total flow volume}
                \end{itemize}
                
                \paragraph{\MCFRLP\ Model:}
                \allowdisplaybreaks[1]
                \begin{alignat}{5}
                \label{equation:MinCoverage:OV}
                    \mbox{min} \ 
                        & \sum_{k \in K} x_k \ 
                            &&
                                &
                                    &&\\
                    \mbox{s.t.} \ 
                    \label{equation:MinCoverage:st:CoverageConstraint}
                        && \frac{\sum_{f \in F} v_f y_f}{\sum_{f \in F} \  v_f}\ 
                            && \geq \ 
                                & C
                                    &&\\                            
                        && \eqref{equation:QP:st:DrivingRangeConstraint} 
                            && - \ 
                                & \hspace*{-0.15cm}\eqref{equation:QP:st:DecVarXY}. \quad
                                    &&\nonumber
                \end{alignat}
                \allowdisplaybreaks[0]
                
The objective function of the \DFRLP{} is replaced with function \eqref{equation:MinCoverage:OV} requiring a minimum number of new charging facilities to be built, while at least a certain proportion of all EVs should be able to complete their round-trip successfully;  this is ensured with constraint \eqref{equation:MinCoverage:st:CoverageConstraint}. Since the number of CSs becomes the objective, equation \eqref{equation:BasicModel:st:NumberOfFacilities} is skipped, but constraints \eqref{equation:QP:st:DrivingRangeConstraint}--\eqref{equation:QP:st:DecVarXY} from the original model are required.


	\subsection{Location-dependent costs per charging station (\LCFRLP)}
		    	\label{subsection:LocationDependentCost}
		    	While planning a network of CSs, one has to consider different one-off costs for construction depending on the CS location. This results mainly from different land costs in urban, sub-urban and rural areas. Thus, this extension, \LCFRLP, involves location-dependent construction costs as an additional input parameter in the optimization model.

                \tabbingItemBeforeParagraph\paragraph{Additional parameters:}
                \begin{itemize}\tabbingItemAfterItemize
                    \tabbingItem{$c_k \geq 0$}{construction costs per charging station at location $k$}
                    \tabbingItem{$B > 0$}{available budget for all stations}
                \end{itemize}

                \paragraph{\LCFRLP\ Model:}
                \allowdisplaybreaks[1]
                \begin{alignat}{5}
                \tag{\ref{equation:BasicModel:OV}}
                \mbox{max} \ 
                        & \sum_{f \in F} v_f y_f \hspace*{-3cm}
                            &&
                                &
                                    && \nonumber \\      
                    \label{equation:LocationDependentCost:st:BudgetConstraint}                         
                    \mbox{s.t.} \ 
                        && \sum_{k \in K} c_k x_k \ 
                            && \leq \ 
                                & B
                                    &&\\
                        && \eqref{equation:QP:st:DrivingRangeConstraint} 
                            && - \ 
                                & \hspace*{-0.15cm}\eqref{equation:QP:st:DecVarXY}. \quad
                                    &&\nonumber         
                \end{alignat}
                \allowdisplaybreaks[0]
 				
 				The objective function \eqref{equation:BasicModel:OV} and constraints \eqref{equation:QP:st:DrivingRangeConstraint}--\eqref{equation:QP:st:DecVarXY} remain the same as in the DFRLP. Again, the number of new facilities to locate is no longer exogenously given, \ie constraint \eqref{equation:BasicModel:st:NumberOfFacilities} is not used. In fact, this restriction is replaced by constraint \eqref{equation:LocationDependentCost:st:BudgetConstraint}, which takes into account that there is only a limited budget to build the CS infrastructure.

		\subsection{Capacitated DFRLPs: determination of the station size}
		    	\label{subsection:StationSize}
Most articles in the field of FRLP are based on the assumption that the capacity of a refuelling station is unlimited. Meaning, the availability of a charging station is sufficient to refuel all flows using this location, regardless of their flow volume and their distance travelled since the last refuelling process.\cite[\citePrefix][p.\ 85]{UpchurchKubyLim:AModelForLocationOfCapacitatedAlternativeFuelStations} In contrast to the early era of e-mobility, with nowadays increasing EV volume, it becomes important to cope with limited CS capacity.\cite[\citePrefix][p.\ 86]{UpchurchKubyLim:AModelForLocationOfCapacitatedAlternativeFuelStations}


The objective of this extension is to simultaneously decide upon the placement and the size of CSs in a network in order to maximize the flow volume covered. Clearly, CS size depends on the number of charging poles installed at a location. In the \CFRLP{} we assume that each charging pole has a limited capacity. We define the pole capacity as energy output in terms of total driving range per pole and period. An optimal location of charging poles is necessary in order to build more charging poles at CSs, where the energy demand is higher. A positive side effect of considering the CS size (number of charging poles per CS) as decision variable, is the avoidance of idle charging facilities due to a low utilisation at charging points.\cite[\citePrefix][p.\ 4]{WuNiu:StudyOnInfluenceFactorsOfElectricVehiclesChargingStationLocationBasedOnISMAndFMICMAC}

\medskip

To further approach reality, \citet{UpchurchKubyLim:AModelForLocationOfCapacitatedAlternativeFuelStations} take limited capacity of charging stations into account by defining capacity as the number of vehicles ``refuelable'' at a station.\cite[\citePrefix][p.\ 1379]{HosseiniMirhassani:AHeuristicAlgorithmForOptimalLocationOfFlowRefuelingCapacitatedStation} 
This leads to the following surreal model assumption that independently of the battery energy level, each EV demands a constant predefined amount of energy at each CS along its flow. Therefore, in our \CFRLP{} the demand for energy at any station depends on the current battery level.


We build on the ideas of \citet{HosseiniMirhassani:AHeuristicAlgorithmForOptimalLocationOfFlowRefuelingCapacitatedStation}, where energy consumption is assumed to be linearly proportional to the distance travelled since the last refuelling process. Moreover we define capacity as the amount of energy available and assume that batteries are always filled to full capacity (``driving range anxiety'').
\cite[\citePrefix][p.\ 1380]{HosseiniMirhassani:AHeuristicAlgorithmForOptimalLocationOfFlowRefuelingCapacitatedStation}

\medskip
In this section, we will discuss two different approaches: first we maximize the flow coverage under an exogenous number of CSs, later on, we minimize the number of installed CSs while satisfying a given coverage level. 


		    \subsubsection{Exogenous given number of charging poles (\CFRLP)}
				\label{subsubsection:Exogenous}
To begin with, we assume that the number of charging poles is exogenous. As construction costs are not explicitly concerned in this model, the exogenous number of CSs (and thus poles as well) can be seen as a proxy for the budget available for deploying an infrastructure.\cite[\citePrefix][p.\ 400]{GimenezGaydouRibeiroGutierrezAntunes:OptimalLocationOfBatteryElectricVehicleChargingStationsInUrbanAreasANewApproach}
Moreover, the assumption from \citet{HosseiniMirhassani:AHeuristicAlgorithmForOptimalLocationOfFlowRefuelingCapacitatedStation} is taken, stating that flows are assumed to be divisible, \ie flow coverage may be lower than 100\%.\cite[\citePrefix][p.\ 1382]{HosseiniMirhassani:AHeuristicAlgorithmForOptimalLocationOfFlowRefuelingCapacitatedStation}  

                \tabbingItemBeforeParagraph\paragraph{Additional decision variables:}
                \begin{itemize}\tabbingItemAfterItemize
                    \tabbingItem{$n_k \in \mathbb{N}_0$}{number of charging poles at location $k$}
                    \tabbingItem{$z_f \in [0, 1]$}{proportion of flow $f$ that is covered}
                    \tabbingItem{$w_{k l f} \in [0, 1]$}{auxiliary variable for linearisation}
                \end{itemize}      
                
                \tabbingItemBeforeParagraph\paragraph{Additional parameters:}
                \begin{itemize}\tabbingItemAfterItemize
                    \tabbingItem{$\mbox{\it Cap} \in \mathbb{N}$}{capacity of charging pole given as the amount of available energy in distance units}\tabbingItemAfterItemize
                    \tabbingItem{$e_f \in (0, 1]$}{positive range-based refuelling proportion $\leq 1$ of flow $f$ to be covered per observation period for flows with $2 \tilde{t}_{O_f D_f} \leq R$; otherwise, $e_f = 1$: \[e_f = \frac{1}{\max\left\{1,\Bigl\lfloor\frac{R}{2 \tilde{t}_{O_f D_f}}\Bigr\rfloor\right\}}, \] where $\tilde{t}_{O_f D_f} > 0$ stays for the real distance between origin and destination (note that $\tilde{t}_{O_f D_f} \neq t_{O_f D_f}$, because we
defined $t_{O_f D_f} = \infty$ in Definition~\ref{definition:def1}); as explained later, this parameter is needed to model short trips properly}\tabbingItemAfterItemize
                    \tabbingItem{$M_k \in \mathbb{N}$}{location-dependent maximum number of charging poles at a charging station (if $M_k = M$ for all $k \in K$, use $M$)}\tabbingItemAfterItemize
                    \tabbingItem{$S \in \mathbb{N}$}{total number of charging poles to locate}
                \end{itemize}
                
                \newpage
                
                \paragraph{\CFRLP\ Model:}
                \allowdisplaybreaks[1]
                \begin{alignat}{5}
                \label{equation:StationSize:Exogenous:OV}
                    \mbox{max} \ 
                        & \sum_{f \in F} v_f z_f \hspace*{10.15cm} 
                            &&
                                &
                                    &&
                \end{alignat}
                \allowdisplaybreaks[0]
                
				\vspace*{-1.2cm}

                \allowdisplaybreaks[1]
                \begin{alignat}{5}
                \label{equation:StationSize:Exogenous:st:NumberOfChargers}
                    \mbox{s.t.} \ 
                        && \sum_{k \in K} n_k \ 
                            && = \ 
                                & S
                                    &&\\
             \label{equation:StationSize:Exogenous:st:StationSizeConstraints}
                    && x_k \
                        && \leq \
                            & n_k \quad
                                    && k \in K \\
            \label{equation:StationSize:Exogenous:st:MaxStationSize}
                    && n_k \
                        && \leq \ 
                            & M_k x_k \quad
                                    && k \in K\\
						\label{equation:StationSize:Exogenous:st:FlowConstraints}
										&& z_f
												&& \leq \
														& y_f \quad
																	&& f \in F\\
\label{equation:StationSize:Exogenous:st:CapacityConstraint_linearized}
                    && \sum_{f \in F\colon k \in K_f} \left ( \sum_{l \in L_{k f}^-} t_{l k} v_f e_f w_{l k f} + \sum_{l \in L_{k f}^+}  t_{k l} v_f e_f w_{k l f} \right) \ 
                        && \leq \ 
                            & Cap \cdot n_k \quad
                                    && k \in K 
                \end{alignat}
                \allowdisplaybreaks[0]
                
				\vspace*{-1.1cm}
                
                \allowdisplaybreaks[1]
                \begin{alignat}{5}            
                \label{equation:StationSize:Exogenous:st:Linearization1}
                    &&  w_{k l f}
                        && \leq \ 
                            & i_{k l f} \quad
                                    && f \in F, k \in K_f, l \in L_{k f}^+ \\
                    && w_{k l f}
                \label{equation:StationSize:Exogenous:st:LinearizationMiddle}
                        && \leq \
                            & z_f \quad
                                    && f \in F, k \in K_f, l \in L_{k f}^+ \\
                \label{equation:StationSize:Exogenous:st:Linearization2}
                    && w_{k l f}
                        && \l\geq \ 
                            & z_f - (1 - i_{k l f}) \quad
                                    && f \in F, k \in K_f, l \in L_{k f}^+ \\
                \label{equation:StationSize:Exogenous:st:def_zf}
                    && z_f
                        && \in \ 
                            & [0, 1] \quad
                                    && f \in F \\
                \label{equation:StationSize:Exogenous:st:def_nk}
                    && n_k
                        && \in \ 
                            & \mathbb{N}_0 \quad
                                    && k \in K\\
                \label{equation:StationSize:Exogenous:st:AuxiliaryVariableDef}
                        &&  w_{k l f}
                            && \in \ 
                                & [0, 1] \quad
                                    && f \in F, k \in K_f, l \in L_{k f}^+ \\
                    && \eqref{equation:QP:st:DrivingRangeConstraint}\hspace*{-0.12cm}
                        && - \ 
                            & \hspace*{-0.15cm}\eqref{equation:QP:st:DecVarXY}. \quad
                                &&\nonumber
                \end{alignat}
                \allowdisplaybreaks[0]

The model above covers two elemental changes from the basic DFRLP: first, we make use of a new integer variable $n_k$, defining the number of charging poles at a potential facility location. Second, by postulating a limited capacity at refuelling stations, it might be impossible for one or more CSs to satisfy the total flow volume (TFV) of the flows using said stations. Thus the continuous variable $z_f$, which indicates the proportion of flow $f$ that can be covered, is introduced.

	        

\medskip

The objective function \eqref{equation:StationSize:Exogenous:OV} maximizes the total flow volume covered. Note that it is the real-valued equivalent of objective function \eqref{equation:BasicModel:OV}, obtained by substituting the binary variable $y_f$ for the proportion $z_f$. While the \DFRLP{} parametrises the total number of CSs in \eqref{equation:BasicModel:st:NumberOfFacilities}, here the total number of charging poles is ensured by \eqref{equation:StationSize:Exogenous:st:NumberOfChargers}. Constraint \eqref{equation:StationSize:Exogenous:st:StationSizeConstraints} states that there is at least one charging pole installed when a location is intended to be a charging station. The maximum number of charging poles allowed at a charging station is defined in constraint \eqref{equation:StationSize:Exogenous:st:MaxStationSize}. \eqref{equation:StationSize:Exogenous:st:FlowConstraints} ensures that the proportional coverage of flow $f$ is zero if the flow is not covered at all, i.e.\ if $y_f = 0$. The limited capacity of charging poles is considered with constraint \eqref{equation:StationSize:Exogenous:st:CapacityConstraint_linearized}. It is important to understand that \eqref{equation:StationSize:Exogenous:st:CapacityConstraint_linearized} is already a linearisation of an originally quadratic context as explained in the following.

To guarantee that capacity is not exceeded at any charging station, it is necessary to know which flows are recharging at a station. This information can be obtained from variable $i_{k l f}$, indicating at which locations ($k$ and $l$) a flow $f$ stops for refueling. The limited capacity at charging poles is expressed by the following non-linear constraint:

            \begin{equation}
            	\label{equation:StationSize:Exogenous:st:CapacityConstraint_quadratic}
                 \sum_{f \in F\colon k \in K_f} \left ( \sum_{l \in L_{k f}^-} t_{l k} i_{l k f} + \sum_{l \in L_{k f}^+}  t_{k l} i_{k l f} \right ) v_f e_f \cdot z_f   
                 \leq 
                 Cap \cdot n_k 
                 \quad  k \in K
            \end{equation}

The left-hand side of this constraint totals the amount of energy demand at location $k$. Assuming that EVs are travelling cyclically, they stop at CS $k$ twice, once on their forward and a second time on their return journey. Consequently, the energy consumption of an EV at location $k$ is linearly proportional to the distance travelled since the last charging stop on the forward trip plus the distance from the last CS used before $k$ on the backward trip (which is expressed by means of forward distance $t_{k l}$ in our symmetric case). The right-hand side of the equation defines the total capacity (in units of distance) available at the potential location site $k$.


The linearisation is done by introducing a new variable $w_{l k f} \defeq i_{l k f} z_f$ and adding linking constraints \eqref{equation:StationSize:Exogenous:st:Linearization1}--\eqref{equation:StationSize:Exogenous:st:Linearization2}. If the right-hand side is $0$, \ie no CS is opened at location $k$, on the first view, there are four possibilities how to fulfil \eqref{equation:StationSize:Exogenous:st:CapacityConstraint_quadratic}: Either $z_f$, $e_f$, $v_f$ or $i_{k l f}$ equals $0$ for all $f \in F$, $k \in L_f$ and $l \in L_{k f}^+$. $z_f = 0$ corresponds to zero flow coverage. $e_f$ cannot be $0$ following its definition. In the special case $v_f = 0$, $z_f$ might take a value greater than zero, but this has no influence on the objective function value, which maximizes the total flow volume covered (CFV), \ie the product of $v_f$ and $z_f$. If $i_{k l f} = 0$, the flow $f$ is not covered at all according to the variable definition, because no CS is used. Finally, if $z_f = 0$ or $i_{k l f} = 0$, $w_{l k f}$ is forced to zero by \eqref{equation:StationSize:Exogenous:st:Linearization1} and \eqref{equation:StationSize:Exogenous:st:LinearizationMiddle}.


\medskip


Let us discuss the parameter $e_f$ in more detail. If $e_f \in (0, 1)$, then the range of the EV allows for a longer round-trip than the flow $f$ demands. The setting $e_f = 1$ reflects two real situations. In the first case, the EV's range exactly corresponds to the length of the round-trip. In the second case, more than one charging stop is needed to finish the whole round-trip. The latter case is upper-bounded with the value $e_f = 1$, because it is not possible to charge more energy then the battery capacity allows and more than one charging stops per round-trip are necessary.
    \begin{example}
    	Imagine an EV with driving range $R = 200$. First, we assume that the round-trip to be covered has length $100$. Then two such trips are possible before recharging is needed; consequently only half of the flow volume needs to be recharged per round-trip, \ie $e_f = 0.5$. Next, assume the round-trip length to be $200$. Now exactly one round-trip can be done without recharging, \ie $e_f = 1$. Finally, assume that our round-trip has length $400$. Then obviously one recharging process per trip is not enough and more than one CS on the trip is needed since the charging volume is bounded the battery capacity, \ie $e_f = 1$. Note, obviously a minimum of $2$ CSs is needed in this case.
		
        \end{example}

		    \subsubsection{Number of charging poles as decision variable (\CMCFRLP)}
				\label{subsubsection:Endogenous}
Allowing the number of charging poles to be no longer exogenously given, the objective is to install a minimum number of charging poles while guaranteeing a certain level of coverage. This objective seems to be reasonable as an incentive to increase consumer acceptance of EVs.

                \paragraph{Additional variables:}
                \begin{itemize}\tabbingItemAfterItemize
                    \tabbingItem{$n_k$, $z_f$}{described in Section~\ref{subsubsection:Exogenous}}
                \end{itemize}      
                
                \paragraph{Additional set and parameters:}
                \begin{itemize}\tabbingItemAfterItemize
                    \tabbingItem{$C \in [0, 1]$}{minimum required flow volume coverage level}
                    \tabbingItem{$\mbox{\it Cap}$, $e_f$, $M_k$ }{described in Section~\ref{subsubsection:Exogenous}}
                \end{itemize}

                \paragraph{\CMCFRLP\ Model:}
                \allowdisplaybreaks[1]
                \begin{alignat}{5}
                \label{equation:StationSize:Endogenous:OV}
                    \mbox{min}\hspace*{-0.08cm} \ 
                        & \sum_{k \in K} n_k \ 
                            &&
                                &
                                    &&\\
%
                    \label{equation:StationSize:Endogenous:st:CoverageConstraint}
                    \mbox{s.t.} \ 
                        && \frac{\sum_{f \in F} v_f z_f}{\sum_{f \in F} v_f} \ 
                            && \geq \ 
                                & C
                                    &&\\
                        && \eqref{equation:StationSize:Exogenous:st:StationSizeConstraints} 
                        && - \ 
                            & \hspace*{-0.15cm}\eqref{equation:StationSize:Exogenous:st:AuxiliaryVariableDef} \quad
                                &&\nonumber \\
                        && \eqref{equation:QP:st:DrivingRangeConstraint} 
                            && - \ 
                                & \hspace*{-0.15cm}\eqref{equation:QP:st:DecVarXY}. \quad
                                    &&\nonumber
                \end{alignat}
                \allowdisplaybreaks[0]
                
The objective function \eqref{equation:StationSize:Endogenous:OV} minimizes the number of charging poles located within the network.
Constraint \eqref{equation:StationSize:Endogenous:st:CoverageConstraint} forces the model to install enough charging poles to cover at least a certain proportion of all EVs driving within the network. 
Furthermore, constraints \eqref{equation:StationSize:Exogenous:st:StationSizeConstraints}--\eqref{equation:StationSize:Exogenous:st:AuxiliaryVariableDef} are borrowed from the model described in Section~\ref{subsubsection:Exogenous}. Constraints \eqref{equation:QP:st:DrivingRangeConstraint}--\eqref{equation:QP:st:DecVarXY} take the driving range explicitly into account and define the decision variables.

	\section{Evaluation}
		\label{section:Evaluation}
In this section, we perform a numerical analysis of the basic model (DFRLP) and its previously described extensions.
AMPL-IDE\footnote{Version: 3.5.0.201802140038.} was used for modelling and calculation. 
The AMPL built-in parameters \verb|_ampl_elapsed_time| and \verb|_solve_elapsed_time| are used to display the solve time in seconds. Solve time is the sum of the both time parameters: the first one measures the elapsed seconds after starting the AMPL process and generating the model, preparing the sets and constraints and the second time parameter determines the solution time in seconds of the solver. 
The models were solved with Gurobi\footnote{version: Gurobi 8.1.0} and the optimization was performed on a macOS Sierra\footnote{version: macOS Sierra 10.12.6} computer with a 2.3 GHz Intel Core i5\footnote{processor type: I5-7360U} processor and an 8 GB/2133 MHz memory.

		\subsection{Benchmark instances}
			\label{subsubsection:BenchmarkInstances}	
All models were tested on four instances randomly generated and described by \citet{DeVriesDuijzer:IncorporatingDrivingRangeVariabilityInNetworkDesignForRefuelingFacilities}; the names are in the same form as in \cite{DeVriesDuijzer:IncorporatingDrivingRangeVariabilityInNetworkDesignForRefuelingFacilities}: ``sXwY'', where ``X'' stays for the number of potential facility locations $|K|$ and ``Y'' corresponds to the number of origin and destination nodes $\lvert OD \rvert$.

			Additionally to test instances taken from \cite{DeVriesDuijzer:IncorporatingDrivingRangeVariabilityInNetworkDesignForRefuelingFacilities} we use data from \citet{CaparKuby:AnEfficientFormulationOfTheFlowRefuelingLocationModelForAlternativeFuelStations} who studied and applied a different model formulation for the FRLP to Florida's highways. This real-world instance is significantly larger than the other ones and has 302 facility locations and \num{2701} flows. Unfortunately, not all of our models can yield meaningful results for such a huge instance in a reasonable amount of time and, moreover, due to missing geographic information, an appropriate classification of particular charging station locations in the \LCFRLP\ is not possible. Thus, we used this test instance only for the \MCFRLP.

\medskip

The driving range of an EV is assumed to be 250 for the basic model and all extensions; this conforms to the assumption of both \citet{DeVriesDuijzer:IncorporatingDrivingRangeVariabilityInNetworkDesignForRefuelingFacilities} and \citet{CaparKuby:AnEfficientFormulationOfTheFlowRefuelingLocationModelForAlternativeFuelStations}.

Table~\ref{table:Characteristics} shows a summary reflecting the characteristics of all test instances and Table~\ref{table:solDFRLP} depicts the flow volume covered (CFV) in the basic DFRLP solutions to all test instances.  

\begin{table}[htb!]
\caption{Test instance characteristics.} 
    \begin{center}
\begin{tabular}[c]{c *{5}{| d{2}}} 
        \label{table:Characteristics}
        Test instance & \mathrm{s40w20} & \mathrm{s60w30} & \mathrm{s80w40} & \mathrm{s100w50} & \mathrm{Florida} \\
        \hline
        $\lvert K \rvert$     & 40         & 60         & 80         & 100         & 302 \\
        $\lvert OD \rvert$     & 20         & 30         & 40         & 50            &  74 \\
        $\lvert F \rvert$     & 190         & 435         & 780         & \num{1225}    & \num{2701} \\
        $\min \lbrace t_{O_f D_f} \rbrace$ 
                                    & 49.43     & 8.08         & 11.48     & 11.56     & 11.40 \\        
        $\max \lbrace t_{O_f D_f} \rbrace$         
                                    & 2071.06 & 2752.80& 2556.34& 2192.80 & 1257.00   \\
        $\widetilde{T}$        & 838.47     & 1152.22    &    767.06    &782.73     &  320.30 \\
        $\overline{T}$        & 907.70 & 1152.77 & 822.52    &     834.48    & 364.43    \\        
        TFV                     & \num{e06} & \num{e06} & \num{e06} & \num{e06} & \num{e12}  \\
        \hline
        \multicolumn{6}{c}{}\\[-0.4cm]
        \multicolumn{6}{c}{\footnotesize{\parbox{12cm}{\centering $\overline{T}$: average travel distance from O to D per flow; $\widetilde{T}$ median travel distance from O to D per flow;  TFV: total flow volume per instance.}}}
        \end{tabular}
    \end{center}
\end{table}

\begin{table}[htb!]
\caption{Flow volume covered in the DFRLP solutions.} 
    \begin{center}
        \begin{tabular}[c]{c *{5}{| r} }
        \label{table:solDFRLP}
        \diagbox{$p$}{Instance}
                & s40w20 & s60w30 & s80w40 & s100w50 & Florida \\
        \hline
        1         & \num{194346}     & \num{71664}         & \num{74765}         & \num{69492}     & \num{2.63e11} \\
        2         & \num{375046}    & \num{116807}     & \num{169974}          & \num{121753}     & \num{4.12e11}\\
        3         & \num{497131}     & \num{167267}     & \num{231627}     & \num{173532}     & \num{5.42e11} \\
        4         & \num{615536}     & \num{221777}     & \num{295026}     & \num{220270}     & \num{6.14e11} \\
        5         & \num{708043}     & \num{237490}     & \num{364453}     & \num{254515}     & \num{6.84e11} \\
        10     & \num{850900}     & \num{426636}     & \num{576111}     & \num{438584}     & \num{8.46e11} \\
        15     & \num{928598}     & \num{559716}     & \num{721745}     & \num{584909}     & \num{9.17e11} \\
        20     & \num{975404}     & \num{683512}     & \num{818067}     &     \num{709314}                    & \num{9.25e11} \\
        25         & \num{1000000} & \num{799192}     & \num{904855}     &     \num{815303}                    & \num{9.48e11} \\
        \hline
        \end{tabular}
    \end{center}
\end{table}

		\subsection{Computational results and performance analysis}
					\label{subsubsection:ComputationalResults}
	\begin{figure}[htb!]
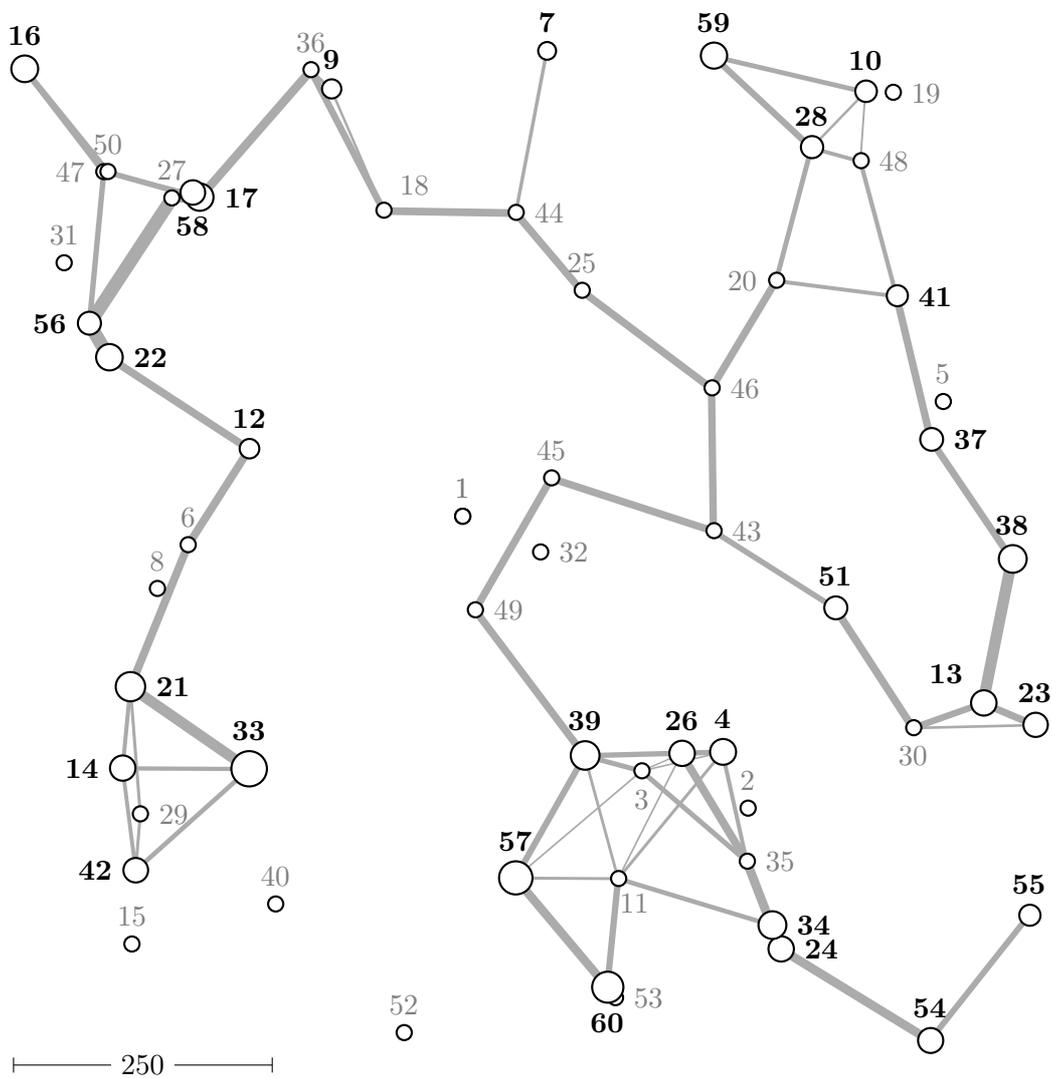

		\centering
		\begin{comment:figures}
			%
		\end{comment:figures}
		\caption{Test instance s60w30.}
		\label{figure:TestInstanceSSixtyWThirty_w/o_category}
	\end{figure}
This section discusses parameters which are required in some model extensions in more detail, followed by a numerical analysis of the performance of the model extensions described in previous sections. As the baseline case we present computational results for the instance s60w30 and refer to other instances in this section in case of results of particular interest. Otherwise, the results for the remaining test instances can be found in Appendix~\ref{Appendix:Results}.

				The test instance s60w30, which is used as the baseline case, is visualised in Figure~\ref{figure:TestInstanceSSixtyWThirty_w/o_category}. Each of the nodes is representing an origin or destination or/and a pure potential facility location. OD nodes are in bold font and the size of the OD nodes represents how many EVs are starting/ending at these locations. The thickness of the road segments between the nodes represents the proportional flow volume travelling along with these nodes.

			\subsubsection{Numerical analysis: minimum flow volume coverage (\MCFRLP)}
				\label{subsubsection:naMinCoverage}
In this section, the \MCFRLP\ is examined for different required minimum coverage levels and applied to all testing instances. Table~\ref{table:Results_MinCoverage} describes the resulting minimum number of charging stations needed to cover a prespecified proportion of all EVs for each test instance.

By reason of the occurrence of a cluster of strongly frequently nodes within the network, the test instance s40w20 differs from the other ones and is therefore described in this section in more detail.

\begin{table}[htb]
\caption{Number of open charging stations in the \MCFRLP\ solutions.} 
    \begin{center}
        \begin{tabular}[c]{c *{5}{| r} }
        \label{table:Results_MinCoverage}
          \diagbox{$C$}{Instance}
                & s40w20 & s60w30 & s80w40 & s100w50 & Florida \\
        \hline
        0.1    &    1        &    2        &    2        &    2        &    1 \\    
        0.2    &    2        &    4        &    3        &    4        &     1\\    
        0.3    &    2        &    7        &    5        &    7        &     2\\    
        0.4    &    3        &    9        &    6        &    10        &     2\\    
        0.5    &    4        &    13        &    8        &    13        &     3\\    
        0.6    &    4        &    17        &    11        &    16        &     4\\    
        0.7    &    5        &    21        &    15        &    20    &     6\\    
        0.8    &    8        &    26    &    20    &    25        &     8\\    
        0.9    &    13        &    33    &    25        &    31        &     14\\    
        1.0    &    23        &    42    &    32        &    45    & 57   \\    
        \hline
        \end{tabular}
    \end{center}
\end{table}

\paragraph{Analysing test instance s40w20.}
Taking a closer look at the resulting graph of the test case s40w20 in Figure~\ref{figure:numberOfOpenedFacilitiesDependentOnMinimumFlowCoverageForTestInstanceSFourtywTwenty}, one can determine an exponentially increasing number of charging stations in order to reach full coverage.
There seems to be a ``threshold level'' at a coverage of $\approx70\%$ of TFV, where further investments in expanding the network of charging stations become extremely expensive in relation to the resulting increase in covered EVs. This can be explained by the gradient of the line segments connecting the neighbouring points. The gradient for each segment is calculated by $\bigtriangleup p / \bigtriangleup C$. Instance s40w20 shows a gradient of at most 10 until a minimum coverage level of 70\% is required. Changing from a coverage of 70\% to 80\% of TFV results in a gradient of 30.

	\begin{figure}[htb!]
		\begin{multicols}{2}
			\begin{figurehere}
				\centering 
				\newcommand*{\maxXNumberOfOpenedFacilitiesDependentOnMinimumFlowCoverageForTestInstanceSFourtywTwenty}{1}
				\FPeval\xScaleNumberOfOpenedFacilitiesDependentOnMinimumFlowCoverageForTestInstanceSFourtywTwenty{\generalXScale / 2 / \maxXNumberOfOpenedFacilitiesDependentOnMinimumFlowCoverageForTestInstanceSFourtywTwenty}
				\newcommand*{\maxYNumberOfOpenedFacilitiesDependentOnMinimumFlowCoverageForTestInstanceSFourtywTwenty}{25}
				\FPeval\yScaleNumberOfOpenedFacilitiesDependentOnMinimumFlowCoverageForTestInstanceSFourtywTwenty{\generalYScale / \maxYNumberOfOpenedFacilitiesDependentOnMinimumFlowCoverageForTestInstanceSFourtywTwenty}
				\begin{comment:figures}
					\begin{tikzpicture}[xscale=\xScaleNumberOfOpenedFacilitiesDependentOnMinimumFlowCoverageForTestInstanceSFourtywTwenty, yscale=\yScaleNumberOfOpenedFacilitiesDependentOnMinimumFlowCoverageForTestInstanceSFourtywTwenty]
						\pgfgettransformentries{\xScaleTikz}{\@tempa}{\@tempa}{\yScaleTikz}{\@tempa}{\@tempa}
						
						\draw[very thin, color=gray, xstep = 0.2, ystep = 5] (0, 0) grid (\maxXNumberOfOpenedFacilitiesDependentOnMinimumFlowCoverageForTestInstanceSFourtywTwenty, \maxYNumberOfOpenedFacilitiesDependentOnMinimumFlowCoverageForTestInstanceSFourtywTwenty);
						
						\def\crossSizeX{\crossSize / \xScaleTikz};
						\def\crossSizeY{\crossSize / \yScaleTikz};
						
						\def\crossOne{(-\crossSizeX,-\crossSizeY) -- (\crossSizeX,\crossSizeY) (-\crossSizeX,\crossSizeY) -- (\crossSizeX,-\crossSizeY)};
						
						\draw[black, shift={(0.1, 1)}] \crossOne;
						\draw[black, shift={(0.2, 2)}] \crossOne;
						\draw[black, shift={(0.3, 2)}] \crossOne;
						\draw[black, shift={(0.4, 3)}] \crossOne;
						\draw[black, shift={(0.5, 4)}] \crossOne;
						\draw[black, shift={(0.6, 4)}] \crossOne;
						\draw[black, shift={(0.7, 5)}] \crossOne;
						\draw[black, shift={(0.8, 8)}] \crossOne;
						\draw[black, shift={(0.9, 13)}] \crossOne;
						\draw[black, shift={(1.0, 23)}] \crossOne;
						
						\draw[black!75, dotted] (0.1, 1) -- (0.2, 2) -- (0.3, 2) -- (0.4, 3) -- (0.5, 4) -- (0.6, 4) -- (0.7, 5) -- (0.8, 8) -- (0.9, 13) -- (1.0, 23);
						
						\def\axisAdditionalLengthPlusTikzX{\axisAdditionalLengthPlus / \xScaleTikz}
						\def\axisAdditionalLengthMinusTikzX{\axisAdditionalLengthMinus / \xScaleTikz}
						\draw[arrow] (-\axisAdditionalLengthMinusTikzX, 0) -- (\maxXNumberOfOpenedFacilitiesDependentOnMinimumFlowCoverageForTestInstanceSFourtywTwenty, 0) -- +(\axisAdditionalLengthPlusTikzX, 0) node[right] {$C$};
						
						\def\axisLabelTikzY{\axisLabel / \yScaleTikz}
						\foreach \pos in {0.0, 0.2, 0.4, 0.6, 0.8, 1.0} \draw[shift={(\pos, 0)}] (0, \axisLabelTikzY) -- (0, -\axisLabelTikzY) node[below] {$\pos$};
						
						\def\axisAdditionalLengthPlusTikzY{\axisAdditionalLengthPlus / \yScaleTikz}
						\def\axisAdditionalLengthMinusTikzY{\axisAdditionalLengthMinus / \yScaleTikz}
						\draw[arrow] (0, -\axisAdditionalLengthMinusTikzY) -- (0, \maxYNumberOfOpenedFacilitiesDependentOnMinimumFlowCoverageForTestInstanceSFourtywTwenty) -- +(0, \axisAdditionalLengthPlusTikzY) node[above] {$p$};
						
						\def\axisLabelTikzX{\axisLabel / \xScaleTikz}
						\foreach \pos in {0, 5, 10, 15, 20, 25} \draw[shift={(0, \pos)}] (\axisLabelTikzX, 0) -- (-\axisLabelTikzX, 0) node[left] {$\pos$};
						
					\end{tikzpicture}
				\end{comment:figures}
				\caption{Number of opened facilities dependent on minimum flow coverage for test instance: s40w20.}
				\label{figure:numberOfOpenedFacilitiesDependentOnMinimumFlowCoverageForTestInstanceSFourtywTwenty}
			\end{figurehere}
			
			\begin{figurehere}
				\centering 
				\newcommand*{\maxXNumberOfOpenedFacilitiesDependentOnMinimumFlowCoverageForTestInstanceSSixtywThirty}{1}
				\FPeval\xScaleNumberOfOpenedFacilitiesDependentOnMinimumFlowCoverageForTestInstanceSSixtywThirty{\generalXScale / 2 / \maxXNumberOfOpenedFacilitiesDependentOnMinimumFlowCoverageForTestInstanceSSixtywThirty}
				\newcommand*{\maxYNumberOfOpenedFacilitiesDependentOnMinimumFlowCoverageForTestInstanceSSixtywThirty}{50}
				\FPeval\yScaleNumberOfOpenedFacilitiesDependentOnMinimumFlowCoverageForTestInstanceSSixtywThirty{\generalYScale / \maxYNumberOfOpenedFacilitiesDependentOnMinimumFlowCoverageForTestInstanceSSixtywThirty}
				\begin{comment:figures}
					\begin{tikzpicture}[xscale=\xScaleNumberOfOpenedFacilitiesDependentOnMinimumFlowCoverageForTestInstanceSSixtywThirty, yscale=\yScaleNumberOfOpenedFacilitiesDependentOnMinimumFlowCoverageForTestInstanceSSixtywThirty]
						\pgfgettransformentries{\xScaleTikz}{\@tempa}{\@tempa}{\yScaleTikz}{\@tempa}{\@tempa}
						
						\draw[very thin, color=gray, xstep = 0.2, ystep = 10] (0, 0) grid (\maxXNumberOfOpenedFacilitiesDependentOnMinimumFlowCoverageForTestInstanceSSixtywThirty, \maxYNumberOfOpenedFacilitiesDependentOnMinimumFlowCoverageForTestInstanceSSixtywThirty);
						
						\def\crossSizeX{\crossSize / \xScaleTikz};
						\def\crossSizeY{\crossSize / \yScaleTikz};
						
						\def\crossOne{(-\crossSizeX,-\crossSizeY) -- (\crossSizeX,\crossSizeY) (-\crossSizeX,\crossSizeY) -- (\crossSizeX,-\crossSizeY)};
						
						\draw[black, shift={(0.1, 2)}] \crossOne;
						\draw[black, shift={(0.2, 4)}] \crossOne;
						\draw[black, shift={(0.3, 7)}] \crossOne;
						\draw[black, shift={(0.4, 9)}] \crossOne;
						\draw[black, shift={(0.5, 13)}] \crossOne;
						\draw[black, shift={(0.6, 17)}] \crossOne;
						\draw[black, shift={(0.7, 21)}] \crossOne;
						\draw[black, shift={(0.8, 26)}] \crossOne;
						\draw[black, shift={(0.9, 33)}] \crossOne;
						\draw[black, shift={(1.0, 42)}] \crossOne;
						
						\draw[black!75, dotted] (0.1, 2) -- (0.2, 4) -- (0.3, 7) -- (0.4, 9) -- (0.5, 13) -- (0.6, 17) -- (0.7, 21) -- (0.8, 26) -- (0.9, 33) -- (1.0, 42);
						
						\def\axisAdditionalLengthPlusTikzX{\axisAdditionalLengthPlus / \xScaleTikz}
						\def\axisAdditionalLengthMinusTikzX{\axisAdditionalLengthMinus / \xScaleTikz}
						\draw[arrow] (-\axisAdditionalLengthMinusTikzX, 0) -- (\maxXNumberOfOpenedFacilitiesDependentOnMinimumFlowCoverageForTestInstanceSSixtywThirty, 0) -- +(\axisAdditionalLengthPlusTikzX, 0) node[right] {$C$};
						
						\def\axisLabelTikzY{\axisLabel / \yScaleTikz}
						\foreach \pos in {0.0, 0.2, 0.4, 0.6, 0.8, 1.0} \draw[shift={(\pos, 0)}] (0, \axisLabelTikzY) -- (0, -\axisLabelTikzY) node[below] {$\pos$};
						
						\def\axisAdditionalLengthPlusTikzY{\axisAdditionalLengthPlus / \yScaleTikz}
						\def\axisAdditionalLengthMinusTikzY{\axisAdditionalLengthMinus / \yScaleTikz}
						\draw[arrow] (0, -\axisAdditionalLengthMinusTikzY) -- (0, \maxYNumberOfOpenedFacilitiesDependentOnMinimumFlowCoverageForTestInstanceSSixtywThirty) -- +(0, \axisAdditionalLengthPlusTikzY) node[above] {$p$};
						
						\def\axisLabelTikzX{\axisLabel / \xScaleTikz}
						\foreach \pos in {0, 10, 20, 30, 40, 50} \draw[shift={(0, \pos)}] (\axisLabelTikzX, 0) -- (-\axisLabelTikzX, 0) node[left] {$\pos$};
						
					\end{tikzpicture}
				\end{comment:figures}
				\caption{Number of opened facilities dependent on minimum flow coverage for test instance: s60w30.}
				\label{figure:numberOfOpenedFacilitiesDependentOnMinimumFlowCoverageForTestInstanceSSixtywThirty}
			\end{figurehere}
		\end{multicols}
	\end{figure}

The reason is straightforward: at a certain level, it becomes necessary to open charging stations at locations with low traffic volume in order to cover additional flows, as locations situated at strongly frequented flows are chosen in the first place. Thus, to ensure additional coverage, an increasing number of locations in ``remoted'' areas is required.

Test instance s40w20 (see Figure \ref{figure:TestInstanceSFourtyWTwenty}; shape of node is irrelevant) differs from the other ones on the basis that there is a cluster of strongly frequently nodes within the network. These nodes are nearby. Thus, a lot of flows can be covered with a small number of charging stations chosen inside of this cluster, while significantly more charging stations are necessary to cover additional flows in more remote areas. This explains the even greater rise in required charging stations between the coverage levels of 90\% and 100\% compared to the other instances. This instance is representative for sparsely populated areas with some densely populated sub-regions/cities. Building charging stations in densely populated areas cover large proportions of TFV, but guaranteeing an adequate infrastructure for all EV drivers within an area leads to a disproportionate burden.

\paragraph{Analysing the baseline case s60w30.}
The reason why instance s60w30 requires more than 50\% of all possible location in order to cover the TFV, lies in the structure of the network. The average travel distance from the origin to destination per flow, which is significantly longer in the case of instance s60w30, contributes markedly to an increasing number of required charging stations per flow. The relation between minimum coverage and required number of charging stations for instance s60w30 is depicted in Figure \ref{figure:numberOfOpenedFacilitiesDependentOnMinimumFlowCoverageForTestInstanceSSixtywThirty}. In this context, it is also important to mention that increasing the required coverage level does not necessarily result in adopting the optimal charging station placement from the previous coverage level plus opening additional charging stations. By way of illustration, take a look at the baseline instance s60w30 (see Figure~\ref{figure:TestInstanceSSixtyWThirty_w/o_category}): to cover 10\% of TFV, it is optimal to install charging stations in the locations 14 and 35. In order to set up a charging station infrastructure that can cover 20\% of all EVs that plan to travel a certain route within a network, charging stations have to be built in locations 21, 33, 35 and 42.

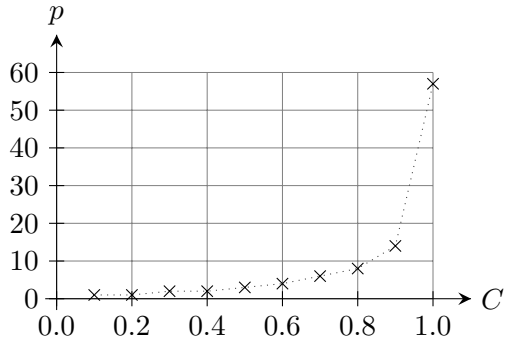
\begin{figure}[htb!]
		\centering 
		\newcommand*{\maxXNumberOfOpenedFacilitiesDependentOnMinimumFlowCoverageForTestInstanceSSixtywThirty}{1}
		\FPeval\xScaleNumberOfOpenedFacilitiesDependentOnMinimumFlowCoverageForTestInstanceSSixtywThirty{\generalXScale / 2 / \maxXNumberOfOpenedFacilitiesDependentOnMinimumFlowCoverageForTestInstanceSSixtywThirty}
		\newcommand*{\maxYNumberOfOpenedFacilitiesDependentOnMinimumFlowCoverageForTestInstanceSSixtywThirty}{60}
		\FPeval\yScaleNumberOfOpenedFacilitiesDependentOnMinimumFlowCoverageForTestInstanceSSixtywThirty{\generalYScale / \maxYNumberOfOpenedFacilitiesDependentOnMinimumFlowCoverageForTestInstanceSSixtywThirty}
		\begin{comment:figures}
			\begin{tikzpicture}[xscale=\xScaleNumberOfOpenedFacilitiesDependentOnMinimumFlowCoverageForTestInstanceSSixtywThirty, yscale=\yScaleNumberOfOpenedFacilitiesDependentOnMinimumFlowCoverageForTestInstanceSSixtywThirty]
				\pgfgettransformentries{\xScaleTikz}{\@tempa}{\@tempa}{\yScaleTikz}{\@tempa}{\@tempa}
				
				\draw[very thin, color=gray, xstep = 0.2, ystep = 10] (0, 0) grid (\maxXNumberOfOpenedFacilitiesDependentOnMinimumFlowCoverageForTestInstanceSSixtywThirty, \maxYNumberOfOpenedFacilitiesDependentOnMinimumFlowCoverageForTestInstanceSSixtywThirty);
				
				\def\crossSizeX{\crossSize / \xScaleTikz};
				\def\crossSizeY{\crossSize / \yScaleTikz};
				
				\def\crossOne{(-\crossSizeX,-\crossSizeY) -- (\crossSizeX,\crossSizeY) (-\crossSizeX,\crossSizeY) -- (\crossSizeX,-\crossSizeY)};
				
				\draw[black, shift={(0.1, 1)}] \crossOne;
				\draw[black, shift={(0.2, 1)}] \crossOne;
				\draw[black, shift={(0.3, 2)}] \crossOne;
				\draw[black, shift={(0.4, 2)}] \crossOne;
				\draw[black, shift={(0.5, 3)}] \crossOne;
				\draw[black, shift={(0.6, 4)}] \crossOne;
				\draw[black, shift={(0.7, 6)}] \crossOne;
				\draw[black, shift={(0.8, 8)}] \crossOne;
				\draw[black, shift={(0.9, 14)}] \crossOne;
				\draw[black, shift={(1.0, 57)}] \crossOne;
				
				\draw[black!75, dotted] (0.1, 1) -- (0.2, 1) -- (0.3, 2) -- (0.4, 2) -- (0.5, 3) -- (0.6, 4) -- (0.7, 6) -- (0.8, 8) -- (0.9, 14) -- (1.0, 57);
				
				\def\axisAdditionalLengthPlusTikzX{\axisAdditionalLengthPlus / \xScaleTikz}
				\def\axisAdditionalLengthMinusTikzX{\axisAdditionalLengthMinus / \xScaleTikz}
				\draw[arrow] (-\axisAdditionalLengthMinusTikzX, 0) -- (\maxXNumberOfOpenedFacilitiesDependentOnMinimumFlowCoverageForTestInstanceSSixtywThirty, 0) -- +(\axisAdditionalLengthPlusTikzX, 0) node[right] {$C$};
				
				\def\axisLabelTikzY{\axisLabel / \yScaleTikz}
				\foreach \pos in {0.0, 0.2, 0.4, 0.6, 0.8, 1.0} \draw[shift={(\pos, 0)}] (0, \axisLabelTikzY) -- (0, -\axisLabelTikzY) node[below] {$\pos$};
				
				\def\axisAdditionalLengthPlusTikzY{\axisAdditionalLengthPlus / \yScaleTikz}
				\def\axisAdditionalLengthMinusTikzY{\axisAdditionalLengthMinus / \yScaleTikz}
				\draw[arrow] (0, -\axisAdditionalLengthMinusTikzY) -- (0, \maxYNumberOfOpenedFacilitiesDependentOnMinimumFlowCoverageForTestInstanceSSixtywThirty) -- +(0, \axisAdditionalLengthPlusTikzY) node[above] {$p$};
				
				\def\axisLabelTikzX{\axisLabel / \xScaleTikz}
				\foreach \pos in {0, 10, 20, 30, 40, 50, 60} \draw[shift={(0, \pos)}] (\axisLabelTikzX, 0) -- (-\axisLabelTikzX, 0) node[left] {$\pos$};
				
			\end{tikzpicture}
		\end{comment:figures}
		\caption{Number of opened facilities dependent on minimum flow coverage for test instance: Florida.}
		\label{figure:numberOfOpenedFacilitiesDependentOnMinimumFlowCoverageForTestInstanceFlorida}
	\end{figure}

\paragraph{Analysing the Florida state highway network.}
A visualisation of the results for this instance can be seen in Figure~\ref{figure:numberOfOpenedFacilitiesDependentOnMinimumFlowCoverageForTestInstanceFlorida}. The large-scale instance representing Florida's highways is characterized by relatively short average travel distances ($\overline{T} = 364.43$). Whereas the driving range is 250, as already mentioned. Moreover, this instance allows positive flow volume for short-distance trips ($\tilde{t}_{O_f D_f} \leq 100$, where $\tilde{t}_{O_f D_f}$ stays for the real distance between origin $O_f$ and destination $D_f$). In total there are 242 out of \num{2791} flows that are characterized by a total travel distance between $O_f$ and $D_f$ less than 100. Both, the short average travel distance and the positive flow volume at short distance trips are reasons for the relatively small number of charging stations in order to cover up to 90\% of TFV.
Covering 100 \% instead of 90\% of TFV within the Florida highway case study requires a tremendously increasing number of charging stations, due to the fact that there are some flows with relatively small flow volume. (8.55\% of TFV travels along flows with total distance less than 100 length units.) Remember the TFV in this test instance is set to \num{e12} EVs.

\bigskip

Therefore, we can summarize that---at least in the early stages of infrastructure development--- it might be insufficient to force an infrastructure that is capable of covering all EVs travelling within an area. The minimum number of charging stations needed depends among other things on the length of the road segments and the average travel distance between origin and destination. Therefore, it is not possible to determine a general threshold coverage level that is applicable throughout all test instances. Figures depicting the relation between minimum coverage and the required number of charging stations for the remaining problem instances are given in Appendix~\ref{Appendix:MCFRLP} and confirm these observations.

\subsubsection{Numerical analysis: location-dependent costs per charging station (\LCFRLP)}
				\label{subsubsection:naLocationDependentCost}	
In the following, we define and discuss the parameters required for the \LCFRLP. On the one hand, location-dependent costs per charging location are introduced and on the other hand a limiting factor, in this case, a budget constraint needs to be defined. Moreover, the main insights obtained from this model extension, when analysing prespecified cost scenarios, are explained. 

\paragraph{Basic test set-up.}

				\begin{enumerate}
			    \item Partitioning facility locations into price classes:\\
The \LCFRLP\ extends the DFRLP by taking costs of installing a charging station for EVs into consideration. These costs are depending on the location of node $k$. This assumption is based on the idea that construction costs for building a charging station located along rural areas are lower than construction costs for locations in dense urban areas. Different costs of land due to the scarcity of land in urban areas justify this assumption.
One can argue, that it is more expensive to build up an infrastructure for a charging station in rural areas, but in real-world cases, when selecting potential facility locations, one criterion is that they have to be close to (smaller) villages at the countryside. Potential facilities are chose to be in locations, where it is possible to build up a station within justifiable budget. Therefore, the costs differences of land outweigh the infrastructure build-up investments. 

\medskip

The potential geographic area, where the DFRLP and its extension should be applied, is represented with a set of nodes, some of them are pure potential facility locations (=junction) and some of them are both, potential facility locations and OD nodes. For testing the model extensions, the test instances generated by \citet{DeVriesDuijzer:IncorporatingDrivingRangeVariabilityInNetworkDesignForRefuelingFacilities} are used. For simplicity and representational purposes, the individual nodes are classified into different cost categories and weighted with category specific costs. Therefore, a kind of density-based cluster algorithm is used to classify the nodes into three different categories, which can be understood as {\em urban, sub-urban} and {\em rural}. Nodes that are in close proximity to several other potential facility locations are assigned to be urban nodes, associated with the highest construction costs. Nodes that are in a single location are associated to be in rural areas, representing the lowest construction costs. Moreover, when classifying potential facility location to cost categories the total flow volume at OD nodes is taken into consideration, meaning a node where a huge number of EVs start or end their round-trip is more likely to be located in an urban area than in a rural one. 

	\begin{figure}[htb!]
		\centering
		\newcommand*{\widthPartitionOfPossibleFacilityLocationsIntoPriceClasses}{4cm}
\newcommand*{\heightPartitionOfPossibleFacilityLocationsIntoPriceClasses}{1.2cm}
		\begin{comment:figures}
			\begin{tikzpicture}[xscale=1.7, yscale=2.1,
					node/.style={black, draw, minimum width=\widthPartitionOfPossibleFacilityLocationsIntoPriceClasses, minimum height=\heightPartitionOfPossibleFacilityLocationsIntoPriceClasses, text width=\widthPartitionOfPossibleFacilityLocationsIntoPriceClasses, text centered}
				]
				\node at (0, 0.5) {For every node $k \in K$};
				\node[node] (node0) at (0, 0) {OD-node \\ $k \in O \cup D$};
				
				\node[node] (node1) at (-3, -1) {Flow volume \\ $\sum_{f \in F} \big|\{k\} \cap K_f\big| \cdot v_f$};
	
				\node[node] (node2) at (3, -1) {Radius urban \\ $P_u = \big|\{l | t_{kl \leq \pi_u}\}\big|$};
	
				\node[node] (node3) at (-1.5, -2) {Radius urban \\ $P_u = \big|\{l | t_{kl \leq \pi_u}\}\big|$};
				\node[node] (node4) at (1.5, -2) {Radius suburban \\ $P_{su} = \big|\{l | t_{kl \leq \pi_{su}}\}\big|$};
				
				\node[node] (node5) at (-3, -3) {Urban category};
				\node[node] (node6) at (0, -3) {Sub-urban category};
				\node[node] (node7) at (3, -3) {Rural category};
				
				\draw[arrow] (node0) -- (node1) node[midway, fill=white] {yes};
				\draw[arrow] (node0) -- (node2) node[midway, fill=white] {no};
				\draw[arrow] (node1) -- (node5) node[midway, fill=white, rotate=90] {$\geq$ 10\% tfv.};
				\draw[arrow] (node1) -- (node3) node[midway, fill=white] {5--10\% tfv.};
				\draw[arrow] (node1) -- (node2) node[midway, fill=white] {$<$ 5\% tfv.};
				\draw (node2) to[out=-165, in=0] (0.5, -1.4) to[out=180, in = 0] (-2, -1.4);
				\draw (-2, -1.4) .. controls (-4, -1.4) .. (node5);
				\draw[arrow] (-3.19, -2.7) -- (node5);
				\node[fill=white] at (0, -1.4) {$\geq N$};
				\draw[arrow] (node2) -- (node4) node[midway, fill=white] {$< N$};
				\draw[arrow] (node3) -- (node5) node[midway, fill=white] {$\geq N$};
				\draw[arrow] (node3) -- (node6) node[midway, fill=white] {$< N$};
				\draw[arrow] (node4) -- (node6) node[midway, fill=white] {$\geq N$};
				\draw[arrow] (node4) -- (node7) node[midway, fill=white] {$< N$};
				
			\end{tikzpicture}
		\end{comment:figures}
		\caption{Partition of possible facility locations into price classes.}
		\label{figure:PartitionOfPossibleFacilityLocationsIntoPrice}
	\end{figure}
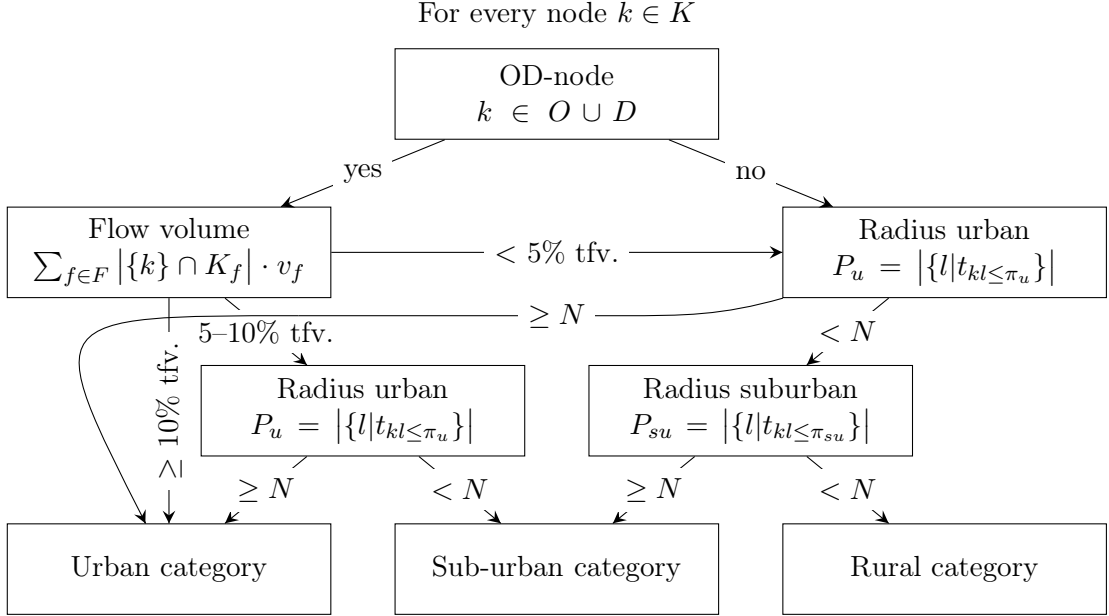

\medskip

Figure~\ref{figure:PartitionOfPossibleFacilityLocationsIntoPrice} depicts an outline of the method used to partition potential facility nodes into price classes. The parameters used for the three different price classes are summarized in Table~\ref{table:CostParameterDefinition}. For each node $k$ it is examined how many other potential locations are within a certain radius $P_u$ and $P_{su}$:
\begin{itemize}
        \item Node $k$ is associated with urban construction costs if the number of other potential facility locations within a radius of 100 length units (denoted $P_u$ in Figure \ref{figure:PartitionOfPossibleFacilityLocationsIntoPrice}) is at least $N=4$. (The threshold is chosen in order to balance the number of nodes in each category within the test instance.) 

        \item Those locations $k \in K$ that have at least $N=4$ within a radius of 150 length units (denoted $P_{su}$ in Figure \ref{figure:PartitionOfPossibleFacilityLocationsIntoPrice}) and have not yet been assigned to the urban category, are sub-urban locations.

        \item Nodes that are not in close proximity to other nodes are defines as rural nodes.
\end{itemize}

		\newpage

For OD nodes, meaning there is at least one flow starting or ending at this node, the flow volume is taken into consideration too:
\begin{itemize}
        \item If there are more than 10\% of the TFV starting or/and ending at node $k$, it is classified as an urban node, independently on the outcome of the procedure described above. 
        
        \item Similarly if the sum of all flows starting or/and ending at node $k$ is greater or equal than 5\% and smaller than 10\% of the TFV, the node is associated at least with sub-urban construction costs. 
\end{itemize} 

\begin{table}[hb!]
\caption{Characteristics of cost parameter definition.} 
    \begin{center}
        \begin{tabular}[c]{c *{4}{| r} }
        \label{table:CostParameterDefinition}
        Test instance         & s40w20 & s60w30 & s80w40 & s100w50  \\
        \hline
        $N$                         & 4            & 4            & 4            & 5         \\
        $\pi_u$                    & 100        & 100        & 100        & 100    \\
        $\pi_{su}$                & 150         & 150        & 150        & 150    \\
        urban (\#)                 & 13            & 9            & 24            & 18        \\
        sub-urban (\#)          & 11            & 27            & 28            & 38        \\
        rural (\#)                 & 16             & 24            & 28            & 48        \\
        \hline
        \end{tabular}
    \end{center}
\end{table}

Figure~\ref{figure:TestInstanceSSixtyWThirtyLCFRLPCostCategories} and 
Figures~\ref{figure:TestInstanceSFourtyWTwenty}--\ref{figure:TestInstanceSHundredWFifty} in Appendix~\ref{Appendix:Visualisation} illustrates a visualisation of the randomly generated test instances after partitioning the potential facility locations into cost categories. In the visualisation black rectangles represent urban potential facility locations, red points denote locations in sub-urban areas and rural potential locations are marked with blue circles.

	\begin{figure}[htb!]
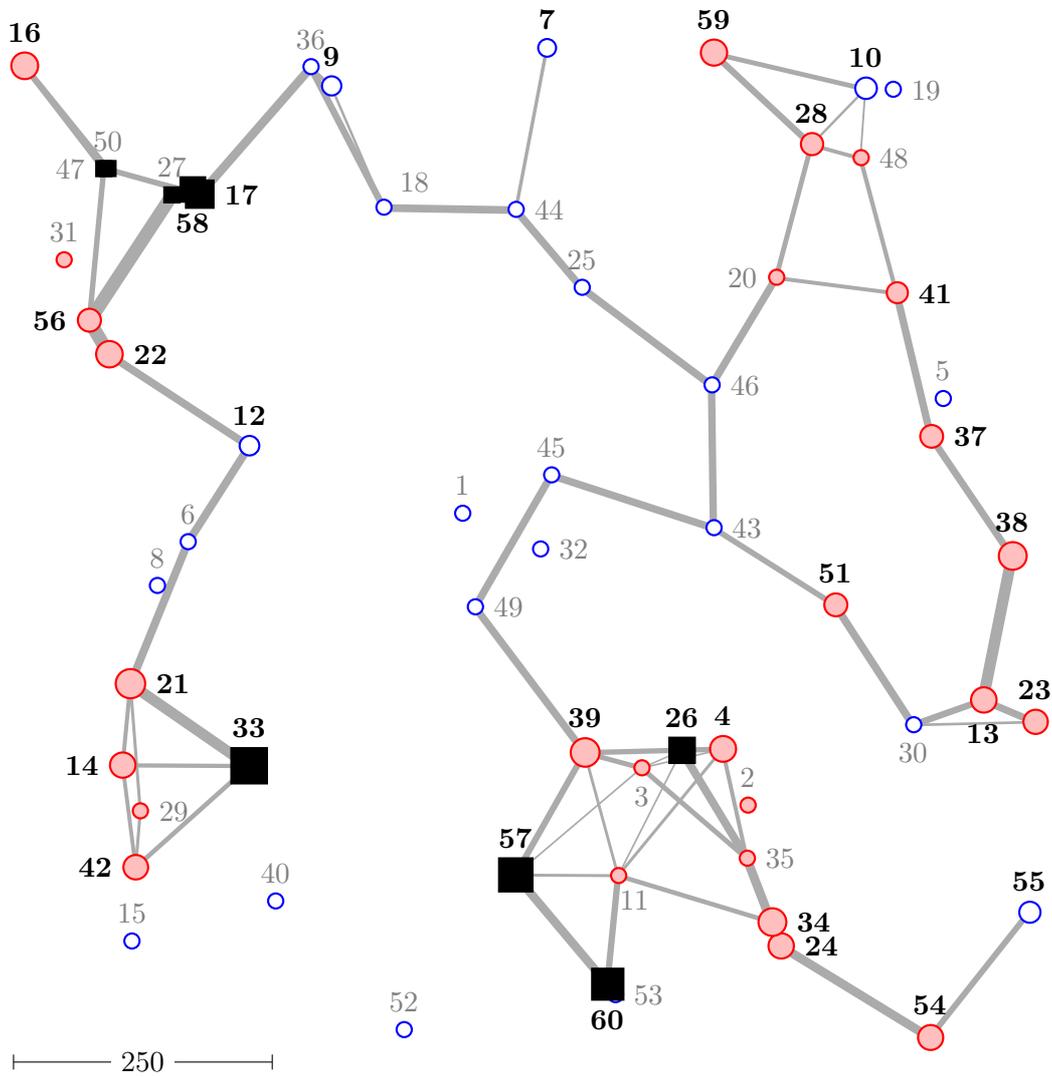

		\centering
		\begin{comment:figures}

		\end{comment:figures}
		\caption{Test instance s60w30: \LCFRLP\ -- cost categories.}
		\label{figure:TestInstanceSSixtyWThirtyLCFRLPCostCategories}
	\end{figure}

In order to avoid results of only limited informative value, we tested different cost structure scenarios by changing the proportions of the construction costs in urban, sub-urban and rural areas. Due to the fact, that construction costs can vary significantly from country to county or even from region to region, one of these scenarios might represent the ``real-world case'' in a specific area. Moreover, this set of scenarios can be interpreted as possible set of options for subsidy systems a government can choose from. The public sector can grant subsidies in order to change a given cost structure and therefore enhance investments in the development of an adequate infrastructure. Depending on their extent, subsidies support the deployment of charging stations in certain areas. 
In the test cases, location-dependent costs are chosen within an interval of $[1,7]$, based on the idea to depict special cases, like installing a charging station in a sub-urban region is more than twice as expensive as in a rural area. While the costs of construction in urban areas are more than twice as expensive as those in a sub-urban area. Tests are applied for 15 combinations of construction cost scenarios which are depicted in Figure~\ref{figure:differentChargingCostScenariosForChargingStations}.

			\begin{figure}[htb!]
				\centering 
				\newcommand*{\maxXDifferentChargingCostScenariosForChargingStations}{15}
				\FPeval\xScaleDifferentChargingCostScenariosForChargingStations{\generalXScale / 0.824 / \maxXDifferentChargingCostScenariosForChargingStations}
				\newcommand*{\maxYDifferentChargingCostScenariosForChargingStations}{7}
				\FPeval\yScaleDifferentChargingCostScenariosForChargingStations{\generalYScale / \maxYDifferentChargingCostScenariosForChargingStations}
				\begin{comment:figures}
					\begin{tikzpicture}[xscale=\xScaleDifferentChargingCostScenariosForChargingStations, yscale=\yScaleDifferentChargingCostScenariosForChargingStations]
						\pgfgettransformentries{\xScaleTikz}{\@tempa}{\@tempa}{\yScaleTikz}{\@tempa}{\@tempa}
						
						\draw[very thin, color=gray, xstep = 1, ystep = 1] (0, 0) grid (\maxXDifferentChargingCostScenariosForChargingStations, \maxYDifferentChargingCostScenariosForChargingStations);
						
						\def\axisLabelTikzY{\axisLabel / \yScaleTikz}
						\def\axisPartitionTikzY{5.5 * \axisLabelTikzY}
						\def\axisLabelTikzX{\axisLabel / \xScaleTikz}
						\draw[thick, color=gray] (5.5, -\axisPartitionTikzY) -- (5.5, 7.0 + \axisPartitionTikzY);
						\draw[thick, color=gray] (9.5, -\axisPartitionTikzY) -- (9.5, 7.0 + \axisPartitionTikzY);
						\draw[thick, color=gray] (12.5, -\axisPartitionTikzY) -- (12.5, 7.0 + \axisPartitionTikzY);
						\draw[thick, color=gray] (14.5, -\axisPartitionTikzY) -- (14.5, 7.0 + \axisPartitionTikzY);
						
						\node[circle, draw=blue!100, fill=white!100, thick, inner sep=0pt, minimum size=2mm] at (1, 1) {};
						\node[circle, draw=red!100, fill=red!25, thick, inner sep=0pt, minimum size=2mm] at (1, 2) {};
						\node[rectangle, draw=black!100, fill=black!100, thick, inner sep=0pt, minimum size=2mm] at (1, 7) {};
						
						\node[circle, draw=blue!100, fill=white!100, thick, inner sep=0pt, minimum size=2mm] at (2, 1) {};
						\node[circle, draw=red!100, fill=red!25, thick, inner sep=0pt, minimum size=2mm] at (2, 3) {};
						\node[rectangle, draw=black!100, fill=black!100, thick, inner sep=0pt, minimum size=2mm] at (2, 7) {};
						
						\node[circle, draw=blue!100, fill=white!100, thick, inner sep=0pt, minimum size=2mm] at (3, 1) {};
						\node[circle, draw=red!100, fill=red!25, thick, inner sep=0pt, minimum size=2mm] at (3, 4) {};
						\node[rectangle, draw=black!100, fill=black!100, thick, inner sep=0pt, minimum size=2mm] at (3, 7) {};
						
						\node[circle, draw=blue!100, fill=white!100, thick, inner sep=0pt, minimum size=2mm] at (4, 1) {};
						\node[circle, draw=red!100, fill=red!25, thick, inner sep=0pt, minimum size=2mm] at (4, 5) {};
						\node[rectangle, draw=black!100, fill=black!100, thick, inner sep=0pt, minimum size=2mm] at (4, 7) {};
						
						\node[circle, draw=blue!100, fill=white!100, thick, inner sep=0pt, minimum size=2mm] at (5, 1) {};
						\node[circle, draw=red!100, fill=red!25, thick, inner sep=0pt, minimum size=2mm] at (5, 6) {};
						\node[rectangle, draw=black!100, fill=black!100, thick, inner sep=0pt, minimum size=2mm] at (5, 7) {};
						
						\node[circle, draw=blue!100, fill=white!100, thick, inner sep=0pt, minimum size=2mm] at (6, 2) {};
						\node[circle, draw=red!100, fill=red!25, thick, inner sep=0pt, minimum size=2mm] at (6, 3) {};
						\node[rectangle, draw=black!100, fill=black!100, thick, inner sep=0pt, minimum size=2mm] at (6, 7) {};
						
						\node[circle, draw=blue!100, fill=white!100, thick, inner sep=0pt, minimum size=2mm] at (7, 2) {};
						\node[circle, draw=red!100, fill=red!25, thick, inner sep=0pt, minimum size=2mm] at (7, 4) {};
						\node[rectangle, draw=black!100, fill=black!100, thick, inner sep=0pt, minimum size=2mm] at (7, 7) {};
						
						\node[circle, draw=blue!100, fill=white!100, thick, inner sep=0pt, minimum size=2mm] at (8, 2) {};
						\node[circle, draw=red!100, fill=red!25, thick, inner sep=0pt, minimum size=2mm] at (8, 5) {};
						\node[rectangle, draw=black!100, fill=black!100, thick, inner sep=0pt, minimum size=2mm] at (8, 7) {};
						
						\node[circle, draw=blue!100, fill=white!100, thick, inner sep=0pt, minimum size=2mm] at (9, 2) {};
						\node[circle, draw=red!100, fill=red!25, thick, inner sep=0pt, minimum size=2mm] at (9, 6) {};
						\node[rectangle, draw=black!100, fill=black!100, thick, inner sep=0pt, minimum size=2mm] at (9, 7) {};
						
						\node[circle, draw=blue!100, fill=white!100, thick, inner sep=0pt, minimum size=2mm] at (10, 3) {};
						\node[circle, draw=red!100, fill=red!25, thick, inner sep=0pt, minimum size=2mm] at (10, 4) {};
						\node[rectangle, draw=black!100, fill=black!100, thick, inner sep=0pt, minimum size=2mm] at (10, 7) {};
						
						\node[circle, draw=blue!100, fill=white!100, thick, inner sep=0pt, minimum size=2mm] at (11, 3) {};
						\node[circle, draw=red!100, fill=red!25, thick, inner sep=0pt, minimum size=2mm] at (11, 5) {};
						\node[rectangle, draw=black!100, fill=black!100, thick, inner sep=0pt, minimum size=2mm] at (11, 7) {};
						
						\node[circle, draw=blue!100, fill=white!100, thick, inner sep=0pt, minimum size=2mm] at (12, 3) {};
						\node[circle, draw=red!100, fill=red!25, thick, inner sep=0pt, minimum size=2mm] at (12, 6) {};
						\node[rectangle, draw=black!100, fill=black!100, thick, inner sep=0pt, minimum size=2mm] at (12, 7) {};
						
						\node[circle, draw=blue!100, fill=white!100, thick, inner sep=0pt, minimum size=2mm] at (13, 4) {};
						\node[circle, draw=red!100, fill=red!25, thick, inner sep=0pt, minimum size=2mm] at (13, 5) {};
						\node[rectangle, draw=black!100, fill=black!100, thick, inner sep=0pt, minimum size=2mm] at (13, 7) {};
						
						\node[circle, draw=blue!100, fill=white!100, thick, inner sep=0pt, minimum size=2mm] at (14, 4) {};
						\node[circle, draw=red!100, fill=red!25, thick, inner sep=0pt, minimum size=2mm] at (14, 6) {};
						\node[rectangle, draw=black!100, fill=black!100, thick, inner sep=0pt, minimum size=2mm] at (14, 7) {};
						
						\node[circle, draw=blue!100, fill=white!100, thick, inner sep=0pt, minimum size=2mm] at (15, 5) {};
						\node[circle, draw=red!100, fill=red!25, thick, inner sep=0pt, minimum size=2mm] at (15, 6) {};
						\node[rectangle, draw=black!100, fill=black!100, thick, inner sep=0pt, minimum size=2mm] at (15, 7) {};
						
						\def\axisAdditionalLengthPlusTikzX{\axisAdditionalLengthPlus / \xScaleTikz}
						\def\axisAdditionalLengthMinusTikzX{\axisAdditionalLengthMinus / \xScaleTikz}
						\draw[-] (-\axisAdditionalLengthMinusTikzX, 0) -- (\maxXDifferentChargingCostScenariosForChargingStations, 0) -- +(0, 0) node[right] {Scenario};
						
						\foreach \pos in {1, 2, 3, 4, 5, 6, 7, 8, 9, 10, 11, 12, 13, 14, 15} \draw[shift={(\pos, 0)}] (0, \axisLabelTikzY) -- (0, -\axisLabelTikzY) node[below] {$\pos$};
						
						\def\axisAdditionalLengthPlusTikzY{\axisAdditionalLengthPlus / \yScaleTikz}
						\def\axisAdditionalLengthMinusTikzY{\axisAdditionalLengthMinus / \yScaleTikz}
						\draw[arrow] (0, -\axisAdditionalLengthMinusTikzY) -- (0, \maxYDifferentChargingCostScenariosForChargingStations) -- +(0, \axisAdditionalLengthPlusTikzY) node[above] {Costs};
						
						\foreach \pos in {0, 1, 2, 3, 4, 5, 6, 7} \draw[shift={(0, \pos)}] (\axisLabelTikzX, 0) -- (-\axisLabelTikzX, 0) node[left] {$\pos$};
						
					\end{tikzpicture}
				\end{comment:figures}
				\caption[Different charging cost scenarios for charging stations.]{Different charging cost scenarios for charging stations. \begin{tikzpicture}\node[rectangle, draw=black!100, fill=black!100, thick, inner sep=0pt, minimum size=2mm] {};\end{tikzpicture} stays for urban, \begin{tikzpicture}\node[circle, draw=red!100, fill=red!25, thick, inner sep=0pt, minimum size=2mm] {};\end{tikzpicture} for suburban, and \begin{tikzpicture}\node[circle, draw=blue!100, fill=white!100, thick, inner sep=0pt, minimum size=2mm] {};\end{tikzpicture} for rural areas.}
				\label{figure:differentChargingCostScenariosForChargingStations}
			\end{figure}

			  \item Definition of the budget:\\
The model takes a limited budget into account. In the testing scenarios for the \LCFRLP\ the available budget $B$ for installing a charging station infrastructure depends, on the one hand on the instance due to the partitioning of potential facility locations into prices classes and, on the other one on the scenario itself due to cost differences concerning the cost categories. Based on preliminary tests, the budget is chosen in order to be sufficient to open 25\% of all potential facilities. The following equation represents the calculation of budget, where $c_u, c_{su}, c_r$ describe the construction costs in urban, sub-urban and rural areas, respectively: 
        $$
        B = (\# urban\ nodes \cdot c_u + \# sub\text{-}urban\ nodes \cdot c_{su} + \# rural\ nodes \cdot c_r) \cdot 0.25
        $$
The restricted budget for investments in the deployment of a charging station infrastructure is depicted in Table~\ref{table:BudgetDefinition} for each scenario and test instance.

\begin{table}[htb]
\caption{Budget for all scenarios and test instances.} 
    \begin{center}
            \begin{tabular}[c]{c *{9}{| r}  }
                \label{table:BudgetDefinition}
        \diagbox{Inst.}{Scen.}
               &    1    &    2    &    3    &    4    &    5    &    6    &    7    &    8    &    9    \\
                \hline
            s40w20&    32.25&    35&    37.75&    40.5    &    43.25&    39&    41.75&    44.5    &47.25\\
            s60w30&    35.25&    42&    48.75&    55.5    &    62.25&    48&    54.75&    61.5    &68.25\\
            s80w40&    63    &    70    &    77    &    84    &    91    &    77    &    84    &    91    &    98\\
            s100w50&62.5    &    72    &    81.5    &    91    &    100.5&    84&    93.5    &    103    &    112.5     \\
        \hline
        \multicolumn{10}{c}{}\\[-0.5em]
        \diagbox{Inst.}{Scen.}
             &    10    &    11    &    12    &    13    &    14    &    \multicolumn{1}{r}{15} \\
         \cline{1-7}
        s40w20     & 45.75&    48.5        &    51.25    &    52.5    &    55.25    &    \multicolumn{1}{r}{59.25} \\
        s60w30     & 60.75&    67.5    &    74.25    &    73.5    &    80.25    &    \multicolumn{1}{r}{86.25} \\
        s80w40     & 91    &98&    105    &    105    &    112    &    \multicolumn{1}{r}{119} \\
        s100w50    & 105.5    &    115&    124.5    &    127    &    136.5&    \multicolumn{1}{r}{148.5} \\
        \cline{1-7}
            \end{tabular}
    \end{center}
\end{table}

		\newpage

    \item Testing for different cost scenarios:\\
    Of course, every two scenarios can be compared. Nevertheless, better insights can be obtained by considering some particular scenario sequences:
    \begin{itemize}
    		\item First, all pairs differing in costs of charging stations (CS) in only one category can give insights into the behaviour of the \LCFRLP. Thus we define the following two categories (S.1 and S.2), each of them containing four sequences.
			\renewcommand{\labelenumii}{S.\arabic{enumii}:}
			\begin{enumerate}
				\item 1-2-3-4-5, 6-7-8-9, 10-11-12, and 13-14: in these sequences, denoted S.\arabic{enumii}.1, S.\arabic{enumii}.2, S.\arabic{enumii}.3 and S.\arabic{enumii}.4, respectively, the costs of urban and rural CSs are constant while the costs of sub-urban CSs gradually increase by 1. 
				\item 5-9-12-14-15, 4-8-11-13, 3-7-10, and 2-6: here, the costs of urban and sub-urban CSs do not change, but the costs of rural CSs gradually increase. These sequences are denoted by S.\arabic{enumii}.1, S.\arabic{enumii}.2, S.\arabic{enumii}.3, and S.\arabic{enumii}.4, respectively.
			\end{enumerate}
			\renewcommand{\labelenumii}{\theenumii. }
		\item Sometimes, it might be informative to compare situations which are ``similar'' with respect to the general cost structure, but differ in more than one cost category.
			\renewcommand{\labelenumii}{S.\arabic{enumii}:}
			\begin{enumerate}
				\setcounter{enumii}{2}
				\item 1-6-10-13-15: all these scenarios have one thing common: the costs of CSs in sub-urban and rural locations are similar while the costs of those in urban areas are much higher. Note that this property gets less significant when moving from scenario 1 towards scenario 15.
			\end{enumerate}
			\renewcommand{\labelenumii}{\theenumii. }
	\end{itemize}
	
	Apart from the sequences described above, some particular scenarios are of special interest:
	\begin{itemize}
		\item Scenario 1 depicts the extreme situation where the sub-urban and rural CS costs are very small compared to those of the urban ones.
		\item In scenario 5 the costs of CSs in urban and sub-urban areas are relatively similar, but the costs in rural locations are much smaller.
		\item Scenario 15 corresponds more or less to the case where all CSs have relatively similar costs.
	\end{itemize}
	
	Finally, scenarios, which lead to the same total number of opened charging stations, can be of special interest.
   
\end{enumerate}

The visualizations shown in this section are to be interpreted as follows: surrounded nodes indicate opened charging stations and dark marked paths indicate covered flows, where the dashed lines represent the proportional coverage. Flows that cannot be covered are depicted as dashed lines in light color.

\paragraph{Analysing baseline case s60w30.}
Before dealing with the particular cost scenarios, one important property should be outlined. Consider \eg the scenario sequence S.1.1: since budget depends on the costs of all CSs which can be placed, the total budget grows by moving from scenario 1 towards scenario 5. Nevertheless, this does not mean that an optimal solution for one scenario is automatically feasible for the next one.

In particular, budget is calculated by weighting each node with costs and taking 25\%. Therefore, the change in budget, by moving from one scenario to the next one in the S.1 scenario sequence, is calculated by taking the amount of sub-urban nodes (\# = 11) multiplied with the cost difference (\# = 1) and taking 25\% of it. Similarly, it works for the S.2 scenario sequences. Thus, if construction costs of a category increase by one unit, there have to be at least 4 possible location nodes of the increasing cost category within the network, otherwise there will not be enough budget to install the very same charging stations in the scenario where costs increased. In the following, we will speak of a {\em threshold} if targeting this issue. Consequently, it can happen that the CFV of an optimal solution decreases while moving in the S.1 and S.2 scenario sequences towards their end.

Of course, this issue becomes less significant if much more possible facility locations than installed CSs exist and it could be overcome by allowing more locations (\eg not only in the junction points). In such a case, a different categorization (based on areas) would be necessary. Nevertheless, budget, which is set to 25\% of total costs, makes it possible to identify the main trends.

  			\paragraph{Discussion of the results for baseline case s60w30.}

\addtolength{\tabcolsep}{-1.5pt}
\begin{table}[htb]
\caption{Results of \LCFRLP\ for all scenarios of s60w30.} 
	\begin{center}
			\begin{tabular}[c]{c *{10}{| r} }
				\label{table:table:ResultsLCFRLPsSixtywThirty}
Scenario  &	\multicolumn{2}{c|}{1}	&	\multicolumn{2}{c|}{2}	&	\multicolumn{2}{c|}{3}	&	\multicolumn{2}{c|}{4}	&	\multicolumn{2}{c}{5}\\
                \hline
$CFV$  &	\multicolumn{2}{r|}{\num{518977}}	&	\multicolumn{2}{r|}{\num{466190}}	&	\multicolumn{2}{r|}{\num{432579}}	&	\multicolumn{2}{r|}{\num{436391}}	&	\multicolumn{2}{r}{\num{443510}}\\

$\sum_{k \in K} x_k$   &	\multicolumn{2}{r|}{16}	&	\multicolumn{2}{r|}{14}	&	\multicolumn{2}{r|}{12}	&	\multicolumn{2}{r|}{11}	&	\multicolumn{2}{r}{13}\\
				\hline
				\hline
urban (\sfrac{\% }{100} | \#)	&	0.0625 & 1	&	0.0714 & 1 &	0.1667 & 2	&	0.3636 & 4	&	0.3077	& 4\\
sub-urban (\sfrac{\% }{100} | \#)	&	0.8125 & 13	&	0.7857 & 11	&	0.6667	 & 8 &	0.4545	& 5 &	0.3846	& 5\\
rural	(\sfrac{\% }{100} | \#) &	0.125 & 2 	&	0.1429 & 2	&	0.1667 & 2	&	0.1818	& 2 &	0.3077& 4\\
				\hline
				\multicolumn{11}{c}{}\\[-0.5em]
Scenario  &	\multicolumn{2}{c|}{6}	&	\multicolumn{2}{c|}{7}	&	\multicolumn{2}{c|}{8}	&	\multicolumn{2}{c|}{9}	&	\multicolumn{2}{c}{10}\\
                \hline
$CFV$  &		\multicolumn{2}{r|}{\num{499044}}	&	\multicolumn{2}{r|}{\num{455130}}	&	\multicolumn{2}{r|}{\num{450888}}	&	\multicolumn{2}{r|}{\num{462837}}	&	\multicolumn{2}{r}{\num{495546}}\\

$\sum_{k \in K} x_k$   &	\multicolumn{2}{r|}{15}	&	\multicolumn{2}{r|}{13}	&	\multicolumn{2}{r|}{12}	&	\multicolumn{2}{r|}{12}	&	\multicolumn{2}{r}{13}\\

				\hline
				\hline
urban (\sfrac{\% }{100} | \#)	&	0.0667 &1 	&	0.1538 & 2	&	0.25	 & 3 &	0.3333 & 4	&	0.2308 & 3\\
sub-urban (\sfrac{\% }{100} | \#)	&	0.8	& 12 &	0.6923	& 9 &	0.5833 & 7	&	0.5 & 6	&	0.6923	& 9\\
rural	(\sfrac{\% }{100} | \#) &	0.1333 & 2	&	0.1538	& 2 &	0.1667 & 2 &	0.1667	& 2 &	0.0769 & 1\\
				\hline
				\multicolumn{11}{c}{}\\[-0.5em]
Scenario  &	\multicolumn{2}{c|}{11}	&	\multicolumn{2}{c|}{12}&	\multicolumn{2}{c|}{13}	&	\multicolumn{2}{c|}{14}	&	\multicolumn{2}{c}{15} \\
                \hline
$CFV$  &		\multicolumn{2}{r|}{\num{479880}}	&	\multicolumn{2}{r|}{\num{472890}}&	\multicolumn{2}{r|}{\num{513453}}	&	\multicolumn{2}{r|}{\num{502720}}	&	\multicolumn{2}{r}{\num{523208}} \\

$\sum_{k \in K} x_k$   &	\multicolumn{2}{r|}{13}	&	\multicolumn{2}{r|}{12}&	\multicolumn{2}{r|}{13}	&	\multicolumn{2}{r|}{13}	&	\multicolumn{2}{r}{14} \\

				\hline
				\hline
urban (\sfrac{\% }{100} | \#)	&	0.2308	 & 3 &	0.3333 & 4	&	0.3077	& 4 &	0.3077& 4 	&	0.2857 & 4\\
sub-urban (\sfrac{\% }{100} | \#)	&	0.6154	& 8 &	0.5833 & 7	&	0.6923	& 9 &	0.6154 & 8 	&	0.5714 & 8\\
rural	(\sfrac{\% }{100} | \#) &	0.1538 & 2 	&	0.0833 & 1	&	0	& 0 &	0.0769	& 1 &	0.1429 & 2\\
				\hline
			\end{tabular}
	\end{center}
\end{table}
\addtolength{\tabcolsep}{+1.5pt}

			\begin{figure}[htb!]
				\centering 
				\newcommand*{\maxXProportionOfUrbanSuburbanRuralStatinsOpenedInDifferentScenariosForTheTestInstanceSSixtyWThirty}{15}
				\FPeval\xScaleProportionOfUrbanSuburbanRuralStatinsOpenedInDifferentScenariosForTheTestInstanceSSixtyWThirty{\generalXScale / 0.838 / \maxXProportionOfUrbanSuburbanRuralStatinsOpenedInDifferentScenariosForTheTestInstanceSSixtyWThirty}
				\newcommand*{\maxYProportionOfUrbanSuburbanRuralStatinsOpenedInDifferentScenariosForTheTestInstanceSSixtyWThirtyNumber}{18.0}
				\FPeval\yScaleProportionOfUrbanSuburbanRuralStatinsOpenedInDifferentScenariosForTheTestInstanceSSixtyWThirtyNumber{\generalYScale / \maxYProportionOfUrbanSuburbanRuralStatinsOpenedInDifferentScenariosForTheTestInstanceSSixtyWThirtyNumber}
				
\newcommand*{\maxYProportionOfUrbanSuburbanRuralStatinsOpenedInDifferentScenariosForTheTestInstanceSSixtyWThirtyProportion}{1.0}
				\FPeval\yScaleProportionOfUrbanSuburbanRuralStatinsOpenedInDifferentScenariosForTheTestInstanceSSixtyWThirtyProportion{\generalYScale / \maxYProportionOfUrbanSuburbanRuralStatinsOpenedInDifferentScenariosForTheTestInstanceSSixtyWThirtyProportion}
				\begin{comment:figures}
					\begin{tikzpicture}[xscale=\xScaleProportionOfUrbanSuburbanRuralStatinsOpenedInDifferentScenariosForTheTestInstanceSSixtyWThirty, yscale=\yScaleProportionOfUrbanSuburbanRuralStatinsOpenedInDifferentScenariosForTheTestInstanceSSixtyWThirtyNumber]
						\pgfgettransformentries{\xScaleTikz}{\@tempa}{\@tempa}{\yScaleTikz}{\@tempa}{\@tempa}
						
						\draw[very thin, color=gray, xstep = 1, ystep = 2] (0, 0) grid (\maxXProportionOfUrbanSuburbanRuralStatinsOpenedInDifferentScenariosForTheTestInstanceSSixtyWThirty, \maxYProportionOfUrbanSuburbanRuralStatinsOpenedInDifferentScenariosForTheTestInstanceSSixtyWThirtyNumber);
						
						\def\axisLabelTikzY{\axisLabel / \yScaleTikz}
						\def\axisNumberTikzY{5.5 * \axisLabelTikzY}
						\def\axisTotalNumberTikzY{3 * \axisLabelTikzY}
						\def\axisLabelTikzX{\axisLabel / \xScaleTikz}
						\draw[thick, color=gray] (5.0, -\axisNumberTikzY) -- (5.0, 18.0 + \axisNumberTikzY);
						\draw[thick, color=gray] (9.0, -\axisNumberTikzY) -- (9.0, 18.0 + \axisNumberTikzY);
						\draw[thick, color=gray] (12.0, -\axisNumberTikzY) -- (12.0, 18.0 + \axisNumberTikzY);
						\draw[thick, color=gray] (14.0, -\axisNumberTikzY) -- (14.0, 18.0 + \axisNumberTikzY);
						
						\begin{scope}[shift={(-0.5, 0)}]
							\draw[black, fill=black] (0.6, 0.0) rectangle (0.8, 1);
							\draw[red, thick, fill=red!25] (0.9, 0.0) rectangle (1.1, 13);
							\draw[blue, thick] (1.2, 0.0) rectangle (1.4, 2);
							
							\draw[black, fill=black] (1.6, 0.0) rectangle (1.8, 1);
							\draw[red, thick, fill=red!25] (1.9, 0.0) rectangle (2.1, 11);
							\draw[blue, thick] (2.2, 0.0) rectangle (2.4, 2);
							
							\draw[black, fill=black] (2.6, 0.0) rectangle (2.8, 2);
							\draw[red, thick, fill=red!25] (2.9, 0.0) rectangle (3.1, 8);
							\draw[blue, thick] (3.2, 0.0) rectangle (3.4, 2);
							
							\draw[black, fill=black] (3.6, 0.0) rectangle (3.8, 4);
							\draw[red, thick, fill=red!25] (3.9, 0.0) rectangle (4.1, 5);
							\draw[blue, thick] (4.2, 0.0) rectangle (4.4, 2);
							
							\draw[black, fill=black] (4.6, 0.0) rectangle (4.8, 4);
							\draw[red, thick, fill=red!25] (4.9, 0.0) rectangle (5.1, 5);
							\draw[blue, thick] (5.2, 0.0) rectangle (5.4, 4);
							
							\draw[black, fill=black] (5.6, 0.0) rectangle (5.8, 1);
							\draw[red, thick, fill=red!25] (5.9, 0.0) rectangle (6.1, 12);
							\draw[blue, thick] (6.2, 0.0) rectangle (6.4, 2);
							
							\draw[black, fill=black] (6.6, 0.0) rectangle (6.8, 2);
							\draw[red, thick, fill=red!25] (6.9, 0.0) rectangle (7.1, 9);
							\draw[blue, thick] (7.2, 0.0) rectangle (7.4, 2);
							
							\draw[black, fill=black] (7.6, 0.0) rectangle (7.8, 3);
							\draw[red, thick, fill=red!25] (7.9, 0.0) rectangle (8.1, 7);
							\draw[blue, thick] (8.2, 0.0) rectangle (8.4, 2);
							
							\draw[black, fill=black] (8.6, 0.0) rectangle (8.8, 4);
							\draw[red, thick, fill=red!25] (8.9, 0.0) rectangle (9.1, 6);
							\draw[blue, thick] (9.2, 0.0) rectangle (9.4, 2);
							
							\draw[black, fill=black] (9.6, 0.0) rectangle (9.8, 3);
							\draw[red, thick, fill=red!25] (9.9, 0.0) rectangle (10.1, 9);
							\draw[blue, thick] (10.2, 0.0) rectangle (10.4, 1);
							
							\draw[black, fill=black] (10.6, 0.0) rectangle (10.8, 3);
							\draw[red, thick, fill=red!25] (10.9, 0.0) rectangle (11.1, 8);
							\draw[blue, thick] (11.2, 0.0) rectangle (11.4, 2);
							
							\draw[black, fill=black] (11.6, 0.0) rectangle (11.8, 4);
							\draw[red, thick, fill=red!25] (11.9, 0.0) rectangle (12.1, 7);
							\draw[blue, thick] (12.2, 0.0) rectangle (12.4, 1);
							
							\draw[black, fill=black] (12.6, 0.0) rectangle (12.8, 4);
							\draw[red, thick, fill=red!25] (12.9, 0.0) rectangle (13.1, 9);
							\draw[blue, thick] (13.2, 0.0) rectangle (13.4, 0);
							
							\draw[black, fill=black] (13.6, 0.0) rectangle (13.8, 4);
							\draw[red, thick, fill=red!25] (13.9, 0.0) rectangle (14.1, 8);
							\draw[blue, thick] (14.2, 0.0) rectangle (14.4, 1);
							
							\draw[black, fill=black] (14.6, 0.0) rectangle (14.8, 4);
							\draw[red, thick, fill=red!25] (14.9, 0.0) rectangle (15.1, 8);
							\draw[blue, thick] (15.2, 0.0) rectangle (15.4, 2);
						\end{scope}
						
						\def\axisAdditionalLengthPlusTikzX{\axisAdditionalLengthPlus / \xScaleTikz}
						\def\axisAdditionalLengthMinusTikzX{\axisAdditionalLengthMinus / \xScaleTikz}
						\draw[-] (-\axisAdditionalLengthMinusTikzX, 0) -- (\maxXProportionOfUrbanSuburbanRuralStatinsOpenedInDifferentScenariosForTheTestInstanceSSixtyWThirty, 0) -- +(0, 0) node[right] {Scenario};
						
						\foreach \pos in {1, 2, 3, 4, 5, 6, 7, 8, 9, 10, 11, 12, 13, 14, 15} \draw[shift={(\pos-0.5, 0)}] (0, -\axisLabelTikzY) node[below] {$\pos$};
						
						\def\axisAdditionalLengthPlusTikzY{\axisAdditionalLengthPlus / \yScaleTikz}
						\def\axisAdditionalLengthMinusTikzY{\axisAdditionalLengthMinus / \yScaleTikz}
						\draw[arrow] (0, -\axisAdditionalLengthMinusTikzY) -- (0, \maxYProportionOfUrbanSuburbanRuralStatinsOpenedInDifferentScenariosForTheTestInstanceSSixtyWThirtyNumber) -- +(0, \axisAdditionalLengthPlusTikzY) node[above] {Total number};
						
						\foreach \pos in {2, 4, 6, 8, 10, 12, 14, 16, 18} \draw[shift={(0, \pos)}] (\axisLabelTikzX, 0) -- (-\axisLabelTikzX, 0) node[left] {$\pos$};
						
					\end{tikzpicture}
					
					
					\hspace*{0.118cm}
					\begin{tikzpicture}[xscale=\xScaleProportionOfUrbanSuburbanRuralStatinsOpenedInDifferentScenariosForTheTestInstanceSSixtyWThirty, yscale=\yScaleProportionOfUrbanSuburbanRuralStatinsOpenedInDifferentScenariosForTheTestInstanceSSixtyWThirtyProportion]
						\pgfgettransformentries{\xScaleTikz}{\@tempa}{\@tempa}{\yScaleTikz}{\@tempa}{\@tempa}
						
						\draw[very thin, color=gray, xstep = 1, ystep = 0.2] (0, 0) grid (\maxXProportionOfUrbanSuburbanRuralStatinsOpenedInDifferentScenariosForTheTestInstanceSSixtyWThirty, \maxYProportionOfUrbanSuburbanRuralStatinsOpenedInDifferentScenariosForTheTestInstanceSSixtyWThirtyProportion);
						
						\def\axisLabelTikzY{\axisLabel / \yScaleTikz}
						\def\axisPartitionTikzY{5.5 * \axisLabelTikzY}
						\def\axisLabelTikzX{\axisLabel / \xScaleTikz}
						\draw[thick, color=gray] (5.0, -\axisPartitionTikzY) -- (5.0, 1.0 + \axisPartitionTikzY);
						\draw[thick, color=gray] (9.0, -\axisPartitionTikzY) -- (9.0, 1.0 + \axisPartitionTikzY);
						\draw[thick, color=gray] (12.0, -\axisPartitionTikzY) -- (12.0, 1.0 + \axisPartitionTikzY);
						\draw[thick, color=gray] (14.0, -\axisPartitionTikzY) -- (14.0, 1.0 + \axisPartitionTikzY);
						
						\begin{scope}[shift={(-0.5, 0)}]
							\draw[black, fill=black] (0.6, 0.0) rectangle (0.8, 0.0625);
							\draw[red, thick, fill=red!25] (0.9, 0.0) rectangle (1.1, 0.8125);
							\draw[blue, thick] (1.2, 0.0) rectangle (1.4, 0.125);
							
							\draw[black, fill=black] (1.6, 0.0) rectangle (1.8, 0.0714285714285714);
							\draw[red, thick, fill=red!25] (1.9, 0.0) rectangle (2.1, 0.785714285714286);
							\draw[blue, thick] (2.2, 0.0) rectangle (2.4, 0.142857142857143);
							
							\draw[black, fill=black] (2.6, 0.0) rectangle (2.8, 0.166666666666667);
							\draw[red, thick, fill=red!25] (2.9, 0.0) rectangle (3.1, 0.666666666666667);
							\draw[blue, thick] (3.2, 0.0) rectangle (3.4, 0.166666666666667);
							
							\draw[black, fill=black] (3.6, 0.0) rectangle (3.8, 0.363636363636364);
							\draw[red, thick, fill=red!25] (3.9, 0.0) rectangle (4.1, 0.454545454545455);
							\draw[blue, thick] (4.2, 0.0) rectangle (4.4, 0.181818181818182);
							
							\draw[black, fill=black] (4.6, 0.0) rectangle (4.8, 0.307692307692308);
							\draw[red, thick, fill=red!25] (4.9, 0.0) rectangle (5.1, 0.384615384615385);
							\draw[blue, thick] (5.2, 0.0) rectangle (5.4, 0.307692307692308);
							
							\draw[black, fill=black] (5.6, 0.0) rectangle (5.8, 0.0666666666666667);
							\draw[red, thick, fill=red!25] (5.9, 0.0) rectangle (6.1, 0.8);
							\draw[blue, thick] (6.2, 0.0) rectangle (6.4, 0.133333333333333);
							
							\draw[black, fill=black] (6.6, 0.0) rectangle (6.8, 0.153846153846154);
							\draw[red, thick, fill=red!25] (6.9, 0.0) rectangle (7.1, 0.692307692307692);
							\draw[blue, thick] (7.2, 0.0) rectangle (7.4, 0.153846153846154);
							
							\draw[black, fill=black] (7.6, 0.0) rectangle (7.8, 0.25);
							\draw[red, thick, fill=red!25] (7.9, 0.0) rectangle (8.1, 0.583333333333333);
							\draw[blue, thick] (8.2, 0.0) rectangle (8.4, 0.166666666666667);
							
							\draw[black, fill=black] (8.6, 0.0) rectangle (8.8, 0.333333333333333);
							\draw[red, thick, fill=red!25] (8.9, 0.0) rectangle (9.1, 0.5);
							\draw[blue, thick] (9.2, 0.0) rectangle (9.4, 0.166666666666667);
							
							\draw[black, fill=black] (9.6, 0.0) rectangle (9.8, 0.230769230769231);
							\draw[red, thick, fill=red!25] (9.9, 0.0) rectangle (10.1, 0.692307692307692);
							\draw[blue, thick] (10.2, 0.0) rectangle (10.4, 0.0769230769230769);
							
							\draw[black, fill=black] (10.6, 0.0) rectangle (10.8, 0.230769230769231);
							\draw[red, thick, fill=red!25] (10.9, 0.0) rectangle (11.1, 0.615384615384615);
							\draw[blue, thick] (11.2, 0.0) rectangle (11.4, 0.153846153846154);
							
							\draw[black, fill=black] (11.6, 0.0) rectangle (11.8, 0.333333333333333);
							\draw[red, thick, fill=red!25] (11.9, 0.0) rectangle (12.1, 0.583333333333333);
							\draw[blue, thick] (12.2, 0.0) rectangle (12.4, 0.0833333333333333);
							
							\draw[black, fill=black] (12.6, 0.0) rectangle (12.8, 0.307692307692308);
							\draw[red, thick, fill=red!25] (12.9, 0.0) rectangle (13.1, 0.692307692307692);
							\draw[blue, thick] (13.2, 0.0) rectangle (13.4, 0);
							
							\draw[black, fill=black] (13.6, 0.0) rectangle (13.8, 0.307692307692308);
							\draw[red, thick, fill=red!25] (13.9, 0.0) rectangle (14.1, 0.615384615384615);
							\draw[blue, thick] (14.2, 0.0) rectangle (14.4, 0.0769230769230769);
							
							\draw[black, fill=black] (14.6, 0.0) rectangle (14.8, 0.285714285714286);
							\draw[red, thick, fill=red!25] (14.9, 0.0) rectangle (15.1, 0.571428571428571);
							\draw[blue, thick] (15.2, 0.0) rectangle (15.4, 0.142857142857143);
						\end{scope}

						\def\axisAdditionalLengthPlusTikzX{\axisAdditionalLengthPlus / \xScaleTikz}
						\def\axisAdditionalLengthMinusTikzX{\axisAdditionalLengthMinus / \xScaleTikz}
						\draw[-] (-\axisAdditionalLengthMinusTikzX, 0) -- (\maxXProportionOfUrbanSuburbanRuralStatinsOpenedInDifferentScenariosForTheTestInstanceSSixtyWThirty, 0) -- +(0, 0) node[right] {Scenario};
						
						\foreach \pos in {1, 2, 3, 4, 5, 6, 7, 8, 9, 10, 11, 12, 13, 14, 15} \draw[shift={(\pos-0.5, 0)}] (0, -\axisLabelTikzY) node[below] {$\pos$};
						
						\def\axisAdditionalLengthPlusTikzY{\axisAdditionalLengthPlus / \yScaleTikz}
						\def\axisAdditionalLengthMinusTikzY{\axisAdditionalLengthMinus / \yScaleTikz}
						\draw[arrow] (0, -\axisAdditionalLengthMinusTikzY) -- (0, \maxYProportionOfUrbanSuburbanRuralStatinsOpenedInDifferentScenariosForTheTestInstanceSSixtyWThirtyProportion) -- +(0, \axisAdditionalLengthPlusTikzY) node[above] {Proportion};
						
						\foreach \pos in {0.0, 0.2, 0.4, 0.6, 0.8, 1.0} \draw[shift={(0, \pos)}] (\axisLabelTikzX, 0) -- (-\axisLabelTikzX, 0) node[left] {$\pos$};
						
					\end{tikzpicture}
				\end{comment:figures}
				\caption[Total number and proportion of urban/suburban/rural statins opened in different scenarios for the test instance s60w30.]{Total number and proportion of urban/suburban/rural statins opened in different scenarios for the test instance s60w30. The left bar \begin{tikzpicture}[xscale=\xScaleProportionOfUrbanSuburbanRuralStatinsOpenedInDifferentScenariosForTheTestInstanceSSixtyWThirty, yscale=1]\draw[black, fill=black] (0.0, 0.0) rectangle (0.2, 0.24);\end{tikzpicture} stays for urban, the middle bar \begin{tikzpicture}[xscale=\xScaleProportionOfUrbanSuburbanRuralStatinsOpenedInDifferentScenariosForTheTestInstanceSSixtyWThirty, yscale=1]\draw[red, thick, fill=red!25] (0.0, 0.0) rectangle (0.2, 0.24);\end{tikzpicture} for suburban, and the right bar \begin{tikzpicture}[xscale=\xScaleProportionOfUrbanSuburbanRuralStatinsOpenedInDifferentScenariosForTheTestInstanceSSixtyWThirty, yscale=1]\draw[blue, thick] (0.0, 0.0) rectangle (0.2, 0.24);\end{tikzpicture} for rural areas.}
				\label{figure:totalNumberAndProportionOfUrbanSuburbanRuralStatinsOpenedInDifferentScenariosForTheTestInstanceSSixtyWThirty}
			\end{figure}

The resulting flow coverage and 
the proportion of opened urban/suburban/rural stations for the test instance s60w30 are described in Table~\ref{table:table:ResultsLCFRLPsSixtywThirty} and depicted in Figure \ref{figure:totalNumberAndProportionOfUrbanSuburbanRuralStatinsOpenedInDifferentScenariosForTheTestInstanceSSixtyWThirty}.

\subparagraph{S.1 scenario sequences.} As sub-urban construction costs increase under the assumption that urban and rural costs remain the same, the proportion of sub-urban CSs is decreasing while the proportion of both, urban and rural CSs is increasing. Reason for it is the enlarged budget (which---as explained above---increases by 25\% of the number of sub-urban locations in every move towards the sequence end).

As sub-urban costs increase, urban locations become relative less expensive. Therefore, the leftover budget, resulting from a bigger budget or from the failure of building a more expensive sub-urban CS, can be used to build urban and/or rural CSs.

Moreover, it might become optimal to give up a sub-urban location that was built in a previous scenario and use the leftover budget to build one or more urban or rural stations, which are relatively cheaper due to the increasing budget.

Finally, note that sometimes it can be optimal to build rural stations even if they are relatively expensive and that the relation ``urban better than sub-urban better than rural'' is not always true. Compared to the other test instances, s60w30 is characterized by longer average travel distances between origins and destinations. Thus, rural locations are sometimes of particular interest because they are essential to cover long distance flows. 

Let us now illustrate that an optimal solution in one scenario is not necessarily feasible in the next one when moving towards the end of a S.1 or S.2 scenario sequence. Changing from scenario 1 to scenario 2 increases budget due to an increase in sub-urban costs, but the number of open CSs decreases and so does the coverage level. In particular, the budget increases by 25\% of the total number of sub-urban locations ($27 \cdot 1 \cdot 0.25 = 6.75$). In scenario 1, CFV is maximized by opening one urban CS, 13 sub-urban and 2 rural ones. The budget required to install the sub-urban CSs is 26 units. In scenario 2, one would need a budget of 39 units to build 13 sub-urban CSs, but budget is limited to 42 units in this scenario. Therefore, it is necessary to renounce the construction of sub-urban CSs in order to stay within budget limitation. Building 8 sub-urban CSs leaves a remaining budget of 9 units which can be used to install one urban and two rural CSs.

	\begin{figure}[htb!]
		\centering
			\begin{comment:figures}

			\end{comment:figures}
		\caption{Test instance s60w30: \LCFRLP\ -- scenario 13.}
		\label{figure:TestInstanceSSixtyWThirtyLCFRLPScenarioThirteen}
	\end{figure}
	\begin{figure}[htb!]
		\centering
			\begin{comment:figures}
				%
			\end{comment:figures}
		\caption{Test instance s60w30: \LCFRLP\ -- scenario 14.}
		\label{figure:TestInstanceSSixtyWThirtyLCFRLPScenarioFourteen}
	\end{figure}
	\begin{figure}[htb!]
		\centering
			\begin{comment:figures}
				%
			\end{comment:figures}
		\caption{Test instance s60w30: \LCFRLP\ -- scenario 15.}
		\label{figure:TestInstanceSSixtyWThirtyLCFRLPScenarioFifteen}
	\end{figure}

\subparagraph{S.2 scenario sequences.} In these scenario sequences cost differences between sub-urban and urban locations are constant, rural costs, however, are significantly lower but increasing when gradually changing towards the sequence end. A trend for a decreasing proportion of rural CS can be observed when looking at absolute and proportional numbers. This behaviour emerges from the fact that building CSs in rural locations becomes more expensive and at rural locations, which are often situated on the shortest path of long distance flows, less frequented flows are usually passing by.

There is another general phenomenon that is specific to these scenario sequences: CFV is increasing continuously when moving towards the sequence end. A continuous increase is possible due to the fact that rural costs are the only changing component and therefore the urban and sub-urban CSs (which are usually preferred) can be taken over from one scenario into the next one.

Consider scenarios 14 and 15, depicted in Figures \ref{figure:TestInstanceSSixtyWThirtyLCFRLPScenarioFourteen} and \ref{figure:TestInstanceSSixtyWThirtyLCFRLPScenarioFifteen}, respectively: looking at Table~\ref{table:table:ResultsLCFRLPsSixtywThirty}, it can be seen that the absolute number of urban and sub-urban CSs remains the same, while one additional rural CS is built in scenario 15. However, by inspecting Figures \ref{figure:TestInstanceSSixtyWThirtyLCFRLPScenarioFourteen} and \ref{figure:TestInstanceSSixtyWThirtyLCFRLPScenarioFifteen} it can be see that the position of one sub-urban CS changed if moving from scenario 14 to scenario 15. \Ie sometimes it might become optimal to change the positioning of CSs in order to increase CFV by covering different flows than in the previous scenario.

\subparagraph{Scenario sequence S.3.} The general trend is an increasing proportion of urban locations and decreasing proportions of sub-urban and rural locations. The reason for this stems from the decreasing cost differences between sub-urban and urban locations, resulting in a trade-off between rural/sub-urban locations and urban locations, as it becomes relatively cheaper to build a CS in an urban location. Due to the fact that more frequented flows are passing by urban locations, CFV will usually increase by installing a CS in an urban area.

There are only small changes between scenarios 1 and 6 and scenarios 10 and 13. The comparison of scenarios 6 and 10, however, shows a far more strongly change in the structure of optimal CS placement. Contrary to scenario 6, in scenario 10 building two sub-urban CSs no longer incurs lower costs than building one urban CS ($4 \cdot 2 = 8 > 7$). Usually urban locations are stronger frequented than sub-urban ones and therefore CFV can be sometimes increased when changing CS placement and using current budget in order to install more CSs in urban locations by giving up sub-urban and rural ones. But still in scenario 10 leftover budget that cannot be used to build urban or sub-urban CSs is used to build rural CSs. This trade-off becomes apparent in Table~\ref{table:table:ResultsLCFRLPsSixtywThirty} when looking at the absolute numbers of CSs in each category in scenario 6 and 10: there are 1 urban, 12 sub-urban and 2 rural CSs in scenario 6 (total costs: $7 + 36 + 4 = 47 \leq B\colon 48$) and 3 urban, 9 sub-urban and 1 rural CSs in scenario 10 (total costs: $21 + 36 + 3 = 60 \leq B\colon 60.75$).

Comparing scenario 13 and 15 (see Figures \ref{figure:TestInstanceSSixtyWThirtyLCFRLPScenarioThirteen} and \ref{figure:TestInstanceSSixtyWThirtyLCFRLPScenarioFifteen}), one would assume the same general trend: increasing proportion of urban CS and decreasing proportion in sub-urban and rural CSs. However this is not the case. In total, there are 13 CSs built in scenario 13 and 14 CSs in scenario 15. In scenario 13, budget is used to build 9 CSs in sub-urban regions and 4 CSs in urban locations. It is easy to check that it would be possible to install the very same infrastructure in scenario 15 as it was in scenario 13. But doing so results in a leftover budget of 4 units, which is not enough to build an additional rural CS. This obviously results in the same coverage level. But by exchanging a sub-urban CS and use the corresponding budget of 6 units to build two rural CSs, the CFV could be increased.

\subparagraph{Scenarios leading to the same total number of opened stations.} The scenario pairs 5-14, 7-10, 8-9, 9-12, 10-11, 11-13, and 13-14 result in the same total number of opened stations in this test instance. Apart from the pair 5-14, these scenarios differ in the costs of CSs in only one category and in the usual case, the increasing budget is used to move a CS into a better location, \ie for moving one CS from a rural location to a sub-urban one or from a sub-urban location to an urban one.

Consider scenarios 8 and 9 more closely: the total number of opened CSs remains the same (\# = 12). The difference between these two scenarios is that sub-urban location costs increase by one unit and the urban locations become relatively cheaper. The optimal infrastructure in scenario 8 causes total costs of 60, where a budget of 61 is available. Installing the very same infrastructure under scenario 9 causes total costs of 67 and therefore covers the very same amount of EVs. In scenario 9 there is an available budget of 68 units. Giving up a sub-urban location and building an urban CS instead is within the budget constraint of scenario 9 and increases CFV.

Sometimes it can also happen that increasing costs of CSs in sub-urban locations does not allow using the old solution as it is the case in the scenario pair 10-11: the same number of CSs is installed (\# = 13). The difference between these two scenarios is that sub-urban location costs increase by one unit. The optimal infrastructure in scenario 10 requires full budget (60 units), as sub-urban costs increase by one unit, it is no longer possible to install the very same infrastructure under scenario 11. Thus, the CFV decreases. The budget is sufficient to install the same amount of urban CSs and one sub-urban CS less than in scenario 10. The resulting leftover budget can be used to build two rural CSs in order to maximize CFV. 

\subparagraph{Scenario 15.}
This scenario is characterized by similar cost differences between rural, sub-urban and urban locations. Moreover, the margin between rural and urban construction costs is small. This results in a CS placement that is close to a placement where no cost differences exist. Let us now compare this scenario with the extreme case, denoted {\em scenario 16} in this paragraph, where every CS costs 7 units. The budget in scenario 16 increases strongly (by 19 units), as rural costs  and sub-urban costs increase by 2 and 1 units, respectively. Building the very same infrastructure of CSs of scenario 15 in scenario 16 causes total costs of 98 units. The leftover budget of 7 units can be used to build an additional CS. But due to the fact, that there are equal costs for all locations, CFV can rise more strongly by changing the structure of CSs placement and situate 15 CSs in locations that maximize CFV. Results for scenario 16 are depicted in Table~\ref{table:table:ResultsLCFRLPsSixtywThirtyScenSixteen}.

\begin{table}[htb]
\caption{Results of \LCFRLP\ for scenarios 15 and 16 of s60w30.} 
	\begin{center}
			\begin{tabular}[c]{c *{4}{| r} }
				\label{table:table:ResultsLCFRLPsSixtywThirtyScenSixteen}
Scenario  &	\multicolumn{2}{c|}{15} & \multicolumn{2}{c}{16}	\\
                \hline
$CFV$  & \multicolumn{2}{r|}{523208} &	\multicolumn{2}{r}{559716} \\
$\sum_{k \in K} x_k$  &	\multicolumn{2}{r|}{14}  & \multicolumn{2}{r}{15} \\
				\hline
				\hline
urban (\sfrac{\%}{100} | \#)	&	0.2857 & 4  & 0.2667 & 4 \\
sub-urban (\sfrac{\%}{100} | \#)	&	0.5714 & 8 &0.6667 & 10  \\
rural	(\sfrac{\%}{100} | \#) &	0.1429 & 2 & 0.0667& 1\\
		\hline
			\end{tabular}
	\end{center}
\end{table}

Finally note that scenario 15 (followed by scenario 1) return the highest CFV.

\medskip

The analysis of the results of the \LCFRLP\ provides main insight that the cost structure of urban, sub-urban and rural costs has a strong influence to the optimal charging station infrastructure. As a result, policymakers are able to enhance investments in deployment of an adequate CS infrastructure by developing and designing efficient subsidy systems. Governments have the possibility to use subsidies to influence the current cost structure of urban, sub-urban and rural location and subsequently change the structure of CS locations.


				\subsubsection{Numerical analysis: determination of the station size (\CFRLP)}
				\label{subsubsection:naCFRLP}	
This section focuses on the analysis of the \CFRLP, while repeating the testing process for different exogenously given numbers of charging poles to locate. This provides an opportunity to observe how charging poles were located in different situations. The results are then analysed in terms of CFV and number of charging station locations.

\paragraph{Basic test set-up:}

\begin{enumerate}
\item Define capacity per charging pole:  \\
        First, the energy consumption at each possible facility location to cover the TFV is calculated. Some potential facility locations are not located on the shortest path of any flow, energy demand at these locations is therefore zero. In order to define a capacity limit for CSs, we calculate the median energy demand per charging location, excluding those CSs that have a zero energy demand. Based on observations, no more than four charging poles are usually installed at a location ($M = 4$). To define the capacity of a charging pole, the median capacity per charging station is divided by the maximum number of charging poles per location which is set to 4. The resulting instance-specific capacity is shown in Table~\ref{table:CapacityPerPole}. (A scaling parameter of 0.001 is used.)
        
        \begin{table}[htb]
        \caption{Capacity per charging pole.} 
            \begin{center}
                \begin{tabular}[c]{c *{1}{| r} }
                \label{table:CapacityPerPole}
                Test instance     & $Cap$    \\
                \hline
                s40w20            & \num{2801.01} \\
                s60w30            & \num{3362.01} \\
                s80w40            & \num{1678.57} \\
                s100w50            & \num{1282.72} \\
                \hline
                \end{tabular}
            \end{center}
        \end{table}

\item Pre-testing: \\
        In order to get information about a reasonable test range for the number of charging poles to locate we do some pre-testing. First, to obtain information about the maximum possible CFV, given the capacity limitation of charging poles, the \CFRLP\ is solved, requiring that a maximum sized charging station ($M = 4$) is built at every possible location. Therefore, the number of charging poles to locate is 4 times the number of potential facility locations. 
        Since the same CFV can be guaranteed with a lower number of charging poles, the model extension \CMCFRLP\ is applied for further pre-testing. In this extension, described in Section~\ref{subsubsection:Endogenous}, the number of charging poles to locate is no longer exogenously given, but the objective is to minimize the number of charging poles needed to cover a certain proportion of EVs. The objective value from the first pre-test, representing the maximum possible CFV, given the capacity limitation, is taken as an input for the second pre-test and defines the minimum required coverage level. The objective value of the second pre-test represents the maximum reasonable number of charging poles given the capacity limits of charging poles. Locating more charging poles, would not increase the number of EVs, able to complete their round-trips successfully.

\begin{example}
Assessing the results of the first pre-test (see Table~\ref{table:Results_CMCFRLP_pre_s60w30}), where the maximum number of charging poles ($60 \cdot 4 = 240$) is located, we can be said that with the previously defined capacity of \num{3362,01} per charging pole it is not possible to cover more than \num{827525} (82.75\%) EVs of TFV. Using this proportion of TFV as an input for the \CMCFRLP, 
the observations by the means of the second pre-test can be summarized, by saying that locating more than 144 charging poles would not be an economically viable decision given the limited capacity of charging poles, as the number of EVs, that can complete their round-trips successfully, cannot be increased by installing further charging poles.
\end{example}

\begin{table}[htb]
\caption{Pre-test (\CMCFRLP) of instance s60w30.} 
    \begin{center}
            \begin{tabular}[c]{c *{1}{| r} }
                \label{table:Results_CMCFRLP_pre_s60w30}
                $C$                            & 0.827525  \\ 
                 \hline            
                 $\sum_{k \in K} n_k$ & 144\\    
                 $\sum_{f \in f} x_k$     & 48\\    
                 PCF                             & 188\\    
                 solve time                 & \num{1598.64}\\ 
                    \hline
			        \multicolumn{2}{c}{}\\[-0.4cm]
                    \multicolumn{2}{c}{\footnotesize{PCF: number of partially covered flows ($z_f > 0$).}}    
            \end{tabular}
    \end{center}
\end{table}

\item Testing for different $S$ values:\\
        In the \CFRLP\ the exogenously given number of charging poles located within the network can be interpreted as kind of a budget constraint.  We repeat the testing process for different numbers of charging poles, representing 25\%, 50\%, 75\% and 100\% of the maximum economically viable number of charging poles (EVCP).
\end{enumerate}

Solve time increases rapidly when the number of charging poles to locate decreases, but it seems to be concavely curved. 

Testing the \CFRLP\ for $S=108$ with Gurobi solver and using an option that depicts the progress of the branch and bound by turning on the Gurobi MIP logfile, shows that the incumbent value did not change in the last \num{3600}~seconds. Therefore, the solving process got interrupted after \num{4600}~seconds of elapsed time 
and a time limit of \num{1800}~seconds for the following testing process with $S = 72$ and $S=36$ is set. In the case of allocating 72 charging poles, the incumbent value did not change within the last 200 seconds when approaching the end ot the time limit. Locating 36 charging poles, the incumbent value remains the same for the last \num{1300}~seconds.

 \footnotetext[5]{Best bound \num{794839,20}, gap 0.6729\%.}
 \footnotetext[6]{Best bound \num{761047,02}, gap 13.1263\%.} 
 \footnotetext[7]{Best bound \num{506922,03}, gap 28.4685\%.}
 
\medskip

In this sections' visualizations, open charging stations are marked with a surrounding circle, where the size of the circle represents the number of installed charging poles per location (1, 2, 3 or 4). Covered paths are indicated with dark dashed lines representing the proportional coverage. Flows that cannot be covered given the charging station infrastructure are dotted in light color.

\paragraph{Analysing the baseline case s60w30.}
In the following, the main goal is to obtain insight into the results of the \CFRLP\ for the baseline case s60w30, where capacity is chosen to be \num{3362,01} (see Table \ref{table:CapacityPerPole}). Computational results are shown in Table~\ref{table:Results_CFRLP_s60w30} and described in the following.

\begin{table}[htb!]
\caption{Results of \CFRLP\ for instance s60w30.} 
    \begin{center}
            \begin{tabular}[c]{c *{4}{| r} }
                \label{table:Results_CFRLP_s60w30}
                $S$ & 144 (100\%) &  108$^*$ (75\%)     &    72$^*$ (50\%)    & 36$^*$ (25\%) \\
                 \hline    
                 CFV        & \num{827524.96} &\num{789526.33} \footnotemark &\num{672741.07} \footnotemark & \num{394588.66} \footnotemark \\                
                 $\sum_{k \in K} x_k$ 	& 48             & 45              & 32             & 21 \\   
                 $\overline{n}$ 			& 3.00			& 2.40			& 2.25		& 1.70\\
                 $\sum_{f \in f} y_f$  	& 416        & 241              & 118            & 38\\    
                 PCF                             	& 188        & 109             & 48                & 23\\    
                 solve time                	& 850.07& \num{5255.77} & \num{1803.48}& \num{1803.40}\\ 
                    \hline
			        \multicolumn{5}{c}{}\\[-0.4cm]
                    \multicolumn{5}{c}{\footnotesize{$\overline{n}$: average station size.}} \\   
                    \multicolumn{5}{c}{\footnotesize{$^*$ indicates no optimal results, interrupted solving process.}}    
            \end{tabular}
    \end{center}
\end{table}    

\begin{figure}[htb!]
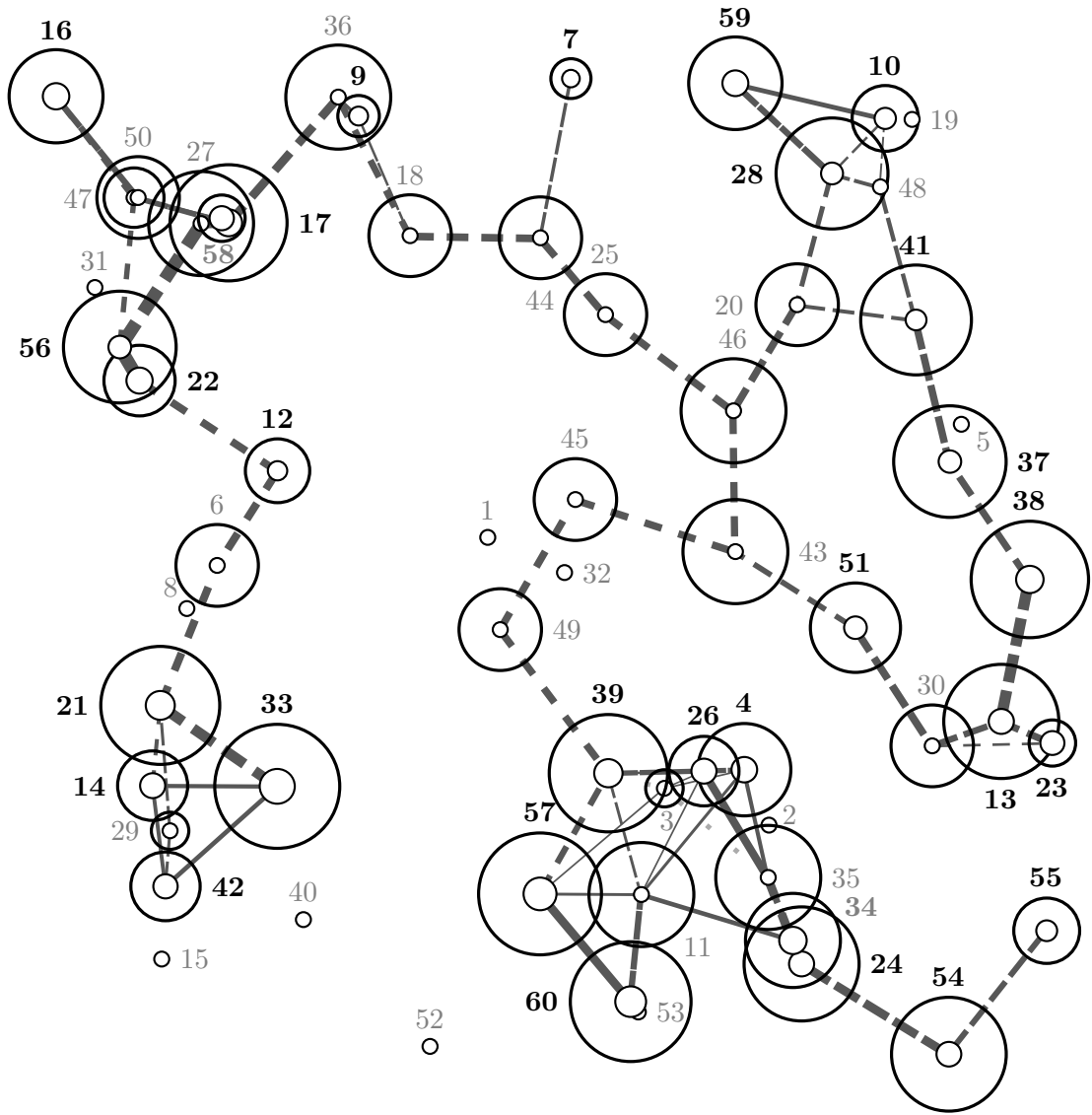

		\centering
			\begin{comment:figures}

			\end{comment:figures}
		\caption{Test instance s60w30: \CFRLP\ -- $S = 144$  (100\%).}
		\label{figure:TestInstanceSSixtyWThirtyCFRLPSOneHundretAndFourtyFour}
	\end{figure}

While in the uncapacitated FRLP the entire TFV within a network can always be covered if enough charging stations are built, this is no longer possible in the capacitated FRLP. 
Pre-tests indicates that given the capacity limitation of charging poles, it is not possible to cover more than \num{827524,96} of \num{e06} (82.75\%) EVs, denoting the TFV of this test instance (see Table~\ref{table:Results_CMCFRLP_pre_s60w30}). Therefore, it is necessary to locate 144 charging poles (see Figure~\ref{figure:TestInstanceSSixtyWThirtyCFRLPSOneHundretAndFourtyFour}), where on average 3 charging poles are installed per station. Given the locations of charging stations, 416 of 435 flows could be covered, however given the limited capacity at charging poles proportions of 188 flows can actually refuel and therefore complete their round-trips. 

Locating 108 (75\% of EVCP) charging poles still covers 95.41\% of the maximum number of EVs that can be covered given the capacity limit of a charging pole. These charging poles are located in 45 possible facility locations. Compared to the previously described case, there are on average 2.4 charging poles per station.

In case of locating 72 charging poles, it is still possible to cover 81.30\% of the maximum number of EVs that can be covered given the capacity limit of a charging pole. This represents an allocation of 50\% of EVCP. 

\begin{figure}[htb!]
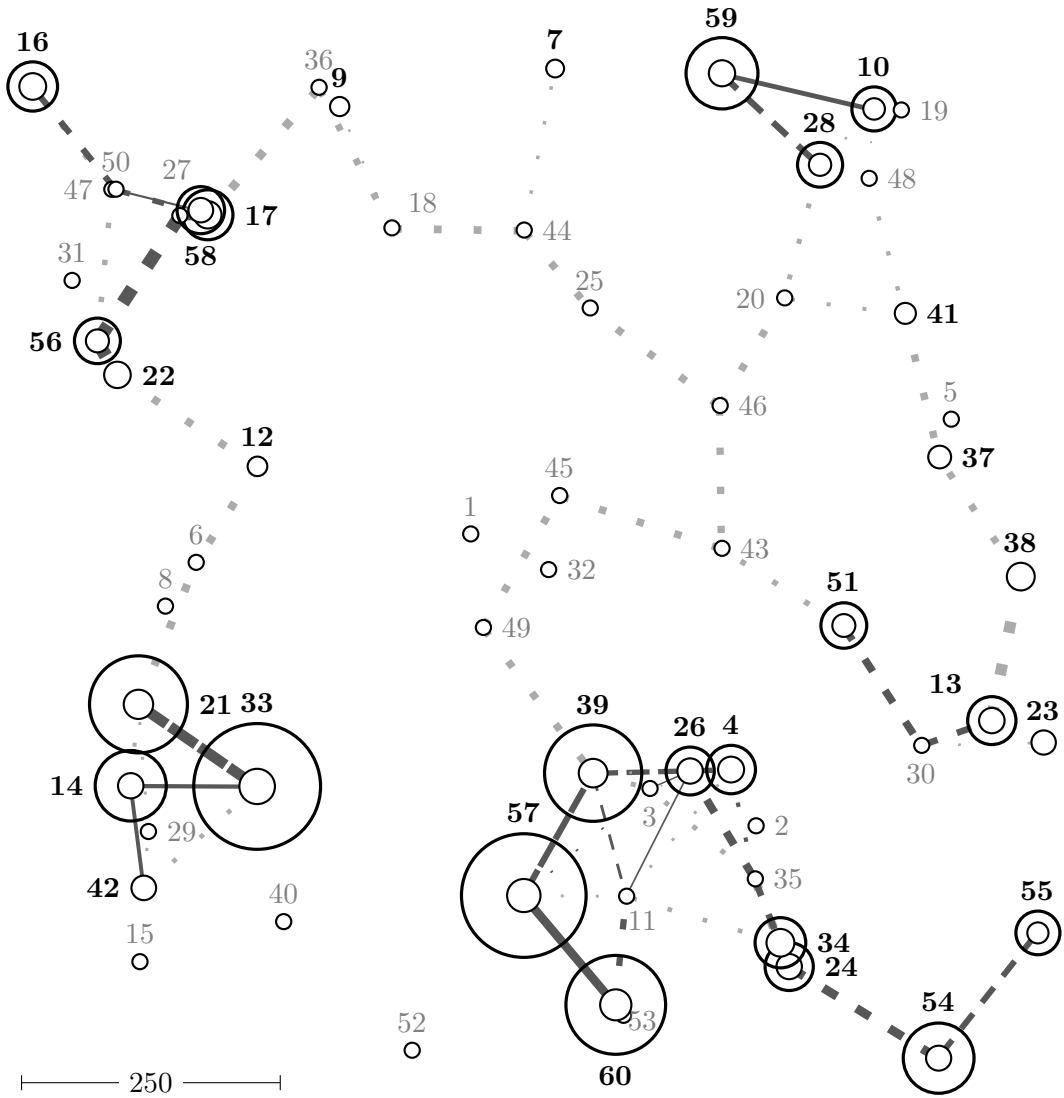

		\centering
			\begin{comment:figures}
				%
			\end{comment:figures}
		\caption{Test instance s60w30: \CFRLP\ -- $S = 36$  (25\%).}
		\label{figure:TestInstanceSSixtyWThirtyCFRLPSThirtySix}
	\end{figure}

Locating 36 charging poles (see Figure~\ref{figure:TestInstanceSSixtyWThirtyCFRLPSThirtySix}), which represents 25\% of EVCP, still covers 47.68\% of the maximum number of EVs that can be covered given the capacity limit of a charging pole. On average there are 1.7 charging poles per station. With the actual charging station placement it would be possible to cover 38 flows, but due to the capacity limitation at charging poles, it is just possible to partially cover 23 flows.

\medskip

In the following, the {\em utilisation} at charging poles is analysed in different scenarios. By pre-testing, the maximum number of EVs which can complete their round-trips successfully, given the capacity limit of charging poles, was evaluated. In the first pre-test, when allocating the maximum possible number of charging poles (S=196), four charging poles are located in 49 potential facility location and this allocation results in average utilisation of 64.49\% per charging station. The second pre-test indicates that the same coverage level can be reached in a case that 144 charging poles are allocated. Following the average utilisation can be increased to 84.24\% and therefore leads to a minimization of idle time per location. 

In summary, it can be ascertained that the number of charging poles is chosen to minimize idle time per location. They are located in order to have high capacity utilisation. Meaning, there are fewer charging poles at less frequented charging locations and the maximum number of charging poles per location is installed at highly frequented locations to maximize capacity utilisation. By increasing the number of charging poles, when the highly frequented charging locations are already equipped with the maximum number of charging poles possible, additional flows with lower flow volumes are covered as well. To cover these less frequented flows, the volume of flows passing by a charging station location can be smaller than the capacity of a single charging pole. Hence, not the entire capacity of this charging pole is utilized. As a consequence, increasing the number of charging poles to locate, leads to a reduction of the average capacity utilisation per charging station. In the case where 36 charging poles are allocated, the average charging station utilisation is 98.40\%. Comparing it with the scenario of installing 144 charging poles within the network, the average utilisation per charging facility is 84.24\%. 

\medskip


The assumption, that energy consumption is linearly proportional to the travel distance, is a simplification, which depicts the duration of the stopover of EVs. Following, EVs travelling short distances between two consecutively passed charging stations have a lower energy demand than those that travel long distances. This results in a reduced energy demand of EVs, which have refuelled at immediately preceding charging stations. Thus, it is possible to reduce energy consumption at a charging station, by locating another one in close proximity. 
In general, an increasing number of charging poles to locate leads to an increasing average number of charging poles per location. As a result, increasing proportions of flows can refuel at these locations and therefore complete their round-trips without running out of fuel.

If the number of charging poles rises, CFV accordingly increases to a smaller degree. This can be explained when looking at the proportions of the covered flows. When 36 charging poles are allocated (see Figure~\ref{figure:TestInstanceSSixtyWThirtyCFRLPSThirtySix}), it is possible to fully cover 10 out of 23 flows and 6 other flows are covered by less than 50\% ($z_f < 0.5$), which results in a CFV of \num{394588.66} (39.46\%) EVs. In case, 144 charging stations are installed (see Figure~\ref{figure:TestInstanceSSixtyWThirtyCFRLPSOneHundretAndFourtyFour}), on average 3 charging poles are built per location. In the following 176 out of 188 flows are fully covered($z_f = 1$), while for a single flow the proportional flow volume covered is less than 50\% ($z_f < 0.5$). This charging pole placement results in total coverage of \num{827524.96} (82.75\%) EVs of TFV. Furthermore, it is important to note, that the number of not fully utilized charging stations increases with an increasing number of charging poles to locate.

Appendix~\ref{Appendix:CFRLP:results} depicts the capacity limitation per charging pole as well as the results of the \CFRLP\ for all other test instances in Tables \ref{table:Results_CMCFRLP_pre_s40w20}--\ref{table:Results_CFRLP_Florida}.

    \section{Conclusions and future research}
        \label{section:ConclusionsAndFutureResearch}
Increasing pressure on reduction of environmental pollution and dependency on petroleum, leads to an increasing acceptance of electric vehicles (EV) towards traditionally fossil-fuel powered automobiles. Consequently, rapid growth in the number of EVs requires an urgent need to develop an adequate charging station infrastructure to stimulate and facilitate their usage. Restricted budget for investments in the development of a charging station infrastructure requires to choose locations deliberately.

In this paper, three extensions considering different objectives and various constraints to the deterministic flow refuelling location problem (DFRLP), described by \citet{DeVriesDuijzer:IncorporatingDrivingRangeVariabilityInNetworkDesignForRefuelingFacilities}, are introduced and implemented. Furthermore, these extensions are analysed using randomly generated problem instances from the literature. 

\medskip

The first model extension, discussed in Section~\ref{subsection:MinCoverage}, addresses the flow volume covered (CFV) as a hard constraint. In contrast, most research focusing on optimal charging station placement for EVs, consider the CFV in the objective function. The numerical analysis clearly showed that---at least in the early stages of infrastructure development---it might be insufficient to force an infrastructure that is capable of covering all EVs travelling within a network. The minimum number of charging stations needed depends, among other things, on the length of the road segments and the average travel distance between origin and destination. 

\medskip

The second model extension, introduced in Section~\ref{subsection:LocationDependentCost}, deals with location-dependent costs and returns useful information on the effects of cost differences concerning the construction costs for charging stations. This research shows that when considering location-dependent costs, results heavily depend on their relation. Tests are carried out for different cost scenarios and policy implications are discussed. 

\medskip

Whereas at the early stage of e-mobility the number of EVs was quite small, the assumption of unlimited capacity at charging stations did not pose any problems. Due to the significantly rising number of EVs, it becomes necessary to think about capacity limits and therefore the size of charging stations. The enhance model, described in Section~\ref{subsection:StationSize}, considers both, the location and size of charging stations. In comparison to existing models, capacity is defined as the quantity of energy available at a charging pole. Available energy at charging stations is not really a limiting condition, but can be taken as implicit measure for the duration of loading. 



            \paragraph{Further Research}
In the model extension \textit{\CFRLP}, dealing with limited capacity per charging pole, we allow only one sequence of opened charging stations per route. It is not possible to split the flow volume of a certain flow $f$ such that different flow proportions use different sequences of opened charging stations. I.e.\ it is not possible to skip some stations by some EVs and to use the remaining ones to refuel a larger proportion of EVs travelling along the same path, but recharging at different (not fully utilised) charging stations. According to this model formulation and the definition of the variable $i_{k l f}$, indicating a cycle segment, substantial changes are necessary to allow different combinations of charging stations within a flow $f$. This could be the objective of further research. 

\medskip

Whereas in the previous model formulation \textit{\CFRLP}, described in Section~\ref{subsection:StationSize}, it is irrelevant whether charging poles are placed separately at different locations or in groups. However, based on economic considerations it would be reasonable to assume decreasing costs per charging pole with increasing scale of charging stations. That is justified by sharing fixed costs among several charging poles.\cite[\citePrefix][p.\ 145]{UpchurchKubyLim:AModelForLocationOfCapacitatedAlternativeFuelStations}
A topic for further research is to consider non-linearly increasing fixed costs of deploying an infrastructure with an increasing number of charging poles installed per location, as a connection to the power grid is available, regardless whether one or more poles are installed. Thus, installing more than one charging pole at a location leads to an increasing capacity per monetary unit. This cost structure could be represented with a piecewise linear function which is concave, as the successive slopes are non-increasing.

\medskip

Finally, it could be also modelled that there exist different types of charging poles: faster ones, which are more expansive, and slower ones, which are cheaper.
	
	\section*{Acknowledgements}
		\label{section:Acknowledgements}
We would like to thank Klaus Ladner for his technical support.

				
	
	\bibliographystyle{IEEEtranSN}
	\bibliography{KastnerGreistorferStanek_AdvancedOptimizationModelsForTheLocationOfChargingStationsInEMobility}

\begin{thebibliography}{16}
\providecommand{\natexlab}[1]{#1}
\providecommand{\url}[1]{#1}
\csname url@samestyle\endcsname
\providecommand{\newblock}{\relax}
\providecommand{\bibinfo}[2]{#2}
\providecommand{\BIBentrySTDinterwordspacing}{\spaceskip=0pt\relax}
\providecommand{\BIBentryALTinterwordstretchfactor}{4}
\providecommand{\BIBentryALTinterwordspacing}{\spaceskip=\fontdimen2\font plus
\BIBentryALTinterwordstretchfactor\fontdimen3\font minus
  \fontdimen4\font\relax}
\providecommand{\BIBforeignlanguage}[2]{{%
\expandafter\ifx\csname l@#1\endcsname\relax
\typeout{** WARNING: IEEEtranSN.bst: No hyphenation pattern has been}%
\typeout{** loaded for the language `#1'. Using the pattern for}%
\typeout{** the default language instead.}%
\else
\language=\csname l@#1\endcsname
\fi
#2}}
\providecommand{\BIBdecl}{\relax}
\BIBdecl

\bibitem[Eur(2021)]{EuropeanEnvironmentAgency:GreenhouseGasEmissionsFromTransportInEurope}
``Greenhouse gas emissions from transport in europe,'' \\European Environment
  Agency (EEA), available online at\\
  \url{https://www.eea.europa.eu/ds_resolveuid/33a71dc1855946d8a6749d0818077e96},\\
  downloaded on 18.09.2021, 2021.

\bibitem[Capar and
  Kuby(2012)]{CaparKuby:AnEfficientFormulationOfTheFlowRefuelingLocationModelForAlternativeFuelStations}
I.~Capar and M.~Kuby, ``An efficient formulation of the flow refueling location
  model for alternative-fuel stations,'' \emph{IIE Transactions}, vol.~44,
  no.~8, pp. 622--636, 2012.

\bibitem[{de Vries} and
  Duijzer(2017)]{DeVriesDuijzer:IncorporatingDrivingRangeVariabilityInNetworkDesignForRefuelingFacilities}
H.~{de Vries} and E.~Duijzer, ``Incorporating driving range variability in
  network design for refueling facilities,'' \emph{Omega}, vol.~69, pp.
  102--114, 2017.

\bibitem[Giménez-Gaydou et~al.(2017)Giménez-Gaydou, Ribeiro, Gutiérrez, and
  Antunes]{GimenezGaydouRibeiroGutierrezAntunes:OptimalLocationOfBatteryElectricVehicleChargingStationsInUrbanAreasANewApproach}
D.~A. Giménez-Gaydou, A.~S.~N. Ribeiro, J.~Gutiérrez, and A.~P. Antunes,
  ``Optimal location of battery electric vehicle charging stations in urban
  areas: A new approach,'' \emph{International Journal of Sustainable
  Transportation}, vol.~10, no.~5, pp. 393--405, 2017.

\bibitem[Hakimi(1964)]{Hakimi:OptimumLocationsOfSwitchingCentersAndTheAbsoluteCentersAndMediansOfAGraph}
S.~L. Hakimi, ``Optimum locations of switching centers and the absolute centers
  and medians of a graph,'' \emph{Operations Research}, vol.~12, no.~3, pp.
  450--459, 1964.

\bibitem[Hodgson(1990)]{Hodgson:AFlowCapturingLocationAllocationModel}
M.~J. Hodgson, ``A flow‐capturing location‐allocation model,''
  \emph{Geographical Analysis}, vol.~22, no.~3, pp. 270--279, 1990.

\bibitem[Hosseini and
  Mirhassani(2017)]{HosseiniMirhassani:AHeuristicAlgorithmForOptimalLocationOfFlowRefuelingCapacitatedStation}
M.~Hosseini and S.~A. Mirhassani, ``A heuristic algorithm for optimal location
  of flow-refueling capacitated stations,'' \emph{International Transactions in
  Operational Research}, vol.~24, no.~6, pp. 1377--1403, 2017.

\bibitem[Kuby and
  Lim(2005)]{KubyLim:TheFlowRefuelingLocationProblemForAlternativeFuelVehicles}
M.~Kuby and S.~Lim, ``The flow-refueling location problem for alternative-fuel
  vehicles,'' \emph{Socio-Economic Planning Sciences}, vol.~39, no.~2, pp.
  125--145, 2005.

\bibitem[Lim and
  Kuby(2010)]{LimKuby:HeuristicAlgorithmsForSitingAlternativeFuelStationsUsingTheFlowRefuelingLocationModel}
S.~Lim and M.~Kuby, ``Heuristic algorithms for siting alternative-fuel stations
  using the flow-refueling location model,'' \emph{European Journal of
  Operational Research}, vol. 204, no.~1, pp. 51--61, 2010.

\bibitem[Mehar and
  Senouci(2013)]{MeharSenouci:AnOptimizationLocationSchemeForElectricChargingStations}
S.~Mehar and S.~M. Senouci, ``An optimization location scheme for electric
  charging stations,'' in \emph{2013 International Conference on Smart
  Communications in Network Technologies (SaCoNeT)}, vol.~01.\hskip 1em plus
  0.5em minus 0.4em\relax Institute of Electrical and Electronics Engineers
  (IEEE), 2013, pp. 1--5.

\bibitem[ReVelle and Swain(1970)]{ReVelleSwain:CentralFacilitiesLocation}
C.~S. ReVelle and R.~W. Swain, ``Central facilities location,''
  \emph{Geographical Analysis}, vol.~2, no.~1, pp. 30--42, 1970.

\bibitem[Toregas et~al.(1971)Toregas, Swain, ReVelle, and
  Bergman]{ToregasSwainReVelleBergman:TheLocationOfEmergencyServiceFacilities}
C.~Toregas, R.~Swain, C.~ReVelle, and L.~Bergman, ``The location of emergency
  service facilities,'' \emph{Operations Research}, vol.~19, no.~6, pp.
  1363--1373, 1971.

\bibitem[Upchurch et~al.(2009)Upchurch, Kuby, and
  Lim]{UpchurchKubyLim:AModelForLocationOfCapacitatedAlternativeFuelStations}
C.~Upchurch, M.~Kuby, and S.~Lim, ``A model for location of capacitated
  alternative-fuel stations,'' \emph{Geographical Analysis}, vol.~41, no.~1,
  pp. 85--106, 2009.

\bibitem[Wang and
  Lin(2013)]{WangLin:LocatingMultipleTypesOfRechargingStationsForBatteryPoweredElectricVehicleTransport}
Y.-W. Wang and C.-C. Lin, ``Locating multiple types of recharging stations for
  battery-powered electric vehicle transport,'' \emph{Transportation Research
  Part E: Logistics and Transportation Review}, vol.~58, pp. 76--87, 2013.

\bibitem[Wu and
  Niu(2017)]{WuNiu:StudyOnInfluenceFactorsOfElectricVehiclesChargingStationLocationBasedOnISMAndFMICMAC}
H.~Wu and D.~Niu, ``Study on influence factors of electric vehicles charging
  station location based on ism and fmicmac,'' \emph{Sustainability}, vol.~9,
  no.~4, pp. 1--19, 2017.

\bibitem[Wu and
  Zhang(2017)]{WuZhang:CanTheDevelopmentOfElectricVehiclesReduceTheEmissionOfAirPollutantsAndGreenhouseGasesInDevelopingCountries}
Y.~Wu and L.~Zhang, ``Can the development of electric vehicles reduce the
  emission of air pollutants and greenhouse gases in developing countries?''
  \emph{Transportation Research Part D: Transport and Environment}, vol.~51,
  pp. 129--145, 2017.

\end{thebibliography}
	
	\medskip
	
	\appendix
	
	\newpage
	
	\section*{Appendix}

\section{Visualization of test instances}
\label{Appendix:Visualisation}
\begin{figure}[htb!]
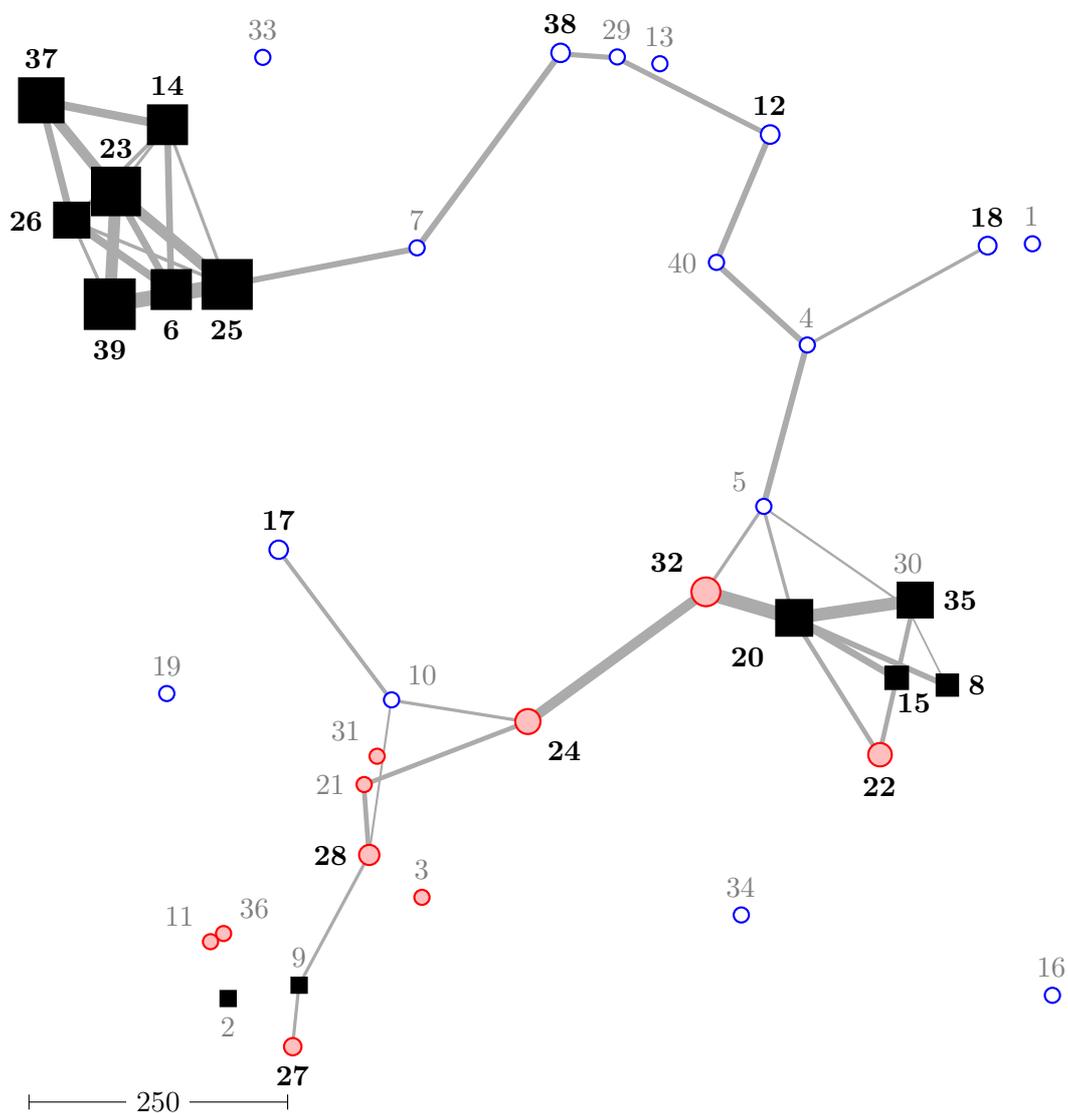

		\centering
			\begin{comment:figures}

			\end{comment:figures}
		\caption{Test instance s40w20.}
		\label{figure:TestInstanceSFourtyWTwenty}
	\end{figure}

\begin{figure}[htb!]
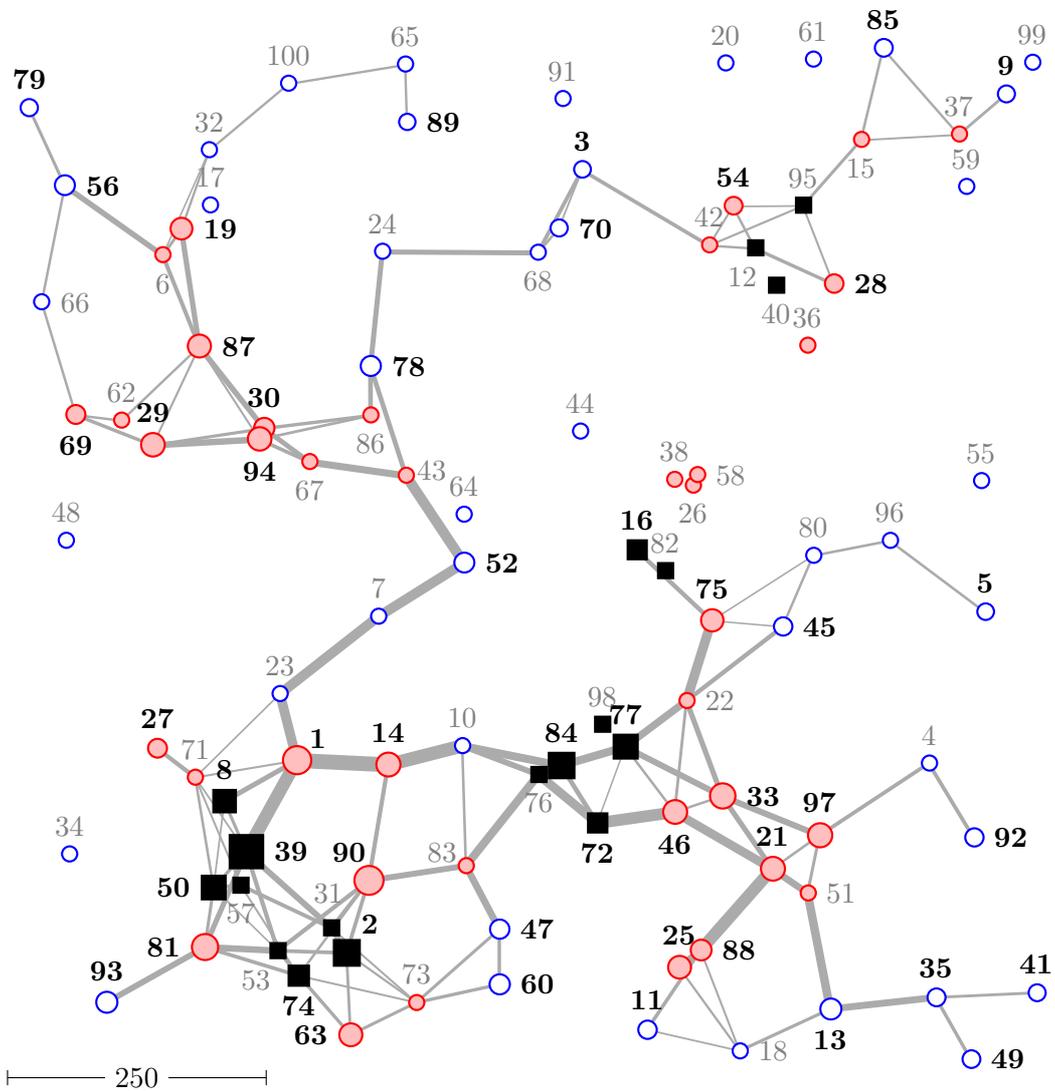

		\centering
		\begin{comment:figures}
			%
		\end{comment:figures}
		\caption{Test instance s100w50.}
		\label{figure:TestInstanceSHundredWFifty}
	\end{figure}

Figures \ref{figure:TestInstanceSFourtyWTwenty}--\ref{figure:TestInstanceSHundredWFifty} visualize the test instances s40w20, s80w40 and s100w50 considering the classification of potential facility locations to cost categories. Black rectangles represent urban potential facility locations, red points denote locations in sub-urban areas and rural potential locations are marked with blue circles.

	\begin{figure}[htb!]
		\centering
		\begin{comment:figures}
			%
		\end{comment:figures}
		\caption{Test instance s80w40.}
		\label{figure:TestInstanceSEightyWFourty}
	\end{figure}

\section{Computational results of remaining test instances}
\label{Appendix:Results}
The following section depicts the results of the test instance s40w20, s80w40 and s100w50 (if not already described), followed by a numerical analyse of the performance of the model extensions described in previous sections.

	\subsection{Minimum flow volume coverage (\MCFRLP)}
	\label{Appendix:MCFRLP}
Figures \ref{figure:numberOfOpenedFacilitiesDependentOnMinimumFlowCoverageForTestInstanceSEightywFourty}
and \ref{figure:numberOfOpenedFacilitiesDependentOnMinimumFlowCoverageForTestInstanceSHundredwFifty} depicts the relation between the number of charging stations required and the minimum flow coverage for test instance s80w40 and s100w50. 
	\begin{figure}[htb!]
		\begin{multicols}{2}
			\begin{figurehere}
				\centering 
				\newcommand*{\maxXNumberOfOpenedFacilitiesDependentOnMinimumFlowCoverageForTestInstanceSEightywFourty}{1}
				\FPeval\xScaleNumberOfOpenedFacilitiesDependentOnMinimumFlowCoverageForTestInstanceSEightywFourty{\generalXScale / 2 / \maxXNumberOfOpenedFacilitiesDependentOnMinimumFlowCoverageForTestInstanceSEightywFourty}
				\newcommand*{\maxYNumberOfOpenedFacilitiesDependentOnMinimumFlowCoverageForTestInstanceSEightywFourty}{50}
				\FPeval\yScaleNumberOfOpenedFacilitiesDependentOnMinimumFlowCoverageForTestInstanceSEightywFourty{\generalYScale / \maxYNumberOfOpenedFacilitiesDependentOnMinimumFlowCoverageForTestInstanceSEightywFourty}
				\begin{comment:figures}
					\begin{tikzpicture}[xscale=\xScaleNumberOfOpenedFacilitiesDependentOnMinimumFlowCoverageForTestInstanceSEightywFourty, yscale=\yScaleNumberOfOpenedFacilitiesDependentOnMinimumFlowCoverageForTestInstanceSEightywFourty]
						\pgfgettransformentries{\xScaleTikz}{\@tempa}{\@tempa}{\yScaleTikz}{\@tempa}{\@tempa}
						
						\draw[very thin, color=gray, xstep = 0.2, ystep = 10] (0, 0) grid (\maxXNumberOfOpenedFacilitiesDependentOnMinimumFlowCoverageForTestInstanceSEightywFourty, \maxYNumberOfOpenedFacilitiesDependentOnMinimumFlowCoverageForTestInstanceSEightywFourty);
						
						\def\crossSizeX{\crossSize / \xScaleTikz};
						\def\crossSizeY{\crossSize / \yScaleTikz};
						
						\def\crossOne{(-\crossSizeX,-\crossSizeY) -- (\crossSizeX,\crossSizeY) (-\crossSizeX,\crossSizeY) -- (\crossSizeX,-\crossSizeY)};
						
						\draw[black, shift={(0.1, 2)}] \crossOne;
						\draw[black, shift={(0.2, 3)}] \crossOne;
						\draw[black, shift={(0.3, 5)}] \crossOne;
						\draw[black, shift={(0.4, 6)}] \crossOne;
						\draw[black, shift={(0.5, 8)}] \crossOne;
						\draw[black, shift={(0.6, 11)}] \crossOne;
						\draw[black, shift={(0.7, 15)}] \crossOne;
						\draw[black, shift={(0.8, 20)}] \crossOne;
						\draw[black, shift={(0.9, 25)}] \crossOne;
						\draw[black, shift={(1.0, 32)}] \crossOne;
						
						\draw[black!75, dotted] (0.1, 2) -- (0.2, 3) -- (0.3, 5) -- (0.4, 6) -- (0.5, 8) -- (0.6, 11) -- (0.7, 15) -- (0.8, 20) -- (0.9, 25) -- (1.0, 32);
						
						\def\axisAdditionalLengthPlusTikzX{\axisAdditionalLengthPlus / \xScaleTikz}
						\def\axisAdditionalLengthMinusTikzX{\axisAdditionalLengthMinus / \xScaleTikz}
						\draw[arrow] (-\axisAdditionalLengthMinusTikzX, 0) -- (\maxXNumberOfOpenedFacilitiesDependentOnMinimumFlowCoverageForTestInstanceSEightywFourty, 0) -- +(\axisAdditionalLengthPlusTikzX, 0) node[right] {$C$};
						
						\def\axisLabelTikzY{\axisLabel / \yScaleTikz}
						\foreach \pos in {0.0, 0.2, 0.4, 0.6, 0.8, 1.0} \draw[shift={(\pos, 0)}] (0, \axisLabelTikzY) -- (0, -\axisLabelTikzY) node[below] {$\pos$};
						
						\def\axisAdditionalLengthPlusTikzY{\axisAdditionalLengthPlus / \yScaleTikz}
						\def\axisAdditionalLengthMinusTikzY{\axisAdditionalLengthMinus / \yScaleTikz}
						\draw[arrow] (0, -\axisAdditionalLengthMinusTikzY) -- (0, \maxYNumberOfOpenedFacilitiesDependentOnMinimumFlowCoverageForTestInstanceSEightywFourty) -- +(0, \axisAdditionalLengthPlusTikzY) node[above] {$p$};
						
						\def\axisLabelTikzX{\axisLabel / \xScaleTikz}
						\foreach \pos in {0, 10, 20, 30, 40, 50} \draw[shift={(0, \pos)}] (\axisLabelTikzX, 0) -- (-\axisLabelTikzX, 0) node[left] {$\pos$};
						
					\end{tikzpicture}
				\end{comment:figures}
				\caption{Number of opened facilities dependent on minimum flow coverage for test instance: s80w40.}
				\label{figure:numberOfOpenedFacilitiesDependentOnMinimumFlowCoverageForTestInstanceSEightywFourty}
			\end{figurehere}

			\begin{figurehere}
				\centering 

				\newcommand*{\maxXNumberOfOpenedFacilitiesDependentOnMinimumFlowCoverageForTestInstanceSHundredwFifty}{1}
				\FPeval\xScaleNumberOfOpenedFacilitiesDependentOnMinimumFlowCoverageForTestInstanceSHundredwFifty{\generalXScale / 2 / \maxXNumberOfOpenedFacilitiesDependentOnMinimumFlowCoverageForTestInstanceSHundredwFifty}
				\newcommand*{\maxYNumberOfOpenedFacilitiesDependentOnMinimumFlowCoverageForTestInstanceSHundredwFifty}{50}
				\FPeval\yScaleNumberOfOpenedFacilitiesDependentOnMinimumFlowCoverageForTestInstanceSHundredwFifty{\generalYScale / \maxYNumberOfOpenedFacilitiesDependentOnMinimumFlowCoverageForTestInstanceSHundredwFifty}
				\begin{comment:figures}
					\begin{tikzpicture}[xscale=\xScaleNumberOfOpenedFacilitiesDependentOnMinimumFlowCoverageForTestInstanceSHundredwFifty, yscale=\yScaleNumberOfOpenedFacilitiesDependentOnMinimumFlowCoverageForTestInstanceSHundredwFifty]
						\pgfgettransformentries{\xScaleTikz}{\@tempa}{\@tempa}{\yScaleTikz}{\@tempa}{\@tempa}
						
						\draw[very thin, color=gray, xstep = 0.2, ystep = 10] (0, 0) grid (\maxXNumberOfOpenedFacilitiesDependentOnMinimumFlowCoverageForTestInstanceSHundredwFifty, \maxYNumberOfOpenedFacilitiesDependentOnMinimumFlowCoverageForTestInstanceSHundredwFifty);
						
						\def\crossSizeX{\crossSize / \xScaleTikz};
						\def\crossSizeY{\crossSize / \yScaleTikz};
						
						\def\crossOne{(-\crossSizeX,-\crossSizeY) -- (\crossSizeX,\crossSizeY) (-\crossSizeX,\crossSizeY) -- (\crossSizeX,-\crossSizeY)};
						
						\draw[black, shift={(0.1, 2)}] \crossOne;
						\draw[black, shift={(0.2, 4)}] \crossOne;
						\draw[black, shift={(0.3, 7)}] \crossOne;
						\draw[black, shift={(0.4, 10)}] \crossOne;
						\draw[black, shift={(0.5, 13)}] \crossOne;
						\draw[black, shift={(0.6, 16)}] \crossOne;
						\draw[black, shift={(0.7, 20)}] \crossOne;
						\draw[black, shift={(0.8, 25)}] \crossOne;
						\draw[black, shift={(0.9, 31)}] \crossOne;
						\draw[black, shift={(1.0, 45)}] \crossOne;
						
						\draw[black!75, dotted] (0.1, 2) -- (0.2, 4) -- (0.3, 7) -- (0.4, 10) -- (0.5, 13) -- (0.6, 16) -- (0.7, 20) -- (0.8, 25) -- (0.9, 31) -- (1.0, 45);
						
						\def\axisAdditionalLengthPlusTikzX{\axisAdditionalLengthPlus / \xScaleTikz}
						\def\axisAdditionalLengthMinusTikzX{\axisAdditionalLengthMinus / \xScaleTikz}
						\draw[arrow] (-\axisAdditionalLengthMinusTikzX, 0) -- (\maxXNumberOfOpenedFacilitiesDependentOnMinimumFlowCoverageForTestInstanceSHundredwFifty, 0) -- +(\axisAdditionalLengthPlusTikzX, 0) node[right] {$C$};
						
						\def\axisLabelTikzY{\axisLabel / \yScaleTikz}
						\foreach \pos in {0.0, 0.2, 0.4, 0.6, 0.8, 1.0} \draw[shift={(\pos, 0)}] (0, \axisLabelTikzY) -- (0, -\axisLabelTikzY) node[below] {$\pos$};
						
						\def\axisAdditionalLengthPlusTikzY{\axisAdditionalLengthPlus / \yScaleTikz}
						\def\axisAdditionalLengthMinusTikzY{\axisAdditionalLengthMinus / \yScaleTikz}
						\draw[arrow] (0, -\axisAdditionalLengthMinusTikzY) -- (0, \maxYNumberOfOpenedFacilitiesDependentOnMinimumFlowCoverageForTestInstanceSHundredwFifty) -- +(0, \axisAdditionalLengthPlusTikzY) node[above] {$p$};
						
						\def\axisLabelTikzX{\axisLabel / \xScaleTikz}
						\foreach \pos in {0, 10, 20, 30, 40, 50} \draw[shift={(0, \pos)}] (\axisLabelTikzX, 0) -- (-\axisLabelTikzX, 0) node[left] {$\pos$};
						
					\end{tikzpicture}
				\end{comment:figures}
				\caption{Number of opened facilities dependent on minimum flow coverage for test instance: s100w50.}
				\label{figure:numberOfOpenedFacilitiesDependentOnMinimumFlowCoverageForTestInstanceSHundredwFifty}
			\end{figurehere}
		\end{multicols}
	\end{figure}

	\subsection{Location-dependent costs per charging station (\LCFRLP)}
		\label{Appendix:LCFRLP:results}
		
		\paragraph{Analysing s40w20.}
\addtolength{\tabcolsep}{-1.5pt}
\begin{table}[htb]
\caption{Results of \LCFRLP\ for all scenarios of s40w20.} 
	\begin{center}
			\begin{tabular}[c]{c *{10}{| r} }
				\label{table:ResultsLCFRLPsFourtywTwenty}
Scenario  &	\multicolumn{2}{c|}{1}	&	\multicolumn{2}{c|}{2}	&	\multicolumn{2}{c|}{3}	&	\multicolumn{2}{c|}{4}	&	\multicolumn{2}{c}{5}\\
                \hline
$CFV$  &	\multicolumn{2}{r|}{\num{656799}}	&	\multicolumn{2}{r|}{\num{708043}}	&	\multicolumn{2}{r|}{\num{722931}}	&	\multicolumn{2}{r|}{\num{730592}}	&	\multicolumn{2}{r}{\num{759839}}\\

$\sum_{k \in K} x_k$   &	\multicolumn{2}{r|}{8}	&	\multicolumn{2}{r|}{6}	&	\multicolumn{2}{r|}{5}	&	\multicolumn{2}{r|}{7}	&	\multicolumn{2}{r}{8}\\

				\hline
				\hline
urban (\sfrac{\% }{100} | \#)	& 	0.6667	& 4 & 	1	& 5 & 	0.7143	& 5 & 	0.625 & 5	& 	0.4167 & 5\\
sub-urban (\sfrac{\% }{100} | \#) 	& 	0.3333 & 2	& 	0 & 0	& 	0 & 0 & 	0	& 0 & 	0 & 0\\
rural (\sfrac{\% }{100} | \#)	& 	0 & 0	& 	0 & 0	& 	0.2857 & 2 	& 	0.375 & 3 	& 	0.5833	& 7\\
				\hline
				\multicolumn{11}{c}{}\\[-0.5em]
Scenario  &	\multicolumn{2}{c|}{6}	&	\multicolumn{2}{c|}{7}	&	\multicolumn{2}{c|}{8}	&	\multicolumn{2}{c|}{9}	&	\multicolumn{2}{c}{10}\\
                \hline
$CFV$  &	\multicolumn{2}{r|}{\num{723842}}	&	\multicolumn{2}{r|}{\num{730592}}	&	\multicolumn{2}{r|}{\num{752703}}	&	\multicolumn{2}{r|}{\num{767591}}	&	\multicolumn{2}{r}{\num{752703}}\\

$\sum_{k \in K} x_k$   &	\multicolumn{2}{r|}{12}	&	\multicolumn{2}{r|}{6}	&	\multicolumn{2}{r|}{8}	&	\multicolumn{2}{r|}{6}	&	\multicolumn{2}{r}{8}\\

				\hline
				\hline
urban (\sfrac{\% }{100} | \#)	& 	0.8333	& 5 & 	0.625	& 5 & 	1 & 6 	& 	0.75	 & 6 & 	1	& 6\\
sub-urban (\sfrac{\% }{100} | \#) 	& 	0.1667 & 1 	& 	0	& 0 & 	0	& 0 & 	0	& 0 & 	0	 & 0\\
rural (\sfrac{\% }{100} | \#)	& 	0	& 0 & 	0.375 & 3 	& 	0	& 0 & 	0.25	& 2 & 	0 & 0\\
				\hline
				\multicolumn{11}{c}{}\\[-0.5em]
Scenario  &	\multicolumn{2}{c|}{11}	&	\multicolumn{2}{c|}{12}&	\multicolumn{2}{c|}{13}	&	\multicolumn{2}{c|}{14}	&	\multicolumn{2}{c}{15}\\
                \hline
$CFV$  &	\multicolumn{2}{r|}{\num{768502}}	&	\multicolumn{2}{r|}{\num{788471}}&	\multicolumn{2}{r|}{\num{793966}}	&	\multicolumn{2}{r|}{\num{804270}}	&	\multicolumn{2}{r}{\num{804270}} \\

$\sum_{k \in K} x_k$   &	\multicolumn{2}{r|}{8}	&	\multicolumn{2}{r|}{7}&	\multicolumn{2}{r|}{7}	&	\multicolumn{2}{r|}{8}	&	\multicolumn{2}{r}{8} \\

				\hline
				\hline
urban (\sfrac{\% }{100} | \#)	& 	0.8571 & 6	& 	1 & 7 	& 	0.75 & 6	& 	0.875 & 7 & 	0.875 & 7 	 \\
sub-urban (\sfrac{\% }{100} | \#) 	& 	0.1429	& 1 & 	0 & 0 	& 	0.25 & 2	& 	0.125 & 1  & 	0.125 & 1 \\
rural (\sfrac{\% }{100} | \#)	& 	0	 & 0 & 	0 & 0	& 	0	 & 0 & 	0 & 0 	& 	0 & 0 \\
				\hline
			\end{tabular}
	\end{center}
\end{table}
\addtolength{\tabcolsep}{+1.5pt}

			\begin{figure}[htb]
				\centering 
				\newcommand*{\maxXProportionOfUrbanSuburbanRuralStatinsOpenedInDifferentScenariosForTheTestInstanceSFourtyWTwenty}{15}
				\FPeval\xScaleProportionOfUrbanSuburbanRuralStatinsOpenedInDifferentScenariosForTheTestInstanceSFourtyWTwenty{\generalXScale / 0.838 / \maxXProportionOfUrbanSuburbanRuralStatinsOpenedInDifferentScenariosForTheTestInstanceSFourtyWTwenty}
				\newcommand*{\maxYProportionOfUrbanSuburbanRuralStatinsOpenedInDifferentScenariosForTheTestInstanceSFourtyWTwentyNumber}{9.0}
				\FPeval\yScaleProportionOfUrbanSuburbanRuralStatinsOpenedInDifferentScenariosForTheTestInstanceSFourtyWTwentyNumber{\generalYScale / \maxYProportionOfUrbanSuburbanRuralStatinsOpenedInDifferentScenariosForTheTestInstanceSFourtyWTwentyNumber}
				
\newcommand*{\maxYProportionOfUrbanSuburbanRuralStatinsOpenedInDifferentScenariosForTheTestInstanceSFourtyWTwentyProportion}{1.0}
				\FPeval\yScaleProportionOfUrbanSuburbanRuralStatinsOpenedInDifferentScenariosForTheTestInstanceSFourtyWTwentyProportion{\generalYScale / \maxYProportionOfUrbanSuburbanRuralStatinsOpenedInDifferentScenariosForTheTestInstanceSFourtyWTwentyProportion}
				\begin{comment:figures}
					\begin{tikzpicture}[xscale=\xScaleProportionOfUrbanSuburbanRuralStatinsOpenedInDifferentScenariosForTheTestInstanceSFourtyWTwenty, yscale=\yScaleProportionOfUrbanSuburbanRuralStatinsOpenedInDifferentScenariosForTheTestInstanceSFourtyWTwentyNumber]
						\pgfgettransformentries{\xScaleTikz}{\@tempa}{\@tempa}{\yScaleTikz}{\@tempa}{\@tempa}
						
						\draw[very thin, color=gray, xstep = 1, ystep = 1] (0, 0) grid (\maxXProportionOfUrbanSuburbanRuralStatinsOpenedInDifferentScenariosForTheTestInstanceSFourtyWTwenty, \maxYProportionOfUrbanSuburbanRuralStatinsOpenedInDifferentScenariosForTheTestInstanceSFourtyWTwentyNumber);
						
						\def\axisLabelTikzY{\axisLabel / \yScaleTikz}
						\def\axisNumberTikzY{5.5 * \axisLabelTikzY}
						\def\axisTotalNumberTikzY{3 * \axisLabelTikzY}
						\def\axisLabelTikzX{\axisLabel / \xScaleTikz}
						\draw[thick, color=gray] (5.0, -\axisNumberTikzY) -- (5.0, 9.0 + \axisNumberTikzY);
						\draw[thick, color=gray] (9.0, -\axisNumberTikzY) -- (9.0, 9.0 + \axisNumberTikzY);
						\draw[thick, color=gray] (12.0, -\axisNumberTikzY) -- (12.0, 9.0 + \axisNumberTikzY);
						\draw[thick, color=gray] (14.0, -\axisNumberTikzY) -- (14.0, 9.0 + \axisNumberTikzY);
						
						\begin{scope}[shift={(-0.5, 0)}]
							\draw[black, fill=black] (0.6, 0.0) rectangle (0.8, 4);
							\draw[red, thick, fill=red!25] (0.9, 0.0) rectangle (1.1, 2);
							\draw[blue, thick] (1.2, 0.0) rectangle (1.4, 0);
							
							
							\draw[black, fill=black] (1.6, 0.0) rectangle (1.8, 5);
							\draw[red, thick, fill=red!25] (1.9, 0.0) rectangle (2.1, 0);
							\draw[blue, thick] (2.2, 0.0) rectangle (2.4, 0);
							
							\draw[black, fill=black] (2.6, 0.0) rectangle (2.8, 5);
							\draw[red, thick, fill=red!25] (2.9, 0.0) rectangle (3.1, 0);
							\draw[blue, thick] (3.2, 0.0) rectangle (3.4, 2);
							
							\draw[black, fill=black] (3.6, 0.0) rectangle (3.8, 5);
							\draw[red, thick, fill=red!25] (3.9, 0.0) rectangle (4.1, 0);
							\draw[blue, thick] (4.2, 0.0) rectangle (4.4, 3);
							
							\draw[black, fill=black] (4.6, 0.0) rectangle (4.8, 5);
							\draw[red, thick, fill=red!25] (4.9, 0.0) rectangle (5.1, 0);
							\draw[blue, thick] (5.2, 0.0) rectangle (5.4, 7);
							
							\draw[black, fill=black] (5.6, 0.0) rectangle (5.8, 5);
							\draw[red, thick, fill=red!25] (5.9, 0.0) rectangle (6.1, 1);
							\draw[blue, thick] (6.2, 0.0) rectangle (6.4, 0);
							
							\draw[black, fill=black] (6.6, 0.0) rectangle (6.8, 5);
							\draw[red, thick, fill=red!25] (6.9, 0.0) rectangle (7.1, 0);
							\draw[blue, thick] (7.2, 0.0) rectangle (7.4, 3);
							
							\draw[black, fill=black] (7.6, 0.0) rectangle (7.8, 6);
							\draw[red, thick, fill=red!25] (7.9, 0.0) rectangle (8.1, 0);
							\draw[blue, thick] (8.2, 0.0) rectangle (8.4, 0);
							
							\draw[black, fill=black] (8.6, 0.0) rectangle (8.8, 6);
							\draw[red, thick, fill=red!25] (8.9, 0.0) rectangle (9.1, 0);
							\draw[blue, thick] (9.2, 0.0) rectangle (9.4, 2);
							
							\draw[black, fill=black] (9.6, 0.0) rectangle (9.8, 6);
							\draw[red, thick, fill=red!25] (9.9, 0.0) rectangle (10.1, 0);
							\draw[blue, thick] (10.2, 0.0) rectangle (10.4, 0);
							
							\draw[black, fill=black] (10.6, 0.0) rectangle (10.8, 6);
							\draw[red, thick, fill=red!25] (10.9, 0.0) rectangle (11.1, 1);
							\draw[blue, thick] (11.2, 0.0) rectangle (11.4, 0);
							
							\draw[black, fill=black] (11.6, 0.0) rectangle (11.8, 7);
							\draw[red, thick, fill=red!25] (11.9, 0) rectangle (12.1, 0);
							\draw[blue, thick] (12.2, 0.0) rectangle (12.4, 0);
							
							\draw[black, fill=black] (12.6, 0.0) rectangle (12.8, 6);
							\draw[red, thick, fill=red!25] (12.9, 0.0) rectangle (13.1, 2);
							\draw[blue, thick] (13.2, 0.0) rectangle (13.4, 0);
							
							\draw[black, fill=black] (13.6, 0.0) rectangle (13.8, 7);
							\draw[red, thick, fill=red!25] (13.9, 0) rectangle (14.1, 1);
							\draw[blue, thick] (14.2, 0.0) rectangle (14.4, 0);
							
							\draw[black, fill=black] (14.6, 0.0) rectangle (14.8, 7);
							\draw[red, thick, fill=red!25] (14.9, 0.0) rectangle (15.1, 1);
							\draw[blue, thick] (15.2, 0.0) rectangle (15.4, 0);
						\end{scope}
						
						\def\axisAdditionalLengthPlusTikzX{\axisAdditionalLengthPlus / \xScaleTikz}
						\def\axisAdditionalLengthMinusTikzX{\axisAdditionalLengthMinus / \xScaleTikz}
						\draw[-] (-\axisAdditionalLengthMinusTikzX, 0) -- (\maxXProportionOfUrbanSuburbanRuralStatinsOpenedInDifferentScenariosForTheTestInstanceSFourtyWTwenty, 0) -- +(0, 0) node[right] {Scenario};
						
						\foreach \pos in {1, 2, 3, 4, 5, 6, 7, 8, 9, 10, 11, 12, 13, 14, 15} \draw[shift={(\pos-0.5, 0)}] (0, -\axisLabelTikzY) node[below] {$\pos$};
						
						\def\axisAdditionalLengthPlusTikzY{\axisAdditionalLengthPlus / \yScaleTikz}
						\def\axisAdditionalLengthMinusTikzY{\axisAdditionalLengthMinus / \yScaleTikz}
						\draw[arrow] (0, -\axisAdditionalLengthMinusTikzY) -- (0, \maxYProportionOfUrbanSuburbanRuralStatinsOpenedInDifferentScenariosForTheTestInstanceSFourtyWTwentyNumber) -- +(0, \axisAdditionalLengthPlusTikzY) node[above] {Total number};
						
						\foreach \pos in {1, 2, 3, 4, 5, 6, 7, 8, 9} \draw[shift={(0, \pos)}] (\axisLabelTikzX, 0) -- (-\axisLabelTikzX, 0) node[left] {$\pos$};
						
					\end{tikzpicture}
					
					
					\hspace*{0.118cm}
					\begin{tikzpicture}[xscale=\xScaleProportionOfUrbanSuburbanRuralStatinsOpenedInDifferentScenariosForTheTestInstanceSFourtyWTwenty, yscale=\yScaleProportionOfUrbanSuburbanRuralStatinsOpenedInDifferentScenariosForTheTestInstanceSFourtyWTwentyProportion]
						\pgfgettransformentries{\xScaleTikz}{\@tempa}{\@tempa}{\yScaleTikz}{\@tempa}{\@tempa}
						
						\draw[very thin, color=gray, xstep = 1, ystep = 0.2] (0, 0) grid (\maxXProportionOfUrbanSuburbanRuralStatinsOpenedInDifferentScenariosForTheTestInstanceSFourtyWTwenty, \maxYProportionOfUrbanSuburbanRuralStatinsOpenedInDifferentScenariosForTheTestInstanceSFourtyWTwentyProportion);
						
						\def\axisLabelTikzY{\axisLabel / \yScaleTikz}
						\def\axisPartitionTikzY{5.5 * \axisLabelTikzY}
						\def\axisLabelTikzX{\axisLabel / \xScaleTikz}
						\draw[thick, color=gray] (5.0, -\axisPartitionTikzY) -- (5.0, 1.0 + \axisPartitionTikzY);
						\draw[thick, color=gray] (9.0, -\axisPartitionTikzY) -- (9.0, 1.0 + \axisPartitionTikzY);
						\draw[thick, color=gray] (12.0, -\axisPartitionTikzY) -- (12.0, 1.0 + \axisPartitionTikzY);
						\draw[thick, color=gray] (14.0, -\axisPartitionTikzY) -- (14.0, 1.0 + \axisPartitionTikzY);
						
						\begin{scope}[shift={(-0.5, 0)}]
							\draw[black, fill=black] (0.6, 0.0) rectangle (0.8, 0.666666666666667);
							\draw[red, thick, fill=red!25] (0.9, 0.0) rectangle (1.1, 0.333333333333333);
							\draw[blue, thick] (1.2, 0.0) rectangle (1.4, 0.0);
							
							\draw[black, fill=black] (1.6, 0.0) rectangle (1.8, 1.0);
							\draw[red, thick, fill=red!25] (1.9, 0.0) rectangle (2.1, 0.0);
							\draw[blue, thick] (2.2, 0.0) rectangle (2.4, 0.0);
							
							\draw[black, fill=black] (2.6, 0.0) rectangle (2.8, 0.714285714285714);
							\draw[red, thick, fill=red!25] (2.9, 0.0) rectangle (3.1, 0.0);
							\draw[blue, thick] (3.2, 0.0) rectangle (3.4, 0.285714285714286);
							
							\draw[black, fill=black] (3.6, 0.0) rectangle (3.8, 0.625);
							\draw[red, thick, fill=red!25] (3.9, 0.0) rectangle (4.1, 0.0);
							\draw[blue, thick] (4.2, 0.0) rectangle (4.4, 0.375);
							
							\draw[black, fill=black] (4.6, 0.0) rectangle (4.8, 0.416666666666667);
							\draw[red, thick, fill=red!25] (4.9, 0.0) rectangle (5.1, 0.0);
							\draw[blue, thick] (5.2, 0.0) rectangle (5.4, 0.583333333333333);
							
							\draw[black, fill=black] (5.6, 0.0) rectangle (5.8, 0.833333333333333);
							\draw[red, thick, fill=red!25] (5.9, 0.0) rectangle (6.1, 0.166666666666667);
							\draw[blue, thick] (6.2, 0.0) rectangle (6.4, 0.0);
							
							\draw[black, fill=black] (6.6, 0.0) rectangle (6.8, 0.625);
							\draw[red, thick, fill=red!25] (6.9, 0.0) rectangle (7.1, 0.0);
							\draw[blue, thick] (7.2, 0.0) rectangle (7.4, 0.375);
							
							\draw[black, fill=black] (7.6, 0.0) rectangle (7.8, 1.0);
							\draw[red, thick, fill=red!25] (7.9, 0.0) rectangle (8.1, 0.0);
							\draw[blue, thick] (8.2, 0.0) rectangle (8.4, 0.0);
							
							\draw[black, fill=black] (8.6, 0.0) rectangle (8.8, 0.75);
							\draw[red, thick, fill=red!25] (8.9, 0.0) rectangle (9.1, 0.0);
							\draw[blue, thick] (9.2, 0.0) rectangle (9.4, 0.25);
							
							\draw[black, fill=black] (9.6, 0.0) rectangle (9.8, 1.0);
							\draw[red, thick, fill=red!25] (9.9, 0.0) rectangle (10.1, 0.0);
							\draw[blue, thick] (10.2, 0.0) rectangle (10.4, 0.0);
							
							\draw[black, fill=black] (10.6, 0.0) rectangle (10.8, 0.857142857142857);
							\draw[red, thick, fill=red!25] (10.9, 0.142857142857143) rectangle (11.1, 0.0);
							\draw[blue, thick] (11.2, 0.0) rectangle (11.4, 0.0);
							
							\draw[black, fill=black] (11.6, 0.0) rectangle (11.8, 1.0);
							\draw[red, thick, fill=red!25] (11.9, 0.0) rectangle (12.1, 0.0);
							\draw[blue, thick] (12.2, 0.0) rectangle (12.4, 0.0);
							
							\draw[black, fill=black] (12.6, 0.0) rectangle (12.8, 0.75);
							\draw[red, thick, fill=red!25] (12.9, 0.25) rectangle (13.1, 0.0);
							\draw[blue, thick] (13.2, 0.0) rectangle (13.4, 0.0);
							
							\draw[black, fill=black] (13.6, 0.0) rectangle (13.8, 0.875);
							\draw[red, thick, fill=red!25] (13.9, 0.125) rectangle (14.1, 0.0);
							\draw[blue, thick] (14.2, 0.0) rectangle (14.4, 0.0);
							
							\draw[black, fill=black] (14.6, 0.0) rectangle (14.8, 0.875);
							\draw[red, thick, fill=red!25] (14.9, 0.0) rectangle (15.1, 0.125);
							\draw[blue, thick] (15.2, 0.0) rectangle (15.4, 0.0);
						\end{scope}
						
						\def\axisAdditionalLengthPlusTikzX{\axisAdditionalLengthPlus / \xScaleTikz}
						\def\axisAdditionalLengthMinusTikzX{\axisAdditionalLengthMinus / \xScaleTikz}
						\draw[-] (-\axisAdditionalLengthMinusTikzX, 0) -- (\maxXProportionOfUrbanSuburbanRuralStatinsOpenedInDifferentScenariosForTheTestInstanceSFourtyWTwenty, 0) -- +(0, 0) node[right] {Scenario};
						
						\foreach \pos in {1, 2, 3, 4, 5, 6, 7, 8, 9, 10, 11, 12, 13, 14, 15} \draw[shift={(\pos-0.5, 0)}] (0, -\axisLabelTikzY) node[below] {$\pos$};
						
						\def\axisAdditionalLengthPlusTikzY{\axisAdditionalLengthPlus / \yScaleTikz}
						\def\axisAdditionalLengthMinusTikzY{\axisAdditionalLengthMinus / \yScaleTikz}
						\draw[arrow] (0, -\axisAdditionalLengthMinusTikzY) -- (0, \maxYProportionOfUrbanSuburbanRuralStatinsOpenedInDifferentScenariosForTheTestInstanceSFourtyWTwentyProportion) -- +(0, \axisAdditionalLengthPlusTikzY) node[above] {Proportion};
						
						\foreach \pos in {0.0, 0.2, 0.4, 0.6, 0.8, 1.0} \draw[shift={(0, \pos)}] (\axisLabelTikzX, 0) -- (-\axisLabelTikzX, 0) node[left] {$\pos$};
						
					\end{tikzpicture}
				\end{comment:figures}
				\caption[Total number and proportion of urban/suburban/rural statins opened in different scenarios for the test instance s40w20.]{Total number and proportion of urban/suburban/rural statins opened in different scenarios for the test instance s40w20. The left bar \begin{tikzpicture}[xscale=\xScaleProportionOfUrbanSuburbanRuralStatinsOpenedInDifferentScenariosForTheTestInstanceSFourtyWTwenty, yscale=1]\draw[black, fill=black] (0.0, 0.0) rectangle (0.2, 0.24);\end{tikzpicture} stays for urban, the middle bar \begin{tikzpicture}[xscale=\xScaleProportionOfUrbanSuburbanRuralStatinsOpenedInDifferentScenariosForTheTestInstanceSFourtyWTwenty, yscale=1]\draw[red, thick, fill=red!25] (0.0, 0.0) rectangle (0.2, 0.24);\end{tikzpicture} for suburban, and the right bar \begin{tikzpicture}[xscale=\xScaleProportionOfUrbanSuburbanRuralStatinsOpenedInDifferentScenariosForTheTestInstanceSFourtyWTwenty, yscale=1]\draw[blue, thick] (0.0, 0.0) rectangle (0.2, 0.24);\end{tikzpicture} for rural areas.}
				\label{figure:totalNumberAndProportionOfUrbanSuburbanRuralStatinsOpenedInDifferentScenariosForTheTestInstanceSFourtyWTwenty}
			\end{figure}

Taking a closer look at the results (see Table~\ref{table:ResultsLCFRLPsFourtywTwenty}  and Figure~\ref{figure:totalNumberAndProportionOfUrbanSuburbanRuralStatinsOpenedInDifferentScenariosForTheTestInstanceSFourtyWTwenty}) and the visualization of instance s40w20 shown in Figure~\ref{figure:TestInstanceSFourtyWTwenty}, it becomes obvious that charging stations are located primarily in urban nodes. Per definition urban nodes are characterized by strongly frequented flows starting and ending their round-trip at these nodes and/or the proximity to other nodes. Moreover, in the test instances taken from \citet{DeVriesDuijzer:IncorporatingDrivingRangeVariabilityInNetworkDesignForRefuelingFacilities} the flow volume is equal to zero for flows with total travel distance between OD-nodes less than 100 (see Section~\ref{subsubsection:BenchmarkInstances}). Checking the percentage of flows that are passing by urban nodes only and summing up their flow volumes, these flows represent 73.456\% of \num{e06} EVs, denoting the TFV in instance s40w20. This is different compared to the other testing instances and therefore explains the large proportion of urban CSs in Figure~\ref{figure:totalNumberAndProportionOfUrbanSuburbanRuralStatinsOpenedInDifferentScenariosForTheTestInstanceSFourtyWTwenty}.  

\begin{figure}[htb!]
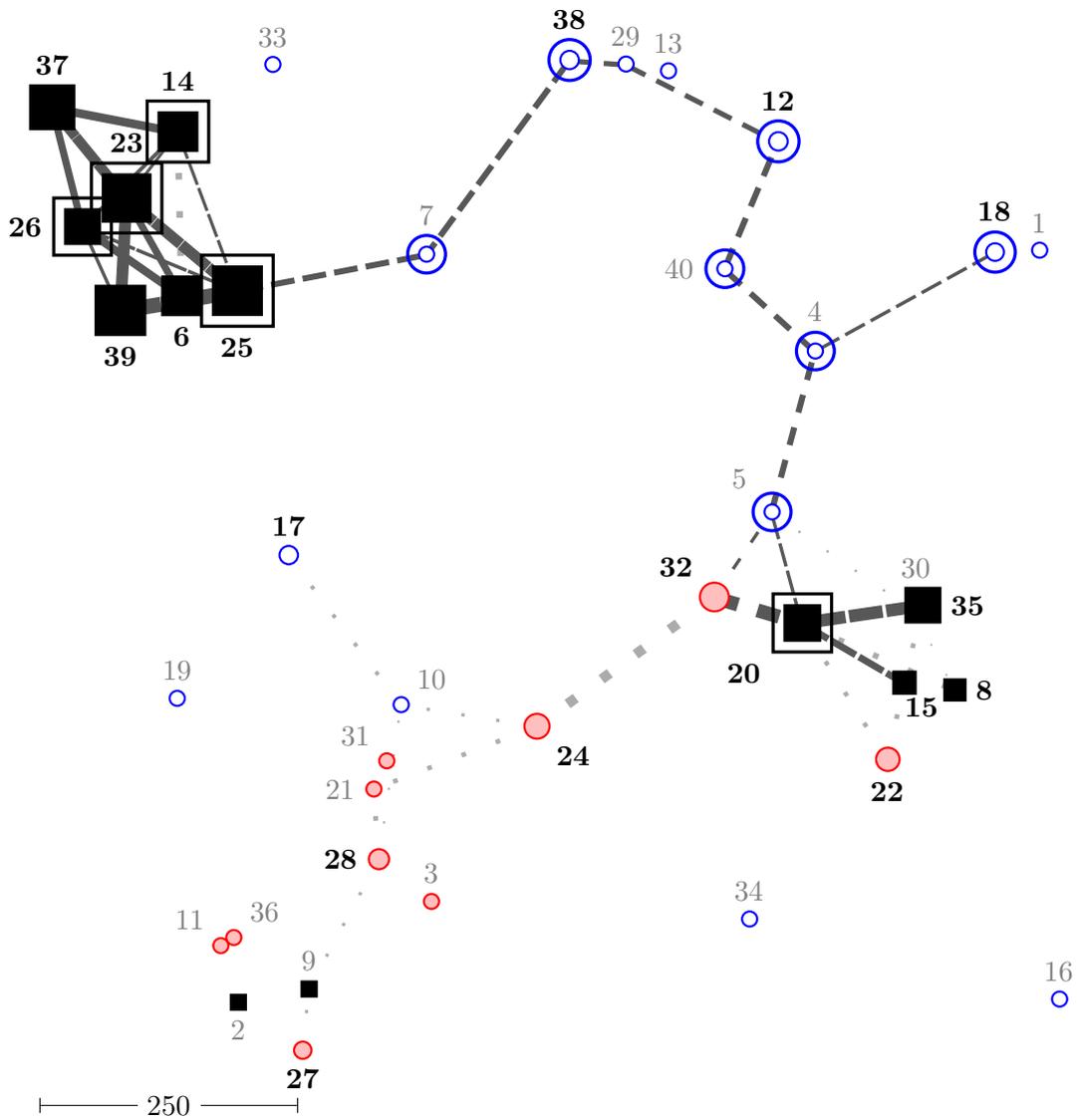

		\centering
			\begin{comment:figures}

			\end{comment:figures}
		\caption{Test instance s40w20: \LCFRLP\ -- scenario 5.}
		\label{figure:TestInstanceSFourtyWTwentyLCFRLPScenarioFive}
	\end{figure}

\paragraph{S.1 scenario sequences.}
Exploring the results of scenarios 1 to 5 (i.e.\ the scenario sequence S.1.1) in more detail, it is worth mentioning that budget increases by 25\% of the number of sub-urban locations in every move towards the sequence end, as sub-urban locations became more expensive under the assumption that urban and rural costs remain the same. In these 5 scenarios, the CFV is continually increasing, due to the fact that the additional budget is used to build further charging stations in urban and/or rural locations. 
Looking at the proportions, it seems at first glance, that in scenarios 3 to 5, the number of urban charging stations is decreasing, but considering the absolute numbers, it becomes apparent that the number of urban stations remains stable, while the number of rural stations increases. In scenario 3 it is not possible to build a sub-urban CS additionally to the five urban locations, because the leftover budget is 2 units and 3 units would be required to build a sub-urban CS. In scenario 4, there is enough budget available to build a sub-urban CS, but CFV increases most when building three rural 
CSs, requiring less budget. This is the same coverage-increasing reason which holds in scenario 5: the additional budget is used to build rural CSs. Installing rural CSs to cover long distance flows, connects two clusters characterized by high flow volumes (see Figure~\ref{figure:TestInstanceSFourtyWTwentyLCFRLPScenarioFive}).

Looking closely at scenario sequence S.1.2, the overall trend is the same as in sequence S.1.1. On the one hand, cost difference between sub-urban and rural locations is increasing, on the other hand, differences between sub-urban and urban location costs are decreasing. In these scenarios under review, when budget increases, the CFV is increasing when charging stations are built in both, urban and rural areas. 

At first glance it seems like the trend is not visible in scenario sequence S.1.3, because in scenario 10 there are urban nodes only and in scenario 11 it is optimal to install the same number of urban charging stations plus an additional CS in a sub-urban area, although prices are increasing for sub-urban locations. This can be explained by comparing scenario 10 and 11 explicitly:  in scenario 10 budget is limited to 45 units, this is sufficient for installing 6 urban CSs and there is still enough leftover budget to build a rural CS. However, it is not possible to cover additional flows by building a single additional rural CS in this particular case.  
In scenario 11, budget is increasing and limited to 48 units, which is enough to install 6 urban CSs and 2 rural or 1 sub-urban CS. Both of these possibilities increase the CFV, but based on the node category definition, a sub-urban location is able to cover more EVs. 

\paragraph{S.2 scenario sequences.}
Following a more detailed analysis of scenario 5, 9, 12, and 14: these scenarios have in common that urban and sub-urban location costs remain stable, while it becomes more expensive to build rural CSs in every move forwards the sequence end. When rural CSs become more expensive, the growing budget is used to build CSs at urban or sub-urban instead of rural locations, as they become relatively cheaper. Figure~\ref{figure:totalNumberAndProportionOfUrbanSuburbanRuralStatinsOpenedInDifferentScenariosForTheTestInstanceSFourtyWTwenty} depicts this trade-off: the proportion of urban locations is increasing, while the proportion of rural stations is decreasing steadily. 
Looking only at the proportional values of scenario 12 and 14, it seems to be counter-intuitive, as one would expect an increasing proportion of urban CSs. However, the change in composition of CSs from scenario 12 to 14 can be explained through the limitation of the budget: the number of urban CSs remains stable and with an increasing budget ($16 \cdot 1 \cdot 0.25 = 4$) under scenario 14, there is enough leftover budget to install an additional sub-urban CS ($ 7 \cdot 7 + 1 \cdot 6 = 55 \leq B\colon 55.25$).
Scenario 12 shows a remaining budget of 2 units, which cannot be used in order to build an additional CS. The budget increase of 4 units, when changing to scenario 14, can be used to build one rural ($c_r = 4$) or one sub-urban ($c_{su} = 6$) CS and therefore increase the CFV.

\paragraph{Scenario sequence S.3.}
In the following, a more in-depth look is taken on the resulting output of scenarios 1, 6, 10, and 13: looking only at the proportions shown in Figure~\ref{figure:totalNumberAndProportionOfUrbanSuburbanRuralStatinsOpenedInDifferentScenariosForTheTestInstanceSFourtyWTwenty}, one could falsely assume that scenario 10 and 13 are different to the other ones, but the absolute numbers indicate that the amount of urban stations remains the same. Due to an increasing number of sub-urban CSs, the proportion of urban CSs is decreasing. This is the same effect as described above, in the comparison of scenario 10 and 11. Changing from scenario 10 to 13 leads to an increase in budget ($27 \cdot 1 \cdot 0.25 = 6.75$). After building the very same 6 urban CSs as in scenario 10, the remaining budget (10 units) is sufficient to build an additional urban CS or two further sub-urban CSs. Looking at scenario 12's  CFV, obtained by building 7 urban CSs, it becomes clear that more EV drivers could be covered by installing two additional CSs in sub-urban location rather than one further CS in an urban location. 

\medskip

Finally, note that this instance has a specific structure. It is characterized by clusters, where a bundle of heavily frequented flows are passing by. This can be seen in Figure~\ref{figure:TestInstanceSFourtyWTwenty}. In scenario 1, building 3 urban CSs (at nodes 23, 25, and 26) within the cluster in the top left hand corner of the network, \num{493451} (49.35\%) EVs of TFV can be covered.

Installing 5 CSs (at nodes 6, 14, 23, 25, and 26) within the top, left cluster (as in scenarios 8, 10, 11, 12, 13, 14, 15) already covers \num{630618} (63.06\%) EVs. The second cluster, on the right hand side of the graph (at nodes 8, 20, and 22, as in scenario 14 and 15) covers \num{173652} (17.37\%) EVs. 

Scenario 5 (see Figure~\ref{figure:TestInstanceSFourtyWTwentyLCFRLPScenarioFive}) is the only one that connects both clusters by installing some rural CSs in order to cover long distance flows, which are travelling from one cluster to another. 

\paragraph{Analysing s80w40.}

Looking at Table~\ref{table:ResultsLCFRLPsEightywFourty} and Figure~\ref{figure:totalNumberAndProportionOfUrbanSuburbanRuralStatinsOpenedInDifferentScenariosForTheTestInstanceSEightyWFourty}, representing the proportional and absolute number of opened charging stations of test instance s80w40 in different scenarios, the existence of the general trend, described in Section~\ref {subsubsection:naLocationDependentCost}, becomes visible. 

\addtolength{\tabcolsep}{-1.5pt}
\begin{table}[htb]
\caption{Results of \LCFRLP\ for all scenarios of s80w40.} 
	\begin{center}
			\begin{tabular}[c]{c *{10}{| r} }
				\label{table:ResultsLCFRLPsEightywFourty}
Scenario  &	\multicolumn{2}{c|}{1}	&	\multicolumn{2}{c|}{2}	&	\multicolumn{2}{c|}{3}	&	\multicolumn{2}{c|}{4}	&	\multicolumn{2}{c}{5}\\
                \hline
$CFV$  &	\multicolumn{2}{r|}{\num{806077}}	&	\multicolumn{2}{r|}{\num{759374}}	&	\multicolumn{2}{r|}{\num{739058}}	&	\multicolumn{2}{r|}{\num{730415}}	&	\multicolumn{2}{r}{\num{736048}}\\

$\sum_{k \in K} x_k$   &	\multicolumn{2}{r|}{24}	&	\multicolumn{2}{r|}{18}	&	\multicolumn{2}{r|}{17}	&	\multicolumn{2}{r|}{20}	&	\multicolumn{2}{r}{23}\\
				\hline
				\hline
urban	(\sfrac{\% }{100} | \#)& 	0.1667 & 4 	& 	0.2778 & 5	& 	0.2941	& 5 & 	0.3 & 6 	& 	0.3478	 & 8\\
sub-urban (\sfrac{\% }{100} | \#)	& 	0.6250 & 15 	& 	0.6111 & 11	& 	0.5882 & 10 	& 	0.3500	& 7 & 	0.1739 & 4\\
rural (\sfrac{\% }{100} | \#)& 	0.2083 & 5 	& 	0.1111 & 2 	& 	0.1176 & 2 	& 	0.3500	 & 7 & 	0.4783 & 11\\
				\hline
				\multicolumn{11}{c}{}\\[-0.5em]
Scenario  &	\multicolumn{2}{c|}{6}	&	\multicolumn{2}{c|}{7}	&	\multicolumn{2}{c|}{8}	&	\multicolumn{2}{c|}{9}	&	\multicolumn{2}{c}{10}\\
                \hline
$CFV$  &	\multicolumn{2}{r|}{\num{783703}}	&	\multicolumn{2}{r|}{\num{759374}}	&	\multicolumn{2}{r|}{\num{747960}}	&	\multicolumn{2}{r|}{\num{741226}}	&	\multicolumn{2}{r}{\num{781279}}\\

$\sum_{k \in K} x_k$   &	\multicolumn{2}{r|}{20}	&	\multicolumn{2}{r|}{18}	&	\multicolumn{2}{r|}{17}	&	\multicolumn{2}{r|}{17}	&	\multicolumn{2}{r}{19}\\
				\hline
				\hline
urban	(\sfrac{\% }{100} | \#)& 	0.25 & 5 	& 	0.2778 & 5 	& 	0.3529 & 6  & 	0.4118 & 7	& 	0.2632	& 5\\
sub-urban (\sfrac{\% }{100} | \#)	& 	0.6000	& 12 & 	0.6111& 11 	& 	0.5294 & 9 	& 	0.4118 & 7 	& 	0.6841 & 13\\
rural (\sfrac{\% }{100} | \#)	& 	0.1500 & 3 	& 	0.1111 & 2 	& 	0.1176 & 2 	& 	0.1765 & 3 	& 	0.0526	 & 1\\
				\hline
				\multicolumn{11}{c}{}\\[-0.5em]
Scenario  &	\multicolumn{2}{c|}{11}	&	\multicolumn{2}{c|}{12}&	\multicolumn{2}{c|}{13}	&	\multicolumn{2}{c|}{14}	&	\multicolumn{2}{c}{15} \\
                \hline
$CFV$  &		\multicolumn{2}{r|}{\num{768235}}	&	\multicolumn{2}{r|}{\num{757415}}&	\multicolumn{2}{r|}{\num{787675}}	&	\multicolumn{2}{r|}{\num{777690}}	&	\multicolumn{2}{r}{\num{794699}} \\

$\sum_{k \in K} x_k$   &	\multicolumn{2}{r|}{18}	&	\multicolumn{2}{r|}{17}&	\multicolumn{2}{r|}{19}	&	\multicolumn{2}{r|}{18}	&	\multicolumn{2}{r}{19} \\
				\hline
				\hline
urban	(\sfrac{\% }{100} | \#)& 	0.3333	& 6 & 	0.3529 & 6 	& 	0.3158 & 6	& 	0.3333 & 6	& 	0.3158 & 6\\
sub-urban (\sfrac{\% }{100} | \#)	& 	0.5556	 & 10 & 	0.5882 & 10 	& 	0.5789 & 11 	& 	0.6111 & 11 	& 	0.6316 & 12  \\
rural (\sfrac{\% }{100} | \#)& 	0.1111 & 2 	& 	0.0588	 & 1 & 	0.1053 & 2 	& 	0.05556 & 1 	& 	0.0526  & 1 \\
				\hline
			\end{tabular}
	\end{center}
\end{table}
\addtolength{\tabcolsep}{+1.5pt}

			\begin{figure}[htb]
				\centering 
				\newcommand*{\maxXProportionOfUrbanSuburbanRuralStatinsOpenedInDifferentScenariosForTheTestInstanceSEightyWFourty}{15}
				\FPeval\xScaleProportionOfUrbanSuburbanRuralStatinsOpenedInDifferentScenariosForTheTestInstanceSEightyWFourty{\generalXScale / 0.838 / \maxXProportionOfUrbanSuburbanRuralStatinsOpenedInDifferentScenariosForTheTestInstanceSEightyWFourty}
				
\newcommand*{\maxYProportionOfUrbanSuburbanRuralStatinsOpenedInDifferentScenariosForTheTestInstanceSEightyWFourtyNumber}{18.0}
				\FPeval\yScaleProportionOfUrbanSuburbanRuralStatinsOpenedInDifferentScenariosForTheTestInstanceSEightyWFourtyNumber{\generalYScale / \maxYProportionOfUrbanSuburbanRuralStatinsOpenedInDifferentScenariosForTheTestInstanceSEightyWFourtyNumber}
				\newcommand*{\maxYProportionOfUrbanSuburbanRuralStatinsOpenedInDifferentScenariosForTheTestInstanceSEightyWFourtyProportion}{1.0}
				\FPeval\yScaleProportionOfUrbanSuburbanRuralStatinsOpenedInDifferentScenariosForTheTestInstanceSEightyWFourtyProportion{\generalYScale / \maxYProportionOfUrbanSuburbanRuralStatinsOpenedInDifferentScenariosForTheTestInstanceSEightyWFourtyProportion}
				\begin{comment:figures}
					\begin{tikzpicture}[xscale=\xScaleProportionOfUrbanSuburbanRuralStatinsOpenedInDifferentScenariosForTheTestInstanceSEightyWFourty, yscale=\yScaleProportionOfUrbanSuburbanRuralStatinsOpenedInDifferentScenariosForTheTestInstanceSEightyWFourtyNumber]
						\pgfgettransformentries{\xScaleTikz}{\@tempa}{\@tempa}{\yScaleTikz}{\@tempa}{\@tempa}
						
						\draw[very thin, color=gray, xstep = 1, ystep = 2] (0, 0) grid (\maxXProportionOfUrbanSuburbanRuralStatinsOpenedInDifferentScenariosForTheTestInstanceSEightyWFourty, \maxYProportionOfUrbanSuburbanRuralStatinsOpenedInDifferentScenariosForTheTestInstanceSEightyWFourtyNumber);
						
						\def\axisLabelTikzY{\axisLabel / \yScaleTikz}
						\def\axisNumberTikzY{5.5 * \axisLabelTikzY}
						\def\axisTotalNumberTikzY{3 * \axisLabelTikzY}
						\def\axisLabelTikzX{\axisLabel / \xScaleTikz}
						\draw[thick, color=gray] (5.0, -\axisNumberTikzY) -- (5.0, 18.0 + \axisNumberTikzY);
						\draw[thick, color=gray] (9.0, -\axisNumberTikzY) -- (9.0, 18.0 + \axisNumberTikzY);
						\draw[thick, color=gray] (12.0, -\axisNumberTikzY) -- (12.0, 18.0 + \axisNumberTikzY);
						\draw[thick, color=gray] (14.0, -\axisNumberTikzY) -- (14.0, 18.0 + \axisNumberTikzY);
						
						\begin{scope}[shift={(-0.5, 0)}]
							\draw[black, fill=black] (0.6, 0.0) rectangle (0.8, 4);
							\draw[red, thick, fill=red!25] (0.9, 0.0) rectangle (1.1, 15);
							\draw[blue, thick] (1.2, 0.0) rectangle (1.4, 5);
							
							\draw[black, fill=black] (1.6, 0.0) rectangle (1.8, 5);
							\draw[red, thick, fill=red!25] (1.9, 0.0) rectangle (2.1, 11);
							\draw[blue, thick] (2.2, 0.0) rectangle (2.4, 2);
							
							\draw[black, fill=black] (2.6, 0.0) rectangle (2.8, 5);
							\draw[red, thick, fill=red!25] (2.9, 0.0) rectangle (3.1, 10);
							\draw[blue, thick] (3.2, 0.0) rectangle (3.4, 2);
							
							\draw[black, fill=black] (3.6, 0.0) rectangle (3.8, 6);
							\draw[red, thick, fill=red!25] (3.9, 0.0) rectangle (4.1, 7);
							\draw[blue, thick] (4.2, 0.0) rectangle (4.4, 7);
							
							\draw[black, fill=black] (4.6, 0.0) rectangle (4.8, 8);
							\draw[red, thick, fill=red!25] (4.9, 0.0) rectangle (5.1, 4);
							\draw[blue, thick] (5.2, 0.0) rectangle (5.4, 11);
							
							\draw[black, fill=black] (5.6, 0.0) rectangle (5.8, 5);
							\draw[red, thick, fill=red!25] (5.9, 0.0) rectangle (6.1, 12);
							\draw[blue, thick] (6.2, 0.0) rectangle (6.4, 3);
							
							\draw[black, fill=black] (6.6, 0.0) rectangle (6.8, 5);
							\draw[red, thick, fill=red!25] (6.9, 0.0) rectangle (7.1, 11);
							\draw[blue, thick] (7.2, 0.0) rectangle (7.4, 2);
							
							\draw[black, fill=black] (7.6, 0.0) rectangle (7.8, 6);
							\draw[red, thick, fill=red!25] (7.9, 0.0) rectangle (8.1, 9);
							\draw[blue, thick] (8.2, 0.0) rectangle (8.4, 2);
							
							\draw[black, fill=black] (8.6, 0.0) rectangle (8.8, 7);
							\draw[red, thick, fill=red!25] (8.9, 0.0) rectangle (9.1, 7);
							\draw[blue, thick] (9.2, 0.0) rectangle (9.4, 3);
							
							\draw[black, fill=black] (9.6, 0.0) rectangle (9.8, 5);
							\draw[red, thick, fill=red!25] (9.9, 0.0) rectangle (10.1, 13);
							\draw[blue, thick] (10.2, 0.0) rectangle (10.4, 1);
							
							\draw[black, fill=black] (10.6, 0.0) rectangle (10.8, 6);
							\draw[red, thick, fill=red!25] (10.9, 0.0) rectangle (11.1, 10);
							\draw[blue, thick] (11.2, 0.0) rectangle (11.4, 2);
							
							\draw[black, fill=black] (11.6, 0.0) rectangle (11.8, 6);
							\draw[red, thick, fill=red!25] (11.9, 0.0) rectangle (12.1, 10);
							\draw[blue, thick] (12.2, 0.0) rectangle (12.4, 1);
							
							\draw[black, fill=black] (12.6, 0.0) rectangle (12.8, 6);
							\draw[red, thick, fill=red!25] (12.9, 0.0) rectangle (13.1, 11);
							\draw[blue, thick] (13.2, 0.0) rectangle (13.4, 2);
							
							\draw[black, fill=black] (13.6, 0.0) rectangle (13.8, 6);
							\draw[red, thick, fill=red!25] (13.9, 0.0) rectangle (14.1, 11);
							\draw[blue, thick] (14.2, 0.0) rectangle (14.4, 1);
							
							\draw[black, fill=black] (14.6, 0.0) rectangle (14.8, 6);
							\draw[red, thick, fill=red!25] (14.9, 0.0) rectangle (15.1, 12);
							\draw[blue, thick] (15.2, 0.0) rectangle (15.4, 1);
						\end{scope}
						
						\def\axisAdditionalLengthPlusTikzX{\axisAdditionalLengthPlus / \xScaleTikz}
						\def\axisAdditionalLengthMinusTikzX{\axisAdditionalLengthMinus / \xScaleTikz}
						\draw[-] (-\axisAdditionalLengthMinusTikzX, 0) -- (\maxXProportionOfUrbanSuburbanRuralStatinsOpenedInDifferentScenariosForTheTestInstanceSEightyWFourty, 0) -- +(0, 0) node[right] {Scenario};
						
						\foreach \pos in {1, 2, 3, 4, 5, 6, 7, 8, 9, 10, 11, 12, 13, 14, 15} \draw[shift={(\pos-0.5, 0)}] (0, -\axisLabelTikzY) node[below] {$\pos$};
						
						\def\axisAdditionalLengthPlusTikzY{\axisAdditionalLengthPlus / \yScaleTikz}
						\def\axisAdditionalLengthMinusTikzY{\axisAdditionalLengthMinus / \yScaleTikz}
						\draw[arrow] (0, -\axisAdditionalLengthMinusTikzY) -- (0, \maxYProportionOfUrbanSuburbanRuralStatinsOpenedInDifferentScenariosForTheTestInstanceSEightyWFourtyNumber) -- +(0, \axisAdditionalLengthPlusTikzY) node[above] {Total number};
						
						\foreach \pos in {2, 4, 6, 8, 10, 12, 14, 16, 18} \draw[shift={(0, \pos)}] (\axisLabelTikzX, 0) -- (-\axisLabelTikzX, 0) node[left] {$\pos$};
						
					\end{tikzpicture}
					
					
					\hspace*{0.118cm}
					\begin{tikzpicture}[xscale=\xScaleProportionOfUrbanSuburbanRuralStatinsOpenedInDifferentScenariosForTheTestInstanceSEightyWFourty, yscale=\yScaleProportionOfUrbanSuburbanRuralStatinsOpenedInDifferentScenariosForTheTestInstanceSEightyWFourtyProportion]
						\pgfgettransformentries{\xScaleTikz}{\@tempa}{\@tempa}{\yScaleTikz}{\@tempa}{\@tempa}
						
						\draw[very thin, color=gray, xstep = 1, ystep = 0.2] (0, 0) grid (\maxXProportionOfUrbanSuburbanRuralStatinsOpenedInDifferentScenariosForTheTestInstanceSEightyWFourty, \maxYProportionOfUrbanSuburbanRuralStatinsOpenedInDifferentScenariosForTheTestInstanceSEightyWFourtyProportion);
						
						\def\axisLabelTikzY{\axisLabel / \yScaleTikz}
						\def\axisPartitionTikzY{5.5 * \axisLabelTikzY}
						\def\axisLabelTikzX{\axisLabel / \xScaleTikz}
						\draw[thick, color=gray] (5.0, -\axisPartitionTikzY) -- (5.0, 1.0 + \axisPartitionTikzY);
						\draw[thick, color=gray] (9.0, -\axisPartitionTikzY) -- (9.0, 1.0 + \axisPartitionTikzY);
						\draw[thick, color=gray] (12.0, -\axisPartitionTikzY) -- (12.0, 1.0 + \axisPartitionTikzY);
						\draw[thick, color=gray] (14.0, -\axisPartitionTikzY) -- (14.0, 1.0 + \axisPartitionTikzY);
						
						\begin{scope}[shift={(-0.5, 0)}]
							\draw[black, fill=black] (0.6, 0.0) rectangle (0.8, 0.166666666666667);
							\draw[red, thick, fill=red!25] (0.9, 0.0) rectangle (1.1, 0.625);
							\draw[blue, thick] (1.2, 0.0) rectangle (1.4, 0.208333333333333);
							
							\draw[black, fill=black] (1.6, 0.0) rectangle (1.8, 0.277777777777778);
							\draw[red, thick, fill=red!25] (1.9, 0.0) rectangle (2.1, 0.611111111111111);
							\draw[blue, thick] (2.2, 0.0) rectangle (2.4, 0.111111111111111);
							
							\draw[black, fill=black] (2.6, 0.0) rectangle (2.8, 0.294117647058823);
							\draw[red, thick, fill=red!25] (2.9, 0.0) rectangle (3.1, 0.588235294117647);
							\draw[blue, thick] (3.2, 0.0) rectangle (3.4, 0.117647058823529);
							
							\draw[black, fill=black] (3.6, 0.0) rectangle (3.8, 0.3);
							\draw[red, thick, fill=red!25] (3.9, 0.0) rectangle (4.1, 0.35);
							\draw[blue, thick] (4.2, 0.0) rectangle (4.4, 0.35);
							
							\draw[black, fill=black] (4.6, 0.0) rectangle (4.8, 0.347826086956522);
							\draw[red, thick, fill=red!25] (4.9, 0.0) rectangle (5.1, 0.173913043478261);
							\draw[blue, thick] (5.2, 0.0) rectangle (5.4, 0.478260869565217);
							
							\draw[black, fill=black] (5.6, 0.0) rectangle (5.8, 0.25);
							\draw[red, thick, fill=red!25] (5.9, 0.0) rectangle (6.1, 0.6);
							\draw[blue, thick] (6.2, 0.0) rectangle (6.4, 0.15);
							
							\draw[black, fill=black] (6.6, 0.0) rectangle (6.8, 0.277777777777778);
							\draw[red, thick, fill=red!25] (6.9, 0.0) rectangle (7.1, 0.611111111111111);
							\draw[blue, thick] (7.2, 0.0) rectangle (7.4, 0.111111111111111);
							
							\draw[black, fill=black] (7.6, 0.0) rectangle (7.8, 0.352941176470588);
							\draw[red, thick, fill=red!25] (7.9, 0.0) rectangle (8.1, 0.529411764705882);
							\draw[blue, thick] (8.2, 0.0) rectangle (8.4, 0.117647058823529);
							
							\draw[black, fill=black] (8.6, 0.0) rectangle (8.8, 0.411764705882353);
							\draw[red, thick, fill=red!25] (8.9, 0.0) rectangle (9.1, 0.411764705882353);
							\draw[blue, thick] (9.2, 0.0) rectangle (9.4, 0.176470588235294);
							
							\draw[black, fill=black] (9.6, 0.0) rectangle (9.8, 0.263157894736842);
							\draw[red, thick, fill=red!25] (9.9, 0.0) rectangle (10.1, 0.68421052631579);
							\draw[blue, thick] (10.2, 0.0) rectangle (10.4, 0.052631578947369);
							
							\draw[black, fill=black] (10.6, 0.0) rectangle (10.8, 0.333333333333333);
							\draw[red, thick, fill=red!25] (10.9, 0.555555555555556) rectangle (11.1, 0.0);
							\draw[blue, thick] (11.2, 0.0) rectangle (11.4, 0.111111111111111);
							
							\draw[black, fill=black] (11.6, 0.0) rectangle (11.8, 0.352941176470588);
							\draw[red, thick, fill=red!25] (11.9, 0.588235294117647) rectangle (12.1, 0.0);
							\draw[blue, thick] (12.2, 0.0) rectangle (12.4, 0.058823529411765);
							
							\draw[black, fill=black] (12.6, 0.0) rectangle (12.8, 0.31578947368421);
							\draw[red, thick, fill=red!25] (12.9, 0.578947368421053) rectangle (13.1, 0.0);
							\draw[blue, thick] (13.2, 0.0) rectangle (13.4, 0.105263157894737);
							
							\draw[black, fill=black] (13.6, 0.0) rectangle (13.8, 0.333333333333333);
							\draw[red, thick, fill=red!25] (13.9, 0.611111111111111) rectangle (14.1, 0.0);
							\draw[blue, thick] (14.2, 0.0) rectangle (14.4, 0.055555555555556);
							
							\draw[black, fill=black] (14.6, 0.0) rectangle (14.8, 0.31578947368421);
							\draw[red, thick, fill=red!25] (14.9, 0.0) rectangle (15.1, 0.631578947368421);
							\draw[blue, thick] (15.2, 0.0) rectangle (15.4, 0.052631578947369);
						\end{scope}

						\def\axisAdditionalLengthPlusTikzX{\axisAdditionalLengthPlus / \xScaleTikz}
						\def\axisAdditionalLengthMinusTikzX{\axisAdditionalLengthMinus / \xScaleTikz}
						\draw[-] (-\axisAdditionalLengthMinusTikzX, 0) -- (\maxXProportionOfUrbanSuburbanRuralStatinsOpenedInDifferentScenariosForTheTestInstanceSEightyWFourty, 0) -- +(0, 0) node[right] {Scenario};
						
						\foreach \pos in {1, 2, 3, 4, 5, 6, 7, 8, 9, 10, 11, 12, 13, 14, 15} \draw[shift={(\pos-0.5, 0)}] (0, -\axisLabelTikzY) node[below] {$\pos$};
						
						\def\axisAdditionalLengthPlusTikzY{\axisAdditionalLengthPlus / \yScaleTikz}
						\def\axisAdditionalLengthMinusTikzY{\axisAdditionalLengthMinus / \yScaleTikz}
						\draw[arrow] (0, -\axisAdditionalLengthMinusTikzY) -- (0, \maxYProportionOfUrbanSuburbanRuralStatinsOpenedInDifferentScenariosForTheTestInstanceSEightyWFourtyProportion) -- +(0, \axisAdditionalLengthPlusTikzY) node[above] {Proportion};
						
						\foreach \pos in {0.0, 0.2, 0.4, 0.6, 0.8, 1.0} \draw[shift={(0, \pos)}] (\axisLabelTikzX, 0) -- (-\axisLabelTikzX, 0) node[left] {$\pos$};
						
					\end{tikzpicture}
				\end{comment:figures}
				\caption[Total number and proportion of urban/suburban/rural statins opened in different scenarios for the test instance s80w40.]{Total number and proportion of urban/suburban/rural statins opened in different scenarios for the test instance s80w40. The left bar \begin{tikzpicture}[xscale=\xScaleProportionOfUrbanSuburbanRuralStatinsOpenedInDifferentScenariosForTheTestInstanceSEightyWFourty, yscale=1]\draw[black, fill=black] (0.0, 0.0) rectangle (0.2, 0.24);\end{tikzpicture} stays for urban, the middle bar \begin{tikzpicture}[xscale=\xScaleProportionOfUrbanSuburbanRuralStatinsOpenedInDifferentScenariosForTheTestInstanceSEightyWFourty, yscale=1]\draw[red, thick, fill=red!25] (0.0, 0.0) rectangle (0.2, 0.24);\end{tikzpicture} for suburban, and the right bar \begin{tikzpicture}[xscale=\xScaleProportionOfUrbanSuburbanRuralStatinsOpenedInDifferentScenariosForTheTestInstanceSEightyWFourty, yscale=1]\draw[blue, thick] (0.0, 0.0) rectangle (0.2, 0.24);\end{tikzpicture} for rural areas.}
				\label{figure:totalNumberAndProportionOfUrbanSuburbanRuralStatinsOpenedInDifferentScenariosForTheTestInstanceSEightyWFourty}
			\end{figure}

\paragraph{S.1 scenario sequences.}
These scenarios are characterized by increasing sub-urban location costs, whereas urban and rural costs remain stable. Consequently, the number of CSs built in sub-urban areas decreases, while the number of urban CSs increases. When sub-urban costs increase, in some scenarios there is not enough budget to install the very same CSs infrastructure as in the previous scenario. Therefore, a trade-off between several sub-urban CSs and one urban CS might lead to an increasing CFV when comparing it to the case that one tries to install as many CSs at the very same locations as in the previous scenario until the whole budget is consumed. However, the CFV reached with the optimal charging station placement of the current scenario is less than in the previous scenario. This effect occurs when changing from scenario 1 to 2. 
Considering the first five scenarios (i.e.\ the scenario sequence S.1.1), CFV is decreasing except from changing from scenario 4 to 5. This can be explained as follows: sub-urban costs increase by one unit and thus, the budget increases by 25\% of the number of sub-urban locations ($28 \cdot 1 \cdot 0.25 = 7$). Following, in order to have enough budget to build an additional, more expensive sub-urban CS, there have to be 4 sub-urban locations (associated with an cost increase of 1 unit) within the network.
Therefore, if the number of opened sub-urban CSs in the previous scenario is below this threshold ($4\ \cdot $ \# opened rural CSs $\leq 28$, see Table~\ref{table:CostParameterDefinition}), it is possible to install the very same CS infrastructure in the current scenario and therefore cover the same flows resulting in the very same CFV.
Covering additional/different flows and therefore increase CFV can be reached by using the growing budget to build rural CSs or substitute some sub-urban CSs (which became  more expensive) for urban and therefore more frequented ones. In scenario 5, installing the very same CSs as in scenario 4 requires the entire budget ($6 \cdot 7 + 7 \cdot 6 + 7 \cdot 1 = 91 \leq B\colon 91$) and therefore covers the very same EVs, but budget allows the substitution of some sub-urban locations for urban and rural locations ($8 \cdot 7 + 4 \cdot 6 + 11 \cdot 1 = 91 \leq B\colon 91$), leading to an increase in CFV. This is also the reason for continuously decreasing CFV when looking at scenarios 6--9, 10--12, and 13--14, as the threshold of 7 opened sub-urban CSs ($4 \cdot 7 \leq 28$) is never undercut. 

\paragraph{Scenario sequence S.2.1.}
Looking more closely the scenarios 5, 9, 12, 14, and 15, where rural costs is the only different price component: CFV increases continuously, as it is possible to install the very same infrastructure of the previous scenario in the current one, because the total number of CSs assigned to the rising cost category (rural locations) account for 28 locations and therefore a threshold of 7 rural CSs results. If this threshold is not exceeded (as it is the case in all the scenarios under review), the budget increase covers the increasing costs for rural CSs. 

A general trend for a decreasing proportion of rural CSs, when gradually changing from from scenario 5 to 14, can be observed in this test instance as well when looking at absolute and proportional numbers. The trend emerges from the fact that building CSs in rural locations becomes more expensive.

\paragraph{Scenario sequence S.3.}
Comparing scenario 1, 6, 10, and 13, all have slight cost difference between rural and sub-urban locations, but significant, even though decreasing cost differences between sub-urban and urban areas when moving towards the sequence end. The general trend of increasing proportion of urban locations and decreasing proportions of sub-urban and rural locations can be identified.
Considering the change from scenario 6 to 10 seems counter-intuitive, but as it is not possible to install the very same CS infrastructure of scenario 6 in scenario 10, CFV decreases as charging station costs increases. As urban costs remain stable, the same number of urban CSs can be build. Resulting from that a leftover budget ($B = 91 - 5 \cdot 7 = 56$) can be used to build at most 14 CSs in sub-urban regions. CFV is maximized by installing 13 CSs in sub-urban regions and one CS in a rural area.

\paragraph{Analysing s100w50.}

Table~\ref{table:ResultsLCFRLPsHundredwFifty} summarizes the results when the 15 different cost scenarios of the \LCFRLP\ are applied to the test instance s100w50. 

\addtolength{\tabcolsep}{-1.5pt}
\begin{table}[htb!]
\caption{Results of \LCFRLP\ for all scenarios of s100w50.} 
	\begin{center}
			\begin{tabular}[c]{c *{10}{| r} }
				\label{table:ResultsLCFRLPsHundredwFifty}
Scenario  &	\multicolumn{2}{c|}{1}	&	\multicolumn{2}{c|}{2}	&	\multicolumn{2}{c|}{3}	&	\multicolumn{2}{c|}{4}	&	\multicolumn{2}{c}{5}\\
                \hline
$CFV$  &	\multicolumn{2}{r|}{\num{810776}}	&	\multicolumn{2}{r|}{\num{770574}}	&	\multicolumn{2}{r|}{\num{742221}}	&	\multicolumn{2}{r|}{\num{732913}	}&	\multicolumn{2}{r}{\num{720819}}\\

$\sum_{k \in K} x_k$   &	\multicolumn{2}{r|}{27}	&	\multicolumn{2}{r|}{28}	&	\multicolumn{2}{r|}{27}	&	\multicolumn{2}{r|}{25}	&	\multicolumn{2}{r}{25}\\
				\hline
				\hline
urban (\sfrac{\% }{100} | \#) 	&	0.1111 & 3 	&	0.0714 & 2 	&	0.1111 & 3 	&	0.12	& 3 &	0.2 & 5\\
sub-urban (\sfrac{\% }{100} | \#)	&	0.6296 & 17 	&	0.5714 & 16 &	0.4444 & 12 	&	0.48 & 12 	&	0.36 & 9\\
rural (\sfrac{\% }{100} | \#)&	0.2593 & 7 	&	0.3571 & 10 	&	0.4444	& 12 &	0.4	& 10 &	0.44 & 11\\
				\hline
				\multicolumn{11}{c}{}\\[-0.5em]
Scenario  &	\multicolumn{2}{c|}{6}	&	\multicolumn{2}{c|}{7}	&	\multicolumn{2}{c|}{8}	&	\multicolumn{2}{c|}{9}	&	\multicolumn{2}{c}{10}\\
                \hline
$CFV$  &	\multicolumn{2}{r|}{\num{795154}	}&	\multicolumn{2}{r|}{\num{759413}}	&	\multicolumn{2}{r|}{\num{741786}}	&	\multicolumn{2}{r|}{\num{732507}}	&	\multicolumn{2}{r}{\num{788278}}\\

$\sum_{k \in K} x_k$   &	\multicolumn{2}{r|}{26}	&	\multicolumn{2}{r|}{24}	&	\multicolumn{2}{r|}{23}	&	\multicolumn{2}{r|}{24}	&	\multicolumn{2}{r}{25}\\
				\hline
				\hline
urban (\sfrac{\% }{100} | \#) 	&	0.1154	 & 3 &	0.125 & 3 	&	0.1739 & 4 	&	0.1667 & 4 	&	0.16 & 4\\
sub-urban (\sfrac{\% }{100} | \#)	&	0.6538 & 17 	&	0.625 & 15 	&	0.5217 & 12 	&	0.4583	 & 11 &	0.56 & 14\\
rural (\sfrac{\% }{100} | \#)&	0.2308 & 6 	&	0.25 & 6 	&	0.3043 & 7 	&	0.375 & 9 	&	0.28	 & 7\\
				\hline
				\multicolumn{11}{c}{}\\[-0.5em]
Scenario  &	\multicolumn{2}{c|}{11}	&	\multicolumn{2}{c|}{12}&	\multicolumn{2}{c|}{13}	&	\multicolumn{2}{c|}{14}	&	\multicolumn{2}{c}{15} \\
                \hline
$CFV$  &	\multicolumn{2}{r|}{\num{766343}}	&	\multicolumn{2}{r|}{\num{753510}}&	\multicolumn{2}{r|}{\num{798888}}	&	\multicolumn{2}{r|}{\num{778154}}	&	\multicolumn{2}{r}{\num{798888}} \\

$\sum_{k \in K} x_k$   &	\multicolumn{2}{r|}{24}	&	\multicolumn{2}{r|}{23}&	\multicolumn{2}{r|}{25}	&	\multicolumn{2}{r|}{24}	&	\multicolumn{2}{r}{25} \\
				\hline
				\hline
urban (\sfrac{\% }{100} | \#) 	&	0.1667 & 4 	&	0.1739	& 4 &	0.16	& 4 &	0.1667 & 4 	&	0.16 & 4 \\
sub-urban (\sfrac{\% }{100} | \#)	&	0.5417 & 13 	&	0.5652	& 13 &	0.6 & 15 	&	0.5833 & 14 	&	0.6  & 15 \\
rural (\sfrac{\% }{100} | \#)&	0.2917	& 7 &	0.2609	 & 6 &	0.24	 & 6 &	0.25	& 6 &	0.24 & 6 \\
				\hline
			\end{tabular}
	\end{center}
\end{table}
\addtolength{\tabcolsep}{+1.5pt}

\paragraph{S.1 scenario sequences.}
In Figure~\ref{figure:totalNumberAndProportionOfUrbanSuburbanRuralStatinsOpenedInDifferentScenariosForTheTestInstanceSHundredWFifty} one can clearly see, if sub-urban construction costs increase, the proportion of sub-urban CSs within the infrastructure network is decreasing. 
Increasing sub-urban costs leads to an increase in budget ($38 \cdot 1 \cdot 0.25 = 9.5$). Thus, if the number of open sub-urban CSs in the previous scenario exceeds a threshold of 9 CSs ($4\ \cdot $ \# opened rural CSs $\leq 38$), it is not possible to install the very same CSs in the current scenario. As a consequence, the number of sub-urban CSs is decreasing when costs are increasing. In order to maximize CFV it might make sense to substitute some sub-urban CSs for an urban CS, this occurs when changing from scenario 4 to 5. As the cost differences between sub-urban and urban locations becomes less significant, the number of sub-urban CSs that have to be exchanged to build an urban charging station is decreasing. E.g.\ in scenario 2 one has to exchange 3 sub-urban CSs in order to have enough budget available to build a CS in an urban area, while in scenario 5 the failure to set up a sub-urban CS results in an available budget of 6 units and it costs 7 units to build an urban CS. 

In this test instance the CFV is continually decreasing when considering scenarios 1--5, 6--9, 10--12, and 13--14. 

			\begin{figure}[htb]
				\centering 
				\newcommand*{\maxXProportionOfUrbanSuburbanRuralStatinsOpenedInDifferentScenariosForTheTestInstanceSHundredWFifty}{15}
				\FPeval\xScaleProportionOfUrbanSuburbanRuralStatinsOpenedInDifferentScenariosForTheTestInstanceSHundredWFifty{\generalXScale / 0.838 / \maxXProportionOfUrbanSuburbanRuralStatinsOpenedInDifferentScenariosForTheTestInstanceSHundredWFifty}
				
\newcommand*{\maxYProportionOfUrbanSuburbanRuralStatinsOpenedInDifferentScenariosForTheTestInstanceSHundredWFiftyNumber}{18.0}
				\FPeval\yScaleProportionOfUrbanSuburbanRuralStatinsOpenedInDifferentScenariosForTheTestInstanceSHundredWFiftyNumber{\generalYScale / \maxYProportionOfUrbanSuburbanRuralStatinsOpenedInDifferentScenariosForTheTestInstanceSHundredWFiftyNumber}
				\newcommand*{\maxYProportionOfUrbanSuburbanRuralStatinsOpenedInDifferentScenariosForTheTestInstanceSHundredWFiftyProportion}{1.0}
				\FPeval\yScaleProportionOfUrbanSuburbanRuralStatinsOpenedInDifferentScenariosForTheTestInstanceSHundredWFiftyProportion{\generalYScale / \maxYProportionOfUrbanSuburbanRuralStatinsOpenedInDifferentScenariosForTheTestInstanceSHundredWFiftyProportion}
				\begin{comment:figures}
					\begin{tikzpicture}[xscale=\xScaleProportionOfUrbanSuburbanRuralStatinsOpenedInDifferentScenariosForTheTestInstanceSHundredWFifty, yscale=\yScaleProportionOfUrbanSuburbanRuralStatinsOpenedInDifferentScenariosForTheTestInstanceSHundredWFiftyNumber]
						\pgfgettransformentries{\xScaleTikz}{\@tempa}{\@tempa}{\yScaleTikz}{\@tempa}{\@tempa}
						
						\draw[very thin, color=gray, xstep = 1, ystep = 2] (0, 0) grid (\maxXProportionOfUrbanSuburbanRuralStatinsOpenedInDifferentScenariosForTheTestInstanceSHundredWFifty, \maxYProportionOfUrbanSuburbanRuralStatinsOpenedInDifferentScenariosForTheTestInstanceSHundredWFiftyNumber);
						
						\def\axisLabelTikzY{\axisLabel / \yScaleTikz}
						\def\axisNumberTikzY{5.5 * \axisLabelTikzY}
						\def\axisTotalNumberTikzY{3 * \axisLabelTikzY}
						\def\axisLabelTikzX{\axisLabel / \xScaleTikz}
						\draw[thick, color=gray] (5.0, -\axisNumberTikzY) -- (5.0, 18.0 + \axisNumberTikzY);
						\draw[thick, color=gray] (9.0, -\axisNumberTikzY) -- (9.0, 18.0 + \axisNumberTikzY);
						\draw[thick, color=gray] (12.0, -\axisNumberTikzY) -- (12.0, 18.0 + \axisNumberTikzY);
						\draw[thick, color=gray] (14.0, -\axisNumberTikzY) -- (14.0, 18.0 + \axisNumberTikzY);
						
						\begin{scope}[shift={(-0.5, 0)}]
							\draw[black, fill=black] (0.6, 0.0) rectangle (0.8, 3);
							\draw[red, thick, fill=red!25] (0.9, 0.0) rectangle (1.1, 17);
							\draw[blue, thick] (1.2, 0.0) rectangle (1.4, 7);
							
							\draw[black, fill=black] (1.6, 0.0) rectangle (1.8, 2);
							\draw[red, thick, fill=red!25] (1.9, 0.0) rectangle (2.1, 16);
							\draw[blue, thick] (2.2, 0.0) rectangle (2.4, 10);
							
							\draw[black, fill=black] (2.6, 0.0) rectangle (2.8, 3);
							\draw[red, thick, fill=red!25] (2.9, 0.0) rectangle (3.1, 12);
							\draw[blue, thick] (3.2, 0.0) rectangle (3.4, 12);
							
							\draw[black, fill=black] (3.6, 0.0) rectangle (3.8, 3);
							\draw[red, thick, fill=red!25] (3.9, 0.0) rectangle (4.1, 12);
							\draw[blue, thick] (4.2, 0.0) rectangle (4.4, 10);
							
							\draw[black, fill=black] (4.6, 0.0) rectangle (4.8, 5);
							\draw[red, thick, fill=red!25] (4.9, 0.0) rectangle (5.1, 9);
							\draw[blue, thick] (5.2, 0.0) rectangle (5.4, 11);
							
							\draw[black, fill=black] (5.6, 0.0) rectangle (5.8, 3);
							\draw[red, thick, fill=red!25] (5.9, 0.0) rectangle (6.1, 17);
							\draw[blue, thick] (6.2, 0.0) rectangle (6.4, 6);
							
							\draw[black, fill=black] (6.6, 0.0) rectangle (6.8, 3);
							\draw[red, thick, fill=red!25] (6.9, 0.0) rectangle (7.1, 15);
							\draw[blue, thick] (7.2, 0.0) rectangle (7.4, 6);
							
							\draw[black, fill=black] (7.6, 0.0) rectangle (7.8, 4);
							\draw[red, thick, fill=red!25] (7.9, 0.0) rectangle (8.1, 12);
							\draw[blue, thick] (8.2, 0.0) rectangle (8.4, 7);
							
							\draw[black, fill=black] (8.6, 0.0) rectangle (8.8, 4);
							\draw[red, thick, fill=red!25] (8.9, 0.0) rectangle (9.1, 11);
							\draw[blue, thick] (9.2, 0.0) rectangle (9.4, 9);
							
							\draw[black, fill=black] (9.6, 0.0) rectangle (9.8, 4);
							\draw[red, thick, fill=red!25] (9.9, 0.0) rectangle (10.1, 14);
							\draw[blue, thick] (10.2, 0.0) rectangle (10.4, 7);
							
							\draw[black, fill=black] (10.6, 0.0) rectangle (10.8, 4);
							\draw[red, thick, fill=red!25] (10.9, 0.0) rectangle (11.1, 13);
							\draw[blue, thick] (11.2, 0.0) rectangle (11.4, 7);
							
							\draw[black, fill=black] (11.6, 0.0) rectangle (11.8, 4);
							\draw[red, thick, fill=red!25] (11.9, 0.0) rectangle (12.1, 13);
							\draw[blue, thick] (12.2, 0.0) rectangle (12.4, 6);
							
							\draw[black, fill=black] (12.6, 0.0) rectangle (12.8, 4);
							\draw[red, thick, fill=red!25] (12.9, 0.0) rectangle (13.1, 15);
							\draw[blue, thick] (13.2, 0.0) rectangle (13.4, 6);
							
							\draw[black, fill=black] (13.6, 0.0) rectangle (13.8, 4);
							\draw[red, thick, fill=red!25] (13.9, 0.0) rectangle (14.1, 14);
							\draw[blue, thick] (14.2, 0.0) rectangle (14.4, 6);
							
							\draw[black, fill=black] (14.6, 0.0) rectangle (14.8, 4);
							\draw[red, thick, fill=red!25] (14.9, 0.0) rectangle (15.1, 15);
							\draw[blue, thick] (15.2, 0.0) rectangle (15.4, 6);
						\end{scope}
						
						\def\axisAdditionalLengthPlusTikzX{\axisAdditionalLengthPlus / \xScaleTikz}
						\def\axisAdditionalLengthMinusTikzX{\axisAdditionalLengthMinus / \xScaleTikz}
						\draw[-] (-\axisAdditionalLengthMinusTikzX, 0) -- (\maxXProportionOfUrbanSuburbanRuralStatinsOpenedInDifferentScenariosForTheTestInstanceSHundredWFifty, 0) -- +(0, 0) node[right] {Scenario};
						
						\foreach \pos in {1, 2, 3, 4, 5, 6, 7, 8, 9, 10, 11, 12, 13, 14, 15} \draw[shift={(\pos-0.5, 0)}] (0, -\axisLabelTikzY) node[below] {$\pos$};
						
						\def\axisAdditionalLengthPlusTikzY{\axisAdditionalLengthPlus / \yScaleTikz}
						\def\axisAdditionalLengthMinusTikzY{\axisAdditionalLengthMinus / \yScaleTikz}
						\draw[arrow] (0, -\axisAdditionalLengthMinusTikzY) -- (0, \maxYProportionOfUrbanSuburbanRuralStatinsOpenedInDifferentScenariosForTheTestInstanceSHundredWFiftyNumber) -- +(0, \axisAdditionalLengthPlusTikzY) node[above] {Total number};
						
						\foreach \pos in {2, 4, 6, 8, 10, 12, 14, 16, 18} \draw[shift={(0, \pos)}] (\axisLabelTikzX, 0) -- (-\axisLabelTikzX, 0) node[left] {$\pos$};
						
					\end{tikzpicture}
					
					
					\hspace*{0.118cm}
					\begin{tikzpicture}[xscale=\xScaleProportionOfUrbanSuburbanRuralStatinsOpenedInDifferentScenariosForTheTestInstanceSHundredWFifty, yscale=\yScaleProportionOfUrbanSuburbanRuralStatinsOpenedInDifferentScenariosForTheTestInstanceSHundredWFiftyProportion]
						\pgfgettransformentries{\xScaleTikz}{\@tempa}{\@tempa}{\yScaleTikz}{\@tempa}{\@tempa}
						
						\draw[very thin, color=gray, xstep = 1, ystep = 0.2] (0, 0) grid (\maxXProportionOfUrbanSuburbanRuralStatinsOpenedInDifferentScenariosForTheTestInstanceSHundredWFifty, \maxYProportionOfUrbanSuburbanRuralStatinsOpenedInDifferentScenariosForTheTestInstanceSHundredWFiftyProportion);
						
						\def\axisLabelTikzY{\axisLabel / \yScaleTikz}
						\def\axisPartitionTikzY{5.5 * \axisLabelTikzY}
						\def\axisLabelTikzX{\axisLabel / \xScaleTikz}
						\draw[thick, color=gray] (5.0, -\axisPartitionTikzY) -- (5.0, 1.0 + \axisPartitionTikzY);
						\draw[thick, color=gray] (9.0, -\axisPartitionTikzY) -- (9.0, 1.0 + \axisPartitionTikzY);
						\draw[thick, color=gray] (12.0, -\axisPartitionTikzY) -- (12.0, 1.0 + \axisPartitionTikzY);
						\draw[thick, color=gray] (14.0, -\axisPartitionTikzY) -- (14.0, 1.0 + \axisPartitionTikzY);
						
						\begin{scope}[shift={(-0.5, 0)}]
							\draw[black, fill=black] (0.6, 0.0) rectangle (0.8, 0.111111111111111);
							\draw[red, thick, fill=red!25] (0.9, 0.0) rectangle (1.1, 0.62962962962963);
							\draw[blue, thick] (1.2, 0.0) rectangle (1.4, 0.259259259259259);
							
							\draw[black, fill=black] (1.6, 0.0) rectangle (1.8, 0.071428571428572);
							\draw[red, thick, fill=red!25] (1.9, 0.0) rectangle (2.1, 0.571428571428571);
							\draw[blue, thick] (2.2, 0.0) rectangle (2.4, 0.357142857142857);
							
							\draw[black, fill=black] (2.6, 0.0) rectangle (2.8, 0.111111111111111);
							\draw[red, thick, fill=red!25] (2.9, 0.0) rectangle (3.1, 0.444444444444444);
							\draw[blue, thick] (3.2, 0.0) rectangle (3.4, 0.444444444444444);
							
							\draw[black, fill=black] (3.6, 0.0) rectangle (3.8, 0.12);
							\draw[red, thick, fill=red!25] (3.9, 0.0) rectangle (4.1, 0.48);
							\draw[blue, thick] (4.2, 0.0) rectangle (4.4, 0.4);
							
							\draw[black, fill=black] (4.6, 0.0) rectangle (4.8, 0.2);
							\draw[red, thick, fill=red!25] (4.9, 0.0) rectangle (5.1, 0.36);
							\draw[blue, thick] (5.2, 0.0) rectangle (5.4, 0.44);
							
							\draw[black, fill=black] (5.6, 0.0) rectangle (5.8, 0.115384615384615);
							\draw[red, thick, fill=red!25] (5.9, 0.0) rectangle (6.1, 0.653846153846154);
							\draw[blue, thick] (6.2, 0.0) rectangle (6.4, 0.230769230769231);
							
							\draw[black, fill=black] (6.6, 0.0) rectangle (6.8, 0.125);
							\draw[red, thick, fill=red!25] (6.9, 0.0) rectangle (7.1, 0.625);
							\draw[blue, thick] (7.2, 0.0) rectangle (7.4, 0.25);
							
							\draw[black, fill=black] (7.6, 0.0) rectangle (7.8, 0.173913043478261);
							\draw[red, thick, fill=red!25] (7.9, 0.0) rectangle (8.1, 0.521739130434783);
							\draw[blue, thick] (8.2, 0.0) rectangle (8.4, 0.304347826086957);
							
							\draw[black, fill=black] (8.6, 0.0) rectangle (8.8, 0.166666666666667);
							\draw[red, thick, fill=red!25] (8.9, 0.0) rectangle (9.1, 0.458333333333333);
							\draw[blue, thick] (9.2, 0.0) rectangle (9.4, 0.375);
							
							\draw[black, fill=black] (9.6, 0.0) rectangle (9.8, 0.16);
							\draw[red, thick, fill=red!25] (9.9, 0.0) rectangle (10.1, 0.56);
							\draw[blue, thick] (10.2, 0.0) rectangle (10.4, 0.28);
							
							\draw[black, fill=black] (10.6, 0.0) rectangle (10.8, 0.166666666666667);
							\draw[red, thick, fill=red!25] (10.9, 0.541666666666667) rectangle (11.1, 0.0);
							\draw[blue, thick] (11.2, 0.0) rectangle (11.4, 0.291666666666667);
							
							\draw[black, fill=black] (11.6, 0.0) rectangle (11.8, 0.173913043478261);
							\draw[red, thick, fill=red!25] (11.9, 0.0) rectangle (12.1, 0.565217391304348);
							\draw[blue, thick] (12.2, 0.0) rectangle (12.4, 0.260869565217391);
							
							\draw[black, fill=black] (12.6, 0.0) rectangle (12.8, 0.16);
							\draw[red, thick, fill=red!25] (12.9, 0.6) rectangle (13.1, 0.0);
							\draw[blue, thick] (13.2, 0.0) rectangle (13.4, 0.24);
							
							\draw[black, fill=black] (13.6, 0.0) rectangle (13.8, 0.166666666666667);
							\draw[red, thick, fill=red!25] (13.9, 0.583333333333333) rectangle (14.1, 0.0);
							\draw[blue, thick] (14.2, 0.0) rectangle (14.4, 0.25);
							
							\draw[black, fill=black] (14.6, 0.0) rectangle (14.8, 0.16);
							\draw[red, thick, fill=red!25] (14.9, 0.0) rectangle (15.1, 0.6);
							\draw[blue, thick] (15.2, 0.0) rectangle (15.4, 0.24);
						\end{scope}

						\def\axisAdditionalLengthPlusTikzX{\axisAdditionalLengthPlus / \xScaleTikz}
						\def\axisAdditionalLengthMinusTikzX{\axisAdditionalLengthMinus / \xScaleTikz}
						\draw[-] (-\axisAdditionalLengthMinusTikzX, 0) -- (\maxXProportionOfUrbanSuburbanRuralStatinsOpenedInDifferentScenariosForTheTestInstanceSHundredWFifty, 0) -- +(0, 0) node[right] {Scenario};
						
						\foreach \pos in {1, 2, 3, 4, 5, 6, 7, 8, 9, 10, 11, 12, 13, 14, 15} \draw[shift={(\pos-0.5, 0)}] (0, -\axisLabelTikzY) node[below] {$\pos$};
						
						\def\axisAdditionalLengthPlusTikzY{\axisAdditionalLengthPlus / \yScaleTikz}
						\def\axisAdditionalLengthMinusTikzY{\axisAdditionalLengthMinus / \yScaleTikz}
						\draw[arrow] (0, -\axisAdditionalLengthMinusTikzY) -- (0, \maxYProportionOfUrbanSuburbanRuralStatinsOpenedInDifferentScenariosForTheTestInstanceSHundredWFiftyProportion) -- +(0, \axisAdditionalLengthPlusTikzY) node[above] {Proportion};
						
						\foreach \pos in {0.0, 0.2, 0.4, 0.6, 0.8, 1.0} \draw[shift={(0, \pos)}] (\axisLabelTikzX, 0) -- (-\axisLabelTikzX, 0) node[left] {$\pos$};
						
					\end{tikzpicture}
				\end{comment:figures}
				\caption[Total number and proportion of urban/suburban/rural statins opened in different scenarios for the test instance s100w50.]{Total number and proportion of urban/suburban/rural statins opened in different scenarios for the test instance s100w50. The left bar \begin{tikzpicture}[xscale=\xScaleProportionOfUrbanSuburbanRuralStatinsOpenedInDifferentScenariosForTheTestInstanceSHundredWFifty, yscale=1]\draw[black, fill=black] (0.0, 0.0) rectangle (0.2, 0.24);\end{tikzpicture} stays for urban, the middle bar \begin{tikzpicture}[xscale=\xScaleProportionOfUrbanSuburbanRuralStatinsOpenedInDifferentScenariosForTheTestInstanceSHundredWFifty, yscale=1]\draw[red, thick, fill=red!25] (0.0, 0.0) rectangle (0.2, 0.24);\end{tikzpicture} for suburban, and the right bar \begin{tikzpicture}[xscale=\xScaleProportionOfUrbanSuburbanRuralStatinsOpenedInDifferentScenariosForTheTestInstanceSHundredWFifty, yscale=1]\draw[blue, thick] (0.0, 0.0) rectangle (0.2, 0.24);\end{tikzpicture} for rural areas.}
				\label{figure:totalNumberAndProportionOfUrbanSuburbanRuralStatinsOpenedInDifferentScenariosForTheTestInstanceSHundredWFifty}
			\end{figure}

\paragraph{S.2 scenario sequences.}
The only changing cost parameter is in rural construction costs. Consequently, it is possible to locate the very same CSs of the previous scenario in the current scenario, as the total number of locations in the cost increasing category exceeds 4 times the number of rural CSs that were located in the previous scenarios. If these number is below the threshold of 12 rural CSs ($4\ \cdot $ \# opened rural CSs $\leq 48$), the cost increase can be covered from the increase in budget. 

\paragraph{Scenario sequence S.3.}
Comparing the results of scenario 1, 6, 10, and 13, there are two changing cost categories (rural and sub-urban). There is a slight change in structure comparing scenario 1 and 6 and scenario 10 and 13: the general trend is that the number or rural CSs is decreasing when changing from scenario 1 to 6, and scenario 10 to 13. But there is a stronger structural change when changing from scenario 6 to 10. The reason therefore is explained in more detail in Section~\ref {subsubsection:naLocationDependentCost}, where the results of the baseline instance s60w30 are discussed. 
Moreover, the proportion of urban CSs is increasing when changing from scenario 1 to 13. 

\medskip

Finally, in Table~\ref{table:LC:SolveTime} the solve times for all test instances and scenarios are summarized.

\begin{table}[htb]
\caption{Solve time in seconds for all scenarios and test instances.} 
	\begin{center}
			\begin{tabular}[c]{c *{5}{| r} }
				\label{table:LC:SolveTime}
		\diagbox{Inst.}{Scen.}
               &	1	&	2	&	3	&	4	&	5\\
                \hline
s40w20 & 3.18 & 2.22 & 2.05  & 2.87  & 	3.81\\
s60w30 & 31.00 & 53.87 & 39.27 & 44.54 & 45.06\\
s80w40 & 182.65 & 236.59 & 199.07 &  204.60 & 208.65\\
s100w50 & 459.34 & 489.97 & 587.81 & 586.97 & 664.42\\
				\hline
				\multicolumn{5}{c}{}\\[-0.5em]
               	&	6	&	7	&	8	&	9	&	10\\
                \hline
s40w20 & 3.48 & 4.43 &  2.68 & 4.43 & 3.72\\
s60w30 & 22.21 & 38.90 & 46.18 & 29.60 & 25.08\\
s80w40 & 166.93 & 154.05 &  158.18 & 160.12 & 167.52\\
s100w50 & 497.42 & 478.30 & 719.57 & 589.73 & 694.66\\
				\hline
				\multicolumn{5}{c}{}\\[-0.5em]
               &	11	&	12	&	13	&	14	&	15 \\
                \hline
s40w20 & 2.56 & 1.80 & 2.08 & 2.12 & 7.74 \\
s60w30 & 30.19 & 47.90 & 18.57 & 21.05 & 27.87 \\
s80w40 & 144.13 & 190.92 & 343.15 & 115.59 & 153.81 \\
s100w50 & 724.88 & 565.75 & 655.57 & 675.72 & 614.66 \\
				\hline
			\end{tabular}
	\end{center}
\end{table}

	\subsection{Determination of the station size (\CFRLP)}
	\label{Appendix:CFRLP:results}
	
	\paragraph{Analysing s40w20.}
Hereafter, the results of the \CFRLP\ for test instance s40w20, where capacity is chosen to be \num{2801,01} (see Table~\ref{table:CapacityPerPole}) are described in detail. Pre-testing, whose result is outlined in Table~\ref{table:Results_CMCFRLP_pre_s40w20}, states that with the previously defined capacity of a charging pole it is not possible to cover more than \num{521918.7265} of \num{e06} (52.19\%) EVs, denoting TFV in this test instance. The second pre-test shows, that this number of EVs can be covered by locating 72 charging poles, representing the EVCP. The placement of further charging poles does not result in an increase of CFV.

\begin{table}[htb]
\caption{Pre-test (\CMCFRLP) of instance s40w20.} 
    \begin{center}
\begin{tabular}[c]{c *{1}{| r} }
                \label{table:Results_CMCFRLP_pre_s40w20}
                $C$                            & 0.521919  \\ 
                 \hline            
                 $\sum_{k \in K} n_k$ & 72\\    
                 $\sum_{f \in f} x_k$     & 26\\    
                 PCF                            & 64\\    
                 solve time                 & 53.55\\ 
                    \hline
			        \multicolumn{2}{c}{}\\[-0.4cm]
                    \multicolumn{2}{c}{\footnotesize{PCF: number of partially covered flows ($z_f > 0$).}}        
            \end{tabular}
    \end{center}
\end{table}

\begin{table}[htb]
\caption{Results of \CFRLP\ for instance s40w20.} 
    \begin{center}
            %
			\end{comment:figures}
		\caption{Test instance s40w20: \CFRLP\ -- $S = 36$  (50\%).}
		\label{figure:TestInstanceSFourtyWTwentyCFRLPSThirtySix}
	\end{figure}
	\begin{figure}[htb!]
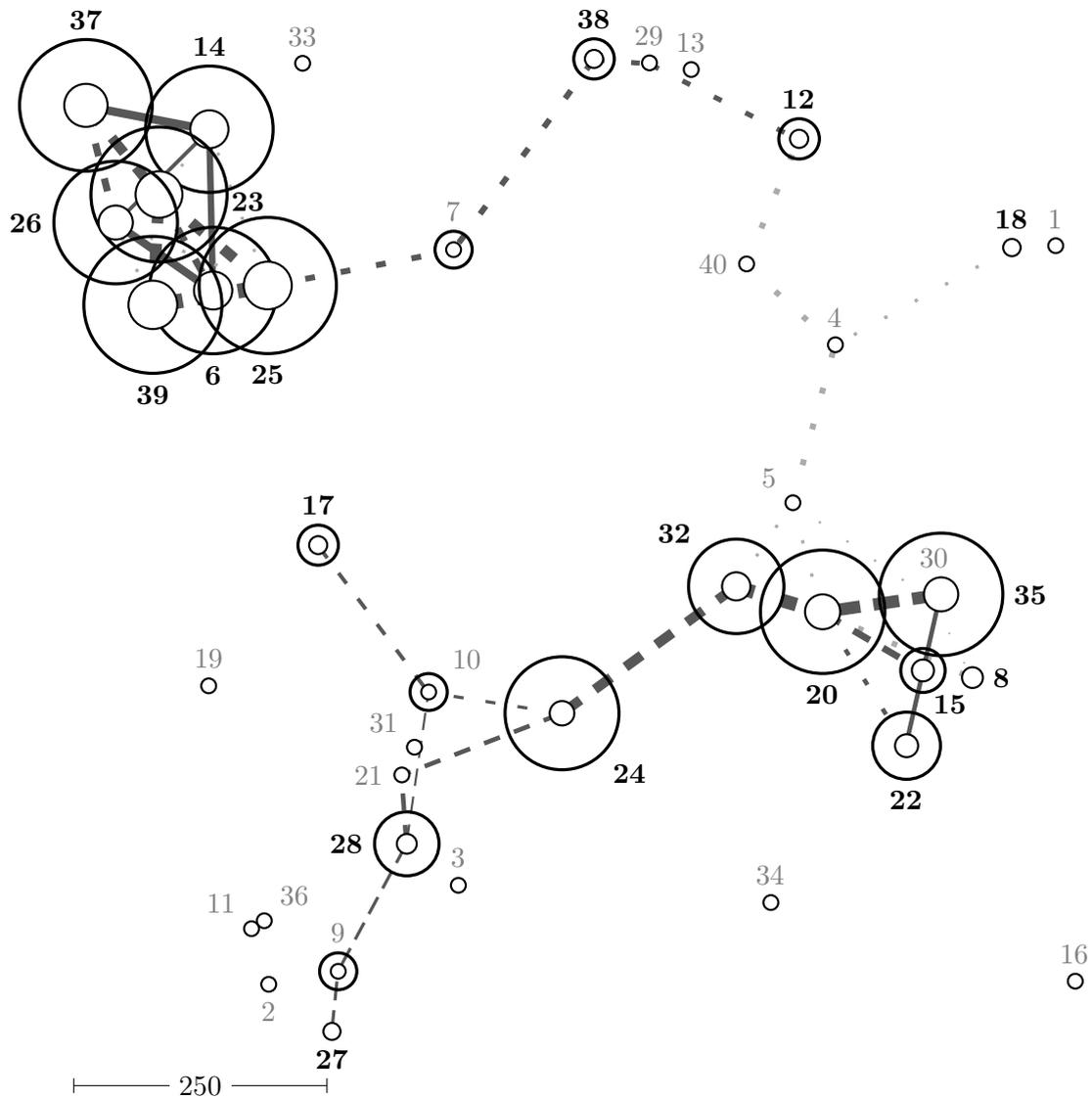

		\centering
			\begin{comment:figures}
				%
			\end{comment:figures}
		\caption{Test instance s40w20: \CFRLP\ -- $S = 54$  (75\%).}
		\label{figure:TestInstanceSFourtyWTwentyCFRLPSFiftyFour}
	\end{figure}

The testing process is repeated for different numbers of charging poles to locate, representing 25\%, 50\%, 75\% and 100\% of EVCP and the results are shown in Table~\ref{table:Results_CFRLP_s40w20}. 
Locating 25\% ($S=18$) of the EVCP covers 42.42\% of the maximum CFV that can be guaranteed given the limited capacity of charging stations. 
In case of locating 75\% ($S=54$) of the EVCP, it is still possible to cover 95.95\% of the CFV, which would be possible if allocating the maximum EVCP.  

Looking at the average number of charging poles per station, it is not continually increasing with a growing number of charging poles to locate, like in the baseline case s60w30, which is described in Section~\ref{subsubsection:naCFRLP}. While in case of allocating 25\% of EVCP, there are on average 3 charging poles per station, locating 100\% of EVCP results on average in 2.77 charging poles per location.
This can be explained by this test instance's characteristic cluster of strongly frequented origin and destination nodes. Within the cluster, maximum-sized charging stations are built, in case there are sufficient charging poles available to locate. If there are not enough charging poles left to install a maximum-sized station, smaller charging station in close proximity to the large ones are build in order to reduce energy demand at the strongly frequented stations. 

In particular, if there are few charging poles to locate, they are primarily located within the cluster in order to build maximum-sized charging stations. As the number of charging poles increases, the number of maximum sized charging stations within the cluster increases. If no more maximum sized charging stations are necessary to cover additional EVs, charging poles are located along less frequented nodes, where in most cases the capacity of a single charging pole is sufficient to cover the passing EVs.
This explains why the average number of charging poles per location is decreasing when comparing the allocation of 50\% (average of 3 charging poles) to 75\% (average of 2.7 charging poles) of EVCP. When 36 charging poles are allocated (see Figure~\ref{figure:TestInstanceSFourtyWTwentyCFRLPSThirtySix}), there are seven maximum-sized charging locations, two stations with three charging poles and two other ones with a single pole. If 54 charging poles are built (see Figure~\ref{figure:TestInstanceSFourtyWTwentyCFRLPSFiftyFour}), the additional charging poles are used to guarantee refuelling for long-distance round-trips in more remote areas. Consequently, the number of charging locations increases considerable, as an increasing number of small charging stations is necessary to cover flows in remote areas. 


\paragraph{Analysing s80w40.}
When the test instance size increases, due to an increasing number of potential facility locations, Gurobi was no longer able to find the optimal solution within a reasonable time. Columns marked with an asterisk point out that the problem was solved with a given time limit. In addition, when the number of charging poles to locate decreases, the reported gap between the best objective and best bound found by Gurobi increases. Whereby, considering the gap, it is important to mention that the best objective value found does not change for a long time when approaching the end of the time limit. However, the best bound found increases steadily. This can be seen in Appendix~\ref{GurobiLogFile}.

The results of the pre-test is shown in Table~\ref{table:Results_CMCFRLP_pre_s80w40}. 
The first pre-testing phase (solving \CFRLP) is interrupted after four hours, the second pre-test (solving \CMCFRLP) after one hour, because the Gurobi MIP logfile shows that the incumbent value did not change in the last \num{13000} and \num{2100}~seconds, respectively when approaching the end ot the time limit. 
The solution for the \CFRLP\ applied to test instance s80w40 is summarized in Table~\ref{table:Results_CFRLP_s80w40}.
Increasing the number of charging poles to locate leads to an increasing average size of charging stations. A detailed analysis concerning the covered flows indicates that the number of fully covered flows ($z_f = 1$) is steadily increasing, while the number of flows that are covered to less than 50\% of their flow volume ($z_f < 0.5$) is decreasing in case the number of charging poles to locate increases. 

\begin{table}[htb]
\caption{Pre-test for \CMCFRLP\ of instance s80w40.} 
    \begin{center}
            \begin{tabular}[c]{c *{1}{| r} }
                \label{table:Results_CMCFRLP_pre_s80w40}
                $C$                            & 0.582610$^*$ \footnotemark \\ 
                 \hline            
                 $\sum_{k \in K} n_k$ & 174 \footnotemark \\    
                 $\sum_{f \in f} x_k$     & 62 \\    
                 PCF                             & 179 \\    
                 solve time                 & \num{4105.60} \\ 
                    \hline
			        \multicolumn{2}{c}{}\\[-0.4cm]
                    \multicolumn{2}{c}{\footnotesize{PCF: number of partially covered flows ($z_f > 0$).}}     \\
                    \multicolumn{2}{c}{\footnotesize{$^*$ indicates no optimal results, interrupted solving process.}}    
            \end{tabular}
    \end{center}
\end{table}
        \footnotetext[8]{Best objective \num{58261.01}, best bound \num{792966.49}, gap 36.1058\% (after \num{14400} sec).}
        \footnotetext{Best objective 174, best bound 20, gap 88.5057\% (after \num{3600} sec).}

\begin{table}[htb]
\caption{Results of \CFRLP\ for instance s80w40.} 
    \begin{center}
            \begin{tabular}[c]{c *{4}{| r} }
                \label{table:Results_CFRLP_s80w40}
                $S$& 174$^*$ (100\%) & 130$^*$ (75\%)    &    87$^*$ (50\%)    &    43$^*$ (25\%) \\
                 \hline    
                 CFV                            & \num{509797.28} \footnotemark &  \num{496324.80} \footnotemark & \num{437704.45} \footnotemark & \num{249329.60} \footnotemark \\                
                 $\sum_{k \in K} x_k$ & 52                 & 41             & 36         & 25 \\  
                 $\overline{n}$			& 3.35 			& 3.39			& 2.42		& 1.72 \\
                 $\sum_{f \in f} y_f$     & 255             & 282        & 221         & 167\\    
                 PCF                            & 90             & 131        & 62            & 21\\    
                 solve time                 & \num{3609.12} & \num{3608.58} & \num{3608.73} & \num{3608.55}\\ 
                    \hline
			        \multicolumn{5}{c}{}\\[-0.4cm]
                     \multicolumn{5}{c}{\footnotesize{$\overline{n}$: average station size	.}} \\   
                    \multicolumn{5}{c}{\footnotesize{$^*$ indicates no optimal results, interrupted solving process.}}
            \end{tabular}
    \end{center}
\end{table}    

\footnotetext[10]{Best bound \num{836727,71}, gap 64.13\%.}
\footnotetext[11]{Best bound \num{835855,04}, gap 68.41\%.}
\footnotetext[12]{best bound \num{848569,23}, gap 93.87\%}
\footnotetext{Best bound \num{705759,43}, gap 183.06\%.}

Due to the fact, that the \CFRLP\ could not be solved to optimality within a reasonable time, a time limit of \num{3600}~seconds is set. That is the reason for different results when comparing the columns of Tables \ref{table:Results_CMCFRLP_pre_s80w40} and \ref{table:Results_CFRLP_s80w40} in case of allocating 174 charging poles.

\paragraph{Analysing s100w50.}
Table~\ref{table:Results_CMCFRLP_pre_s100w50} depicts the results obtained from pre-testing and the solution of the \CFRLP\ is summarized in Table~\ref{table:Results_CFRLP_s100w50}. Because the solving processed ended after reaching a time limit, differences in the CFV again occurs in case of testing the \CFRLP\ with the results of the pre-testing model \CMCFRLP, representing the EVCP to locate (S=219) .

\begin{table}[htb]
\caption{Pre-test (\CMCFRLP) of instance s100w50.} 
    \begin{center}
            \begin{tabular}[c]{c *{1}{| r} }
                \label{table:Results_CMCFRLP_pre_s100w50}
                $C$                            & 0.604579$^*$ \footnotemark \\ 
                 \hline            
                 $\sum_{k \in K} n_k$ & 219 \footnotemark\\    
                 $\sum_{f \in f} x_k$     & 82\\    
                 PCF                            & 273\\    
                 solve time                 & 4019.16 \\ 
                    \hline
			        \multicolumn{2}{c}{}\\[-0.4cm]
                    \multicolumn{2}{c}{\footnotesize{PCV: number of partially covered flows ($z_f > 0$).}} \\
                    \multicolumn{2}{c}{\footnotesize{$^*$ indicates no optimal results, interrupted solving process.}}    
            \end{tabular}
    \end{center}
\end{table}
\footnotetext[14]{Best objective \num{604579,47}, best bound \num{834075,08}, gap 37.96\% (after \num{14400} sec).}
\footnotetext{Best objective 219, best bound 20, gap 90.87\% (after \num{3600} sec).}

\begin{table}[htb]
\caption{Results of \CFRLP\ for instance s100w50.} 
    \begin{center}
            \begin{tabular}[c]{c *{4}{| r} }
                \label{table:Results_CFRLP_s100w50}
                $S$ & 219$^*$ (100\%) & 164$^*$ (75\%)    &    109$^*$ (50\%)    &    54$^*$ (25\%) \\
                 \hline    
                 CFV                            & \num{524261.32} \footnotemark    & \num{513501.34} \footnotemark & \num{433967.79} \footnotemark & \num{247072.03} \footnotemark \\                
                 $\sum_{k \in K} x_k$ & 63             & 59                & 50            & 36 \\    
                 $\overline{n}$			&	3.48		& 2.78				& 2.18		& 1.50 \\
                 $\sum_{f \in f} y_f$     & 1132         & 1132            & 391        & 301 \\    
                 PCF                            & 198         & 180            & 105        & 48 \\    
                 solve time                 & 3619.27    & 3621.02     & 3620.95    & 3621.99 \\ 
                    \hline
			        \multicolumn{5}{c}{}\\[-0.4cm]
                    \multicolumn{5}{c}{\footnotesize{$\overline{n}$: average station size	.}} \\   
                    \multicolumn{5}{c}{\footnotesize{$^*$ indicates no optimal results, interrupted solving process.}}        
            \end{tabular}
    \end{center}
\end{table}

A detailed analysis of the results indicates that allocating 400 charging poles (installing four charging poles at every potential facility location) results in an average utilisation of 57.03\% per charging location. The second pre-testing phase indicates that the same coverage level can be reached when allocating less charging poles (S=219, representing the EVCP), which consequently results in a higher utilisation per charging location.

\footnotetext[16]{Best bound \num{869785.76}, gap 65.91\%.}
\footnotetext[17]{Best bound \num{869785.76}, gap 69.38\%.}
\footnotetext[18]{Best bound \num{858707.29}, gap 97.87\%.}
\footnotetext{Best bound \num{793058.56}, gap 220.98\%.}    

It can be summarised that a decreasing number of charging poles to locate results in an increasing average utilisation per charging location. Locating 54 charging poles (25\% of EVCP), placed in 35 potential facility locations, leads to an average utilisation of 91.24\%. Having enough budget to install 164 charging poles (75\% of EVCP), most EVs can be covered if these locations are installed in 59 possible facility locations. The average utilisation per charging location based on this charging pole allocation is 82.64\%.

\paragraph{Analysing Florida.}

The capacity per charging pole is \num{23740,35286}. Whereby, consistent with the other test instances, a scaling parameter of 0.001 is used. The TFV of this test instance is \num{1e+12}.

\begin{table}[htb]
\caption{Pre-test (\CMCFRLP) of instance Florida.} 
    \begin{center}
            \begin{tabular}[c]{c *{1}{| r} }
                \label{table:Results_CMCFRLP_pre_Florida}
                $C$                            & 0.000828\\ 
                 \hline            
                 $\sum_{k \in K} n_k$ & 870 \footnotemark \\    
                 $\sum_{f \in f} x_k$     & 293\\    
                 PCF                            & 130\\    
                 solve time                 & \num{3854.78} \\ 
                    \hline
                    \multicolumn{2}{c}{\footnotesize{PCV: number of partially covered flows ($z_f > 0$).}} \\
                    \multicolumn{2}{c}{\footnotesize{$^*$ indicates no optimal results, interrupted solving process.}}    
            \end{tabular}
    \end{center}
\end{table}
\footnotetext{Best bound \num{1}, gap 99.8851\%}

\begin{table}[htb]
\caption{Results of \CFRLP\ for instance Florida.} 
    \begin{center}
            \begin{tabular}[c]{c *{4}{| r} }
                \label{table:Results_CFRLP_Florida}
                $S$ & 870$^*$ (100\%) & 652$^*$ (75\%)    &    435$^*$ (50\%)    &    217$^*$ (25\%) \\
                 \hline    
                 CFV                            & \num{6,60E+08} \footnotemark & \num{5,39e+08} \footnotemark & \num{4,06e+08} \footnotemark & \num{2.69e+08} \footnotemark  \\                
                 $\sum_{k \in K} x_k$ & 218 & 163 & 109 & 55\\    
                 $\overline{n}$			&	4 & 4 & 3.99 & 3.95 \\
                 $\sum_{f \in f} y_f$     & \num{2530} & \num{2187} & \num{1066} & 469  \\    
                 PCF                           & 107 &87  & 60 & 37 \\    
                 solve time             & \num{3847.80} & \num{3867.11} & \num{3871.50} &   \num{3878.55} \\ 
                    \hline
                    \multicolumn{5}{c}{\footnotesize{$\overline{n}$: average station size	.}} \\   
                    \multicolumn{5}{c}{\footnotesize{$^*$ indicates no optimal results, interrupted solving process.}}        
            \end{tabular}
    \end{center}
\end{table}    
\footnotetext[21]{Best bound \num{9.5773e+11}, gap \num{144939.62}\%.}
\footnotetext[22]{Best bound \num{9.5773e+11}, gap \num{177709.08}\%.}
\footnotetext[23]{Best bound \num{9.5773e+11}, gap \num{235905.02}\%.}
\footnotetext{Best bound \num{9.5773e+11}, gap \num{356169.92}\%.}

\section{Gurobi MIP log file}
\label{GurobiLogFile}
To obtain information about the branch-and-bound tree  and therefore the process of the optimization, the Gurobi option \verb|outlev=1| is used. 
Briefly describing the columns, the ``Nodes'' columns give the current number of explored and unexplored nodes in the branch-and-bound tree. 
The ``Current Node'' columns list information about the current node, where ``Obj'' indicates the optimal value of the LP relaxation at the current node. If an integer solution is found at a node, it is indicated with a ``H'' in the first column. The ``BestBd'' column gives the value of the best integer solution found so far. The ``Gap'' column represents the relative gap between the best integer solution found and the current upper bound if an integer solution has been found.

\begin{center}
\small
\begin{verbatim}
CFRLP for s80w40 with S = 174 
Gurobi 8.1.0: outlev=1
timelim=3600

Root relaxation: objective 8.872952e+05, 201679 iterations, 1182.99 seconds
Total elapsed time = 1222.65s
Total elapsed time = 1225.98s

    Nodes    |    Current Node    |     Objective Bounds      |     Work
 Expl Unexpl |  Obj  Depth IntInf | Incumbent    BestBd   Gap | It/Node Time

     0     0 887295.163    0 18548 486490.066 887295.163  82.4%     - 1235s
     0     0 874583.919    0 18302 486490.066 874583.919  79.8%     - 1432s
     0     0 871990.641    0 18552 486490.066 871990.641  79.2%     - 1485s
     0     0 871223.741    0 18603 486490.066 871223.741  79.1%     - 1507s
     0     0 870967.512    0 18585 486490.066 870967.512  79.0%     - 1522s
     0     0 870843.577    0 18477 486490.066 870843.577  79.0%     - 1531s
     0     0 870760.377    0 18533 486490.066 870760.377  79.0%     - 1540s
     0     0 870683.381    0 18605 486490.066 870683.381  79.0%     - 1548s
     0     0 870674.098    0 18604 486490.066 870674.098  79.0%     - 1552s
     0     0 870664.907    0 18598 486490.066 870664.907  79.0%     - 1553s
     0     0 870555.192    0 18581 486490.066 870555.192  78.9%     - 1556s
     0     0 870545.617    0 18540 486490.066 870545.617  78.9%     - 1559s
     0     0 870544.164    0 18535 486490.066 870544.164  78.9%     - 1561s
     0     0 863255.026    0 18793 486490.066 863255.026  77.4%     - 1756s
H    0     0                    509797.27780 863255.026  69.3%     - 1756s
     0     0 858161.471    0 18990 509797.278 858161.471  68.3%     - 1906s
     0     0 855617.300    0 19154 509797.278 855617.300  67.8%     - 2032s
     0     0 854943.086    0 19115 509797.278 854943.086  67.7%     - 2083s
     0     0 854853.820    0 19189 509797.278 854853.820  67.7%     - 2103s
     0     0 854788.538    0 19189 509797.278 854788.538  67.7%     - 2121s
     0     0 854565.594    0 19135 509797.278 854565.594  67.6%     - 2147s
     0     0 854505.514    0 18939 509797.278 854505.514  67.6%     - 2162s
     0     0 854490.623    0 18946 509797.278 854490.623  67.6%     - 2169s
     0     0 854482.530    0 18912 509797.278 854482.530  67.6%     - 2172s
     0     0 854455.832    0 18870 509797.278 854455.832  67.6%     - 2179s
     0     0 854453.779    0 18901 509797.278 854453.779  67.6%     - 2182s
     0     0 849404.235    0 18899 509797.278 849404.235  66.6%     - 2394s
     0     0 848137.493    0 19076 509797.278 848137.493  66.4%     - 2484s
     0     0 847840.466    0 19118 509797.278 847840.466  66.3%     - 2525s
     0     0 847771.492    0 18968 509797.278 847771.492  66.3%     - 2551s
     0     0 847678.926    0 19053 509797.278 847678.926  66.3%     - 2581s
     0     0 847610.530    0 19128 509797.278 847610.530  66.3%     - 2599s
     0     0 847503.941    0 19055 509797.278 847503.941  66.2%     - 2623s
     0     0 847487.687    0 19008 509797.278 847487.687  66.2%     - 2632s
     0     0 847418.223    0 18988 509797.278 847418.223  66.2%     - 2647s
     0     0 847405.285    0 19116 509797.278 847405.285  66.2%     - 2657s
     0     0 847352.632    0 19161 509797.278 847352.632  66.2%     - 2677s
     0     0 847340.687    0 19051 509797.278 847340.687  66.2%     - 2686s
     0     0 847329.147    0 18980 509797.278 847329.147  66.2%     - 2692s
     0     0 847313.166    0 19039 509797.278 847313.166  66.2%     - 2708s
     0     0 847307.230    0 19043 509797.278 847307.230  66.2%     - 2715s
     0     0 843477.349    0 19019 509797.278 843477.349  65.5%     - 2901s
     0     0 842437.503    0 19131 509797.278 842437.503  65.2%     - 3006s
     0     0 842208.248    0 19091 509797.278 842208.248  65.2%     - 3043s
     0     0 842090.928    0 18988 509797.278 842090.928  65.2%     - 3070s
     0     0 841930.201    0 19002 509797.278 841930.201  65.1%     - 3098s
     0     0 841869.615    0 19056 509797.278 841869.615  65.1%     - 3119s
     0     0 841850.829    0 18932 509797.278 841850.829  65.1%     - 3128s
     0     0 841827.107    0 18810 509797.278 841827.107  65.1%     - 3138s
     0     0 841735.417    0 18780 509797.278 841735.417  65.1%     - 3151s
     0     0 841721.997    0 18830 509797.278 841721.997  65.1%     - 3159s
     0     0 841717.815    0 18855 509797.278 841717.815  65.1%     - 3165s
     0     0 837929.264    0 19182 509797.278 837929.264  64.4%     - 3390s
     0     0 837074.367    0 19344 509797.278 837074.367  64.2%     - 3502s
     0     0 836851.104    0 19373 509797.278 836851.104  64.2%     - 3547s
     0     0 836770.348    0 19368 509797.278 836770.348  64.1%     - 3577s
     0     0 836727.707    0 19426 509797.278 836727.707  64.1%     - 3596s
     0     0          -    0      509797.278 836727.707  64.1%     - 3600s
\end{verbatim}
\normalsize
\end{center}

\end{document}